\documentclass[11pt,a4paper]{amsart}

\usepackage[T1]{fontenc}
 \usepackage[english]{babel}
 \usepackage[latin1]{inputenc}
\usepackage{fullpage}
\usepackage[dvips]{graphics}
\usepackage{color}
\usepackage{url}
\usepackage{stmaryrd}
\usepackage{rotating}
\usepackage{longtable}

\usepackage{mathabx}
\usepackage{amsmath}
\usepackage{amsthm}
\usepackage{amssymb}
\usepackage{bbm}
\usepackage{enumerate}
\usepackage[all,poly,necula]{xy}

\DeclareSymbolFont{bbdold}{U}{bbold}{m}{n}
\DeclareSymbolFontAlphabet{\mathbbd}{bbdold}
\newcommand{\un}{\mathbbd 1}

\theoremstyle{plain}
\newtheorem{thm}{Theorem}[section]

\newtheorem{prop}[thm]{Proposition}

\newtheorem{lem}[thm]{Lemma}

\newtheorem{cor}[thm]{Corollary}

\newtheorem{conj}[thm]{Conjecture}

\theoremstyle{definition}
\newtheorem{df}[thm]{Definition}

\newtheorem{nt}[thm]{Notation}

\newtheorem{ex}[thm]{Example}

\newtheorem{exs}[thm]{Examples}

\theoremstyle{remark}
\newtheorem*{demo}{Proof}
\newtheorem{rem}[thm]{Remark}

\newtheorem{rems}[thm]{Remarks}

\newcounter{thmnct}

\newenvironment{thmn}[1]{
 \addtocounter{thmnct}{1}
 \theoremstyle{plain}
 \newtheorem*{thmntmp\roman{thmnct}}{Theorem #1}
 \begin{thmntmp\roman{thmnct}}
}{\end{thmntmp\roman{thmnct}}}

\DeclareMathOperator{\GL}{GL}

\DeclareMathOperator{\id}{Id}

\DeclareMathOperator{\coker}{coker}

\DeclareMathOperator{\End}{End}
\DeclareMathOperator{\Aut}{Aut}

\DeclareMathOperator{\Gr}{Gr}

\DeclareMathOperator{\Hom}{Hom}
\DeclareMathOperator{\uHom}{\mathbf{Hom}}
\DeclareMathOperator{\Ext}{Ext}
\DeclareMathOperator{\uExt}{\mathbf{Ext}}

\DeclareMathOperator{\add}{add}

\DeclareMathOperator{\md}{mod}

\DeclareMathOperator{\gldim}{gl. dim}
\DeclareMathOperator{\prdim}{proj. dim}
\DeclareMathOperator{\op}{op}
\DeclareMathOperator{\rad}{rad}
\DeclareMathOperator{\ind}{ind}
\DeclareMathOperator{\irr}{irr}

\DeclareMathOperator{\Lie}{Lie}
\DeclareMathOperator{\Add}{\mathfrak{Add}}
\DeclareMathOperator{\adds}{\mathfrak{add}}
\DeclareMathOperator{\udim}{\mathbf{dim}}
\DeclareMathOperator{\Sub}{Sub}

\DeclareMathOperator{\gr}{gr}

\newcommand{\vect}{\md}

\renewcommand{\leq}{\leqslant}
\renewcommand{\geq}{\geqslant}
\newcommand{\Ar}{\mathcal{A}}
\newcommand{\Cr}{\mathcal{C}}
\newcommand{\Dr}{\mathcal{D}}
\newcommand{\Er}{\mathcal{E}}

\newcommand{\Mr}{\mathcal{M}}
\newcommand{\Tr}{\mathcal{T}}
\newcommand{\Erg}{\underline{\mathcal{E}}}

\newcommand{\N}{\mathbb{N}}
\newcommand{\Z}{\mathbb{Z}}
\newcommand{\Q}{\mathbb{Q}}

\newcommand{\C}{\mathbb{C}}

\newcommand{\G}{\mathbf{{\Gamma}}}
\newcommand{\g}{\mathbf{g}}
\newcommand{\h}{\mathbf{h}}
\renewcommand{\H}{\mathbf{H}}

\renewcommand{\epsilon}{\varepsilon}

\newcommand{\T}{\mathcal{T}}

\newcommand{\inter}{\cap}

\newcommand{\restr}[1]{\!|_{#1}}

\newcommand{\ngo}{\mathfrak{n}}
\newcommand{\ggo}{\mathfrak{g}}
\newcommand{\bgo}{\mathfrak{b}}
\newcommand{\ag}{\mathbf{a}}
\newcommand{\ig}{\mathbf{i}}
\newcommand{\jg}{\mathbf{j}}

\newcommand{\dg}{\mathbf{d}}
\newcommand{\tg}{\mathbf{t}}

\newcommand{\tens}{\otimes}
\newcommand{\cqfd}{\qed}
\renewcommand{\phi}{\varphi}

\renewcommand{\bar}[1]{\overline{#1}}

\newcommand{\mat}[4]{\left(\!\!
    \begin{array}{cc}
      #1 & #2 \\
      #3 & #4
    \end{array}
    \!\!\right)}

\newcommand{\inj}{\hookrightarrow}
\newcommand{\surj}{\twoheadrightarrow}
\newcommand{\tr}[1]{\,{\vphantom{#1}}^{\text{t}}\!{#1}}
\def\citeb#1#2{\cite[#1]{#2}}
\newcommand{\obullet}{\odot \hspace{-.63em} \bullet}
\renewcommand{\tilde}[1]{\widetilde{#1}}

\definecolor{orange}{rgb}{1,0.90,0.5}
\let\oldmarginpar\marginpar
\renewcommand\marginpar[1]{\-\oldmarginpar[\raggedleft\footnotesize \fcolorbox{blue}{orange}{\parbox{\marginparwidth}{\color{blue}{#1}}}]%
{\raggedright\footnotesize \fcolorbox{blue}{orange}{\parbox{\marginparwidth}{\color{blue}{#1}}}}}

\newcommand{\forloop}[5][1] { \setcounter{#2}{#3} \ifthenelse{#4} { #5 \addtocounter{#2}{#1} \forloop[#1]{#2}{\value{#2}}{#4}{#5} }{ } } 
\newcounter{i}
\newcounter{j}
\newcommand{\dets}[4]{
    \displaystyle \left|\;
    \begin{matrix}
      \forloop{i}{1}{\value{i} < 7}{
       \forloop{j}{1}{\value{j} < 7}{
        \ifthenelse{\value{i} < #1 \or \value{i} > #2 \or \value{j} < #3 \or \value{j} > #4}{\ifthenelse{\value{i}<\value{j}}{\cdot}{\ifthenelse{\value{i} = \value{j}}{1}{}}}{\bullet} \ifthenelse{\value{j}<6}{&}{}
       } \ifthenelse{\value{i}<6}{\\}{}
      }
    \end{matrix}
    \;\right|^{\vphantom{M^{M^M}}}_{\vphantom{M_{M_M}}}}

\begin{document}
\setlength{\parindent}{5mm}
\title{Categorification of skew-symmetrizable cluster algebras}
\author{Laurent Demonet}
\address{LMNO, Universit\'e de Caen, Esplanade de la Paix, 14000 Caen}
\email{Laurent.Demonet@normalesup.org}

\date{}

\begin{abstract}
 We propose a new framework for categorifying skew-symmetrizable cluster algebras. Starting from an exact stably $2$-Calabi-Yau category $\Cr$ endowed with the action of a finite group ${\Gamma}$, we construct a ${\Gamma}$-equivariant mutation on the set of maximal rigid ${\Gamma}$-invariant objects of $\Cr$. Using an appropriate cluster character, we can then attach to these data an explicit skew-symmetrizable cluster algebra. As an application we prove the linear independence of the cluster monomials in this setting. Finally, we illustrate our construction with examples associated with partial flag varieties and unipotent subgroups of Kac-Moody groups, generalizing to the non simply-laced case several results of Gei\ss-Leclerc-Schröer.
\end{abstract}

\maketitle
\tableofcontents

\section{Introduction}

\subsection{Cluster algebras}
In 2001, Fomin and Zelevinsky introduced a new class
of algebras called  cluster algebras \cite{FoZe02}, \cite{FoZe03} motivated by canonical bases and total positivity \cite{Lu91}, \cite{Lu97}. By construction, a cluster algebra is a commutative ring endowed with distinguished generators (\emph{cluster variables}) grouped in subsets of the same cardinality (\emph{clusters}). The clusters are not disjoint. On the contrary, each cluster $A$ has neighbours obtained by replacing each of its variables $x_i$ by a new variable $x'_i$. The new cluster $\mu_i(A) = A \setminus \{x_i\} \cup \{x'_i\}$ is called the mutation of $A$ in the direction $x_i$. Moreover, the mutations are always of the form
\[
x_i x'_i = M_i + M_i',
\]
where $M_i$ and $M_i'$ are monomials in the variables of $A$ other than $x_i$. The axioms imply strong compatibility relations between the monomials of the exchange relations. In particular, a cluster algebra is fully determined by a \emph{seed}, that is, a single cluster and its exchange relations with all its neighbours. In practice, one usually defines a cluster algebra by giving an initial seed. By iterating the exchange relations, one can express every cluster variable in terms of the variables of the initial seed.

Berenstein, Fomin and Zelevinsky have shown that the coordinate rings of many algebraic varieties attached to complex semi-simple Lie groups were endowed with the structure of a cluster algebra \cite{BeFoZe05}. Other examples have been given by 
Gei\ss, Leclerc and Schröer \cite{GeLeSc08}, \cite{GeLeSc}.
 
Since their emergence, cluster algebras have aroused a lot of interest, coming in particular from their links with many other subjects: combinatorics (see for instance \cite{ChFoZe02}, \cite{FoShTh08}), Poisson geometry \cite{GeShVa03}, \cite{GeShVa05}, integrable systems \cite{FoZe03-1}, Teichm\"uller spaces \cite{FoGo06}, and, last but not least, representations of finite-dimensional algebras.

Unfortunately, because of the inductive description of cluster algebras, many properties of the cluster variables which might seem elementary are in fact very hard to prove. For instance:

\begin{conj}[Fomin-Zelevinsky] \label{cFZ}
 Cluster monomials (that is products of cluster variables of a single cluster) are linearly independent.
\end{conj}

In seminal articles, Marsh, Reineke, Zelevinsky \cite{MaReZe03},
Buan, Marsh, Reineke, Reiten, Todorov \cite{BuMaReReTo06} and Caldero, Chapoton \cite{CaCh06}
have shown that the important class of acyclic cluster algebras could be modelled with categories constructed from representations of quivers. This gives in particular a global (i.e. non inductive) understanding of these algebras, and gives new tools for studying them. For example, this allowed Fu and Keller \cite{FuKe} to prove the previous conjecture for a family of cluster algebras containing acyclic cluster algebras.

At the same time, Gei\ss, Leclerc and Schröer have studied cluster algebras associated with Lie groups of type $A$, $D$, $E$,
and have modelled them by categories of modules over Gelfand-Ponomarev preprojective algebras of the same type. 
They have shown that cluster monomials form a subset of the dual semi-canonical basis \cite{GeLeSc06} introduced by
Lusztig \cite{Lu00}, proving the above conjecture in this other context.

More recently, Derksen, Weyman and Zelevinsky \cite{DeWeZe08}, \cite{DeWeZe} have obtained a far-reaching generalization of \cite{MaReZe03}, which also contains all the above examples. They have shown that one can control $F$-polynomials and $\mathbf{g}$-vectors of every cluster algebra whose initial seed is encoded by a skew-symmetric matrix, using representations of quivers with potentials. This enabled them to prove the linear independence conjecture, as well as many other conjectures on $F$-polynomials and $\mathbf{g}$-vectors formulated in \cite{FoZe07}.

But the theory of Fomin and Zelevinsky includes more general seeds given by skew-symme\-tri\-za\-ble matrices
(i.e. products of a skew-symmetric matrix by a diagonal matrix). For example, cluster algebras associated to Lie groups of type $B$, $C$, $F$, $G$ are only skew-symme\-tri\-za\-ble. The aim of this article is to extend the results of Gei\ss, Leclerc, Schröer and Fu, Keller to the skew-symme\-tri\-za\-ble case.

\subsection{Actions of groups on categories} 
    
It is helpful to view a skew-symmetric matrix
$\tilde M=[\tilde m_{ij}]\in M_n(\Z)$
as an oriented graph~$Q$ (i.e. a quiver) with vertex set $Q_0=\{1, 2, \ldots, n\}$ and $\tilde m_{ij}$ arrows
from $i$ to $j$ if $\tilde m_{ij}>0$ (resp. from $j$ to $i$ if $\tilde m_{ij}<0$). If a group ${\Gamma}$ acts on $Q$, one can associate with it a new matrix $M$ indexed by the orbit set $Q_0/{\Gamma}$, by defining $m_{\ig \jg}$ as the number of arrows of $Q$ between a fixed vertex $j$ of the orbit $\jg$ and any vertex of the orbit~$\ig$ (counted positively if the arrows go from  $\ig$ to $\jg$ and negatively if they go from $\jg$ to $\ig$). It is easy to see that $M$ is skew-symme\-tri\-za\-ble. For example, if
$$\tilde M = \left( \begin{matrix}
			0 & 0 & 0 & 1 \\  
			0 & 0 & 0 & 1 \\
			0 & 0 & 0 & 1 \\
			-1 & -1 & -1 & 0
                    \end{matrix}
		    \right)$$
then the quiver $Q$ is of type $D_4$
$$\xymatrix{
 1 \ar[dr] & \\
 2 \ar[r] & 4 \\
 3. \ar[ur] & 
}$$
There is an arrow from $1$ to the orbit of $4$ hence $m_{\mathbf{4} \mathbf{1}} = -1$ and there are three arrows from the orbit of $1$ to $4$ hence $m_{\mathbf{1} \mathbf{4}} = 3$. Thus we obtain the matrix
$$M = \left( \begin{matrix}
 		0 & 3 \\
		-1 & 0
             \end{matrix}
	\right)$$
of type $G_2$
$$\xymatrix@R=.5cm@C=.5cm{
     \mathbf{1} \ar@3{-}[rr] & \text{\LARGE <} & \mathbf{4}.  
}$$

Hence the action of a group ${\Gamma}$ on a skew-symmetric matrix $\tilde M$ gives rise to a skew-symme\-tri\-za\-ble matrix $M$. If $\tilde M$ is the initial seed of a cluster algebra $\tilde \Ar$ categorified as before by a category $\Cr$, it is natural to try to categorify the cluster algebra $\Ar$ with seed $M$ by a category $\Cr'$ constructed from $\Cr$ and the group~${\Gamma}$.

This leads to study $k$-additive categories $\Cr$ on which a group ${\Gamma}$ acts by auto-equivalences. In this situation, one can form a category $\Cr {\Gamma}$ whose objects are pairs $(X,(\psi_g)_{g\in {\Gamma}})$ consisting of an object $X$ of $\Cr$ isomorphically invariant under $\Gamma$, together with a family of isomorphisms $\psi_g$ from $X$ to each of its images by the elements $g$ of ${\Gamma}$. One also requires that the $\psi_g$ satisfy natural compatibility conditions. The category $\Cr {\Gamma}$ will be called the ${\Gamma}$-equivariant category.
 
One then shows useful results of transfer. For example:
\begin{itemize}
\item if $\Cr$ is abelian, then $\Cr {\Gamma}$ is abelian; 
\item if $\Cr$ is exact and if for all $g\in {\Gamma}$ the auto-equivalence of $\Cr$ associated to $g$ is exact, then $\Cr {\Gamma}$ is exact; 
\item if $H$ is a normal subgroup of ${\Gamma}$ then ${\Gamma}/H$ acts on $\Cr H$ and one has an equivalence of categories $(\Cr H)({\Gamma}/H) \simeq \Cr {\Gamma}$. 
\end{itemize} 
We also prove that $\Cr {\Gamma}$ can be endowed with a natural action of the category $\md k[{\Gamma}]$ of representations of ${\Gamma}$ over $k$.

The categories $\Cr$ used by Gei\ss, Leclerc, Schröer and Fu, Keller for modelling cluster algebras always have the following essential properties. They are Frobenius categories (i.e. exact categories with enough injectives and projectives, and the injectives and projectives are the same), and they satisfy
\[
\Ext^1_\Cr(X,Y) \simeq \Ext^1_\Cr(Y,X)^*,
\]
functorially in $X$ and $Y$. To summarize, such a category $\Cr$ is said to be $2$-Calabi-Yau. In this framework, the notion of cluster-tilting object introduced by \cite{Iy07} is very useful. An object $X$ is cluster-tilting if it is rigid, that is, if $\Ext^1_\Cr(X,X)= 0$ and if every object $Y$ satisfying $\Ext^1_\Cr(X,Y)= 0$ is in the additive envelope of $X$. If $\Cr$ categorifies a skew-symmetric cluster algebra $\tilde \Ar$, the cluster-tilting objects model the clusters of $\tilde \Ar$, and their indecomposable direct summands correspond to cluster variables.

Our principal transfer result shows that if $\Cr$ is $2$-Calabi-Yau, then $\Cr {\Gamma}$ is also $2$-Calabi-Yau. Moreover, the two natural adjoint functors linking $\Cr$ and $\Cr {\Gamma}$ induce reciprocal bijections between isomorphism classes of ${\Gamma}$-stable cluster-tilting objects of $\Cr$ and isomorphism classes of $\md k[{\Gamma}]$-stable cluster-tilting objects of $\Cr {\Gamma}$.

In order to apply these general results to the examples studied by Gei\ss, Leclerc and Schröer, we have computed explicitly the category $\Cr {\Gamma}$ in several cases \cite{De-1}. This includes in particular the case when $\Cr$ is the module category of a preprojective algebra.

\subsection{Categorification of skew-symme\-tri\-za\-ble cluster algebras}

Consider a $2$-Calabi-Yau category $\Cr$ on which acts a finite group ${\Gamma}$. For completing the categorification, one needs to develop a theory of mutations of $\md k[{\Gamma}]$-stable cluster-tilting objects of $\Cr {\Gamma}$, or, equivalently, of  ${\Gamma}$-stable cluster-tilting objects of $\Cr$. In the skew-symmetric case (which can be seen as the case where ${\Gamma}$ is trivial), it is known that such a theory is possible as soon as there exists a cluster-tilting object whose associated quiver has neither loops nor $2$-cycles. We introduce in the general case the concepts of $\md k[{\Gamma}]$-loops and of $\md k[{\Gamma}]$-$2$-cycles for a $\md k[{\Gamma}]$-stable cluster-tilting object. One then shows that if $\Cr {\Gamma}$ admits a $\md k[{\Gamma}]$-stable cluster-tilting object having neither $\md k[{\Gamma}]$-loops nor $\md k[{\Gamma}]$-$2$-cycles, all $\md k[{\Gamma}]$-stable cluster-tilting objects also have this property. Under this hypothesis, one can define a mutation operation. More precisely, if $T$ is $\md k[{\Gamma}]$-stable cluster-tilting and if $\bar{X}$ is the $\md k[{\Gamma}]$-orbit of an indecomposable non projective direct summand $X$ of $T$, one constructs another $\md k[{\Gamma}]$-stable cluster-tilting object $T'$ obtained by replacing $\bar{X}$ by the $\md k[{\Gamma}]$-orbit $\bar{Y}$ of another indecomposable object $Y$. One denotes $\mu_{\bar{X}}(T) = T'$. One can also associate to $T$ a skew-symme\-tri\-za\-ble matrix $B(T)$ whose rows are indexed by the $\md k[{\Gamma}]$-orbits $\bar{X}$ of indecomposable summands of $T$ and the columns by the $\md k[{\Gamma}]$-orbits $\bar{X}$ of indecomposable non projective factors of $T$.
The coefficients $b_{\bar{X} \bar{Y}}$ are the numbers of arrows in the Gabriel quiver of $\End_\Cr(T)$ between a fixed indecomposable object $Y$ of $\bar Y$ and any indecomposable object $X$ of $\bar X$ (the arrows from $\bar X$ to $\bar Y$ being counted positively and the arrows from $\bar Y$ to $\bar X$ being counted negatively). We then show (see theorem \ref{cfmutfz}):
\begin{thmn}{A}
 The mutation of $\Gamma$-stable cluster-tilting objects of $\Cr$ agrees with the mutation defined combinatorially by Fomin and Zelevinsky for skew-symme\-tri\-za\-ble matrices. That is,
\[
B(\mu_{\bar X}(T)) = \mu_{\bar X}(B (T)),
\]
where, in the right-hand side, by abuse of notation, $\mu_{\bar X}$ is the matrix mutation of Fomin and Zelevinsky.
\end{thmn}
Via the above-mentioned bijection between $\md k[{\Gamma}]$-stable cluster-tilting objects of $\Cr {\Gamma}$ and ${\Gamma}$-stable cluster-tilting objects of $\Cr$, one can associate to each ${\Gamma}$-stable cluster-tilting object $T$ of $\Cr$ a matrix which will be also denoted by $B(T)$.
  
Finally, in order to attach to $\Cr$ and ${\Gamma}$ a cluster algebra, we introduce a notion of ${\Gamma}$-equivariant cluster character. In the skew-symmetric case, according to the work of Caldero-Chapoton \cite{CaCh06}, of Caldero-Keller \cite{CaKe08}, \cite{CaKe06}, of Palu \cite{Pa08}, of Dehy-Keller \cite{DeKe08} and of Fu-Keller \cite{FuKe} one can assign to every object $X$ of $\Cr$ a Laurent polynomial in the cluster variables of an initial seed of the cluster algebra $\Ar$ categorified by~$\Cr$. If this initial seed is ${\Gamma}$-stable, one can identify the cluster variables which belong to the same ${\Gamma}$-orbit. This specialization of the Laurent polynomial associated to $X$ only depends on the ${\Gamma}$-orbit $\bar{X}$ of $X$, and is denoted by $P_{\bar{X}}$. One deduces from the previous construction that if $T$ is a ${\Gamma}$-stable cluster-tilting object in $\Cr$, and if $\Ar$ is the cluster algebra whose initial seed has skew-symme\-tri\-za\-ble matrix $B(T)$, then the cluster variables of $\Ar$ are of the form $P_{\bar{X}}$ where $\bar{X}$ is the orbit of an indecomposable summand of a ${\Gamma}$-stable cluster-tilting object of $\Cr$. In this situation, we say that the pair $(\Cr, \Gamma)$ is a categorification of $\Ar$. One can then generalize the result of Fu and Keller (see corollary \ref{indlin}):

\begin{thmn}{B}
 Let $B$ be an $m \times n$ matrix with skew-symme\-tri\-za\-ble principal part. If $B$ has full rank and if the cluster algebra $\Ar(B)$ has a categorification $(\Cr, \Gamma)$, then the cluster monomials are linearly independent.
\end{thmn}

As a result, we obtain a proof of conjecture \ref{cFZ} for a large family of skew-symme\-tri\-za\-ble cluster algebras.

\subsection{Applications}

Finally, we give new families of examples of categorification of cluster algebras. Let $G$ be a semi-simple connected and simply-connected Lie group of simply-laced Dynkin diagram $\Delta$, and let $\Lambda$ be the associated preprojective algebra. Gei\ss, Leclerc and Schröer have shown that the subcategories $\Sub I_J$ of $\md \Lambda$ induce cluster structures on the multi-homogeneous coordinate rings of partial flag varieties associated to $G$ \cite{GeLeSc08}. Our work allows to extend this result to the case where $G$ corresponds to a non simply-laced Dynkin diagram. In particular, one obtains a proof of the conjecture \ref{cFZ} for these cluster algebras, and one can complete the classification of partial flag varieties whose cluster structure is of finite type  (i.e. admit a finite number of clusters). In particular, this proves the conjecture of \cite[\S 14]{GeLeSc08}.

Let $G$ be a Kac-Moody group of symmetric Cartan matrix, and let $\Lambda$ be the associated preprojective algebra. Gei\ss, Leclerc and Schröer have introduced certain subcategories $\Cr_M$ of $\md \Lambda$ and shown that they induce cluster structures on the coordinate ring of some unipotent subgroups and unipotent cells of $G$ \cite{GeLeSc} (see also \cite{BuIyReSc} which gives a different definition of similar subcategories). Our work allows to extend these results to the Kac-Moody groups $G$ with symme\-tri\-za\-ble Cartan matrices. In particular, one obtains for all these examples a proof of the conjecture \ref{cFZ}. As a particular case of this construction, we get (see theorem \ref{acyc}):
\begin{thmn}{C}
 For every acyclic cluster algebra without coefficient $\Ar$, there is a category $\Cr$ and a finite group ${\Gamma}$ acting on $\Cr$ which categorify $\Ar$ up to specialization of coefficients to $1$. This holds in particular for cluster algebras of finite type.
\end{thmn}
Note that the works of Gei\ss, Leclerc and Schröer use as a crucial fact the existence of the dual semicanonical basis constructed by Lusztig for the coordinate ring of a maximal unipotent subgroup of $G$. But, when $G$ is not of simply-laced type, there is no available construction of semicanonical bases. Our result can be interpreted as giving a part of the dual semicanonical basis in the non simply-laced case, namely the set of cluster monomials.

\section{Equivariant categories}

\label{equiv}
For references about monoidal categories and module categories over a monoidal category, see for example \cite{BaKi01}, \cite{ChPr94}, \cite{Ka95} and \cite{Os03}.

\subsection{Definitions and first properties}
Let $k$ be a field, $\Cr$ a $k$-category, $\Hom$-finite and Krull-Schmidt (which means that the endomorphism rings of indecomposable objects are local, or equivalently, that every idempotent splits). Let ${\Gamma}$ be a finite group whose cardinality is not divisible by the characteristic of $k$. Let $\G = \md k({\Gamma})$ be the monoidal category of $k({\Gamma})$-modules, where $k({\Gamma})$ is the Hopf algebra of $k$-valued functions on the group ${\Gamma}$. Remark that the simple objects in $\G$ are the one-dimensional $k(\Gamma)$-modules given by evaluation maps at each element $g$ of $\Gamma$. If $g \in {\Gamma}$, the corresponding simple object in $\G$ will be denoted by $\g$. With this notation, it is easy to check that the monoidal structure is simply $\g \tens \h = \g\h$, where for $g,h \in \Gamma$, we denote by $\g\h$ the simple $k(\Gamma)$-module corresponding to $gh \in \Gamma$.

\begin{df}
 An \emph{action} of ${\Gamma}$ on $\Cr$ is a structure of $\G$-module category on $\Cr$.
\end{df}

\begin{rem}
 If one considers, as in \cite[p. 254]{ReRi85}, a group morphism $\rho$ from ${\Gamma}$ to the group of autofunctors of $\Cr$, one obtains a strict $\G$-module structure by setting $\g \tens - = \rho(g)$. 
\end{rem}

We now introduce a category of ${\Gamma}$-invariant objects of $\Cr$. The naive idea of considering the full subcategory of $\Cr$ of invariant objects does not work because almost none of the desired properties are preserved.

\begin{df}
 Let $\Cr$ be endowed with an action of ${\Gamma}$. The \emph{${\Gamma}$-equivariant category of $\Cr$} is the category whose objects are pairs $(X, \psi)$, where $X \in \Cr$, and $\psi = (\psi_g)_g \in {\Gamma}$ is a family of isomorphisms $\psi_g: \g \tens X \rightarrow X$ such that, for every $g, h \in {\Gamma}$, the following diagram commutes:
 $$\xymatrix{
  \g \tens (\h \tens X) \ar[rr]^{\id_\g \tens \psi_h} & & \g \tens X \ar[d]_{\psi_g} \\
  \g\h \tens X \ar[u]^\alpha \ar[rr]^{\psi_{gh}} & & X.
 }$$
 Here, $\alpha$ denotes the $\Gamma$-module structural isomorphism. We also assume that $\psi_e: \un \tens X \rightarrow X$ is the structural isomorphism of the $\G$-module category whenever $e$ is the neutral element of ${\Gamma}$.
 
 The morphisms from an object $(X, \psi)$ to an object $(Y, \chi)$ are the morphisms $f$ from $X$ to $Y$ such that for every $g \in {\Gamma}$, the following diagram commutes:
 $$\xymatrix{
  \g \tens X \ar[d]_{\id_\g \tens f} \ar[r]^{\psi_g} & X \ar[d]_f \\
  \g \tens Y \ar[r]^{\chi_g} & Y.
 }$$
\end{df}

\begin{nt}
 \begin{itemize}
  \item In the sequel, in particular in diagrams, we will denote by $\psi$ every arrow of the form $\id \tens \psi_g$.
  \item Every $\G$-module structural isomorphism will be denoted by $\alpha$ in diagrams.
  \item The ${\Gamma}$-equivariant category of $\Cr$ will be denoted by $\Cr {\Gamma}$.
 \end{itemize}
\end{nt}

\begin{rem}
 This category is equivalent to the ``skew group category'' considered by Reiten and Riedtmann in \cite[p. 254]{ReRi85}. We have found our definition easier to handle because it does not require to use a Karoubi envelope. Moreover, it permits to deal with non-strict actions, which is more practical in certain cases. For more details about this problem, in particular for the proof of the equivalence, see \cite{De08-1}.
\end{rem}

\begin{exs}
 \label{exequiv}
 \begin{enumerate}
  \item If $\Cr = \vect k$ is the category of finite dimensional $k$-vector spaces, and the action of ${\Gamma}$ on $\Cr$ is trivial, then $\Cr {\Gamma} \simeq \md k[{\Gamma}]$.
  \item If $A$ is a $k$-algebra and $\Gamma$ acts on $A$, ${\Gamma}$ acts naturally on $\md A$. Then, one has $(\md A) {\Gamma} \simeq \md (AG)$ where $AG$ is the skew group algebra defined in \cite{ReRi85} (see also section \ref{acc}).
  \item If ${\Gamma}$ is a cyclic group, the category $\Cr \Gamma$ is the same as the one considered in \cite[chapter 11]{Lu93}.
 \end{enumerate}
\end{exs}

The following proposition is an easy generalization of a proposition of Gabriel \cite[p. 94--95]{Ga81}:

\begin{prop}
 \label{relcycl}
 If ${\Gamma} = \langle g_0 \rangle$ is cyclic of order $n \in \N$, and if every element of $k$ has an $n$-th root in $k$, then every $X \in \Cr$ such that $X \simeq \g_0 \tens X$ has a lift $(X, \psi)$ in $\Cr {\Gamma}$.
\end{prop}

Note that proposition \ref{relcycl} does not generalize to a non cyclic group (see e.g. \cite[ex. 2.1.19]{De08-1}). 

% \begin{ex}
%  Let $A$ the path algebra of the quiver
%  $$\xymatrix{
%   & 1 \ar[d]^a & \\
%   4 \ar[r]^b & 5 & 2 \ar[l]^d \\
%   & 3 \ar[u]^c &
%  }$$
%  and ${\Gamma} = \Z/2\Z \times \Z/2\Z = \langle g,h \rangle$. One lets ${\Gamma}$ act on $A$ in the following way: $g$ exchanges $1$, $a$ and $3$, $c$, and it leaves invariant other vertices; $h$ exchanges $2$, $b$ and $4$, $d$ and leaves invariant other vertices. Following example \ref{exequiv} (the second one), ${\Gamma}$ acts on $\md A$. Consider the representation $X$:
%  $$\xymatrix{ & k \ar[d]^{f_a}& \\
%  k \ar[r]^{f_b} & k^2 & k \ar[l]^{f_d} \\
%  & k \ar[u]^{f_c} &
%  }$$
%  with $f_a = \vech{1}{0}$, $f_b = \vech{1}{1}$, $f_c = \vech{0}{1}$, $f_d = \vech{1}{-1}$. One remarks that there are two isomorphisms of order $2$ from $\g \tens X$ to $X$, the restrictions of which to vertex $5$ are
%  $$\psi_1 = \mat{0}{1}{1}{0} \quad \text{and} \quad \psi_2 = \mat{0}{-1}{-1}{0}$$
%  and two isomorphisms of order $2$ from $\h \tens X$ to $X$, the restrictions of which to vertex $5$ are
%  $$\phi_1 = \mat{1}{0}{0}{-1} \quad \text{and} \quad \phi_2 = \mat{-1}{0}{0}{1}.$$
%  It follows that, for all $g \in {\Gamma}$, $\g \tens X \simeq X$, but there is no lift $(X, \psi) \in \Cr {\Gamma}$ because $\psi_i$ and $\phi_j$ do not commute for $i, j \in \{1,2\}$. 
% \end{ex}

\begin{lem}
 \label{propcrg}
 \begin{enumerate}
  \item The category $\Cr {\Gamma}$ is $k$-additive, $\Hom$-finite and Krull-Schmidt.
  \item If $\Cr$ is abelian, then $\Cr {\Gamma}$ is also abelian.
  \item If $\Cr$ is exact and if for all $g \in {\Gamma}$, the functor $\g \tens -: \Cr \rightarrow \Cr$ is exact (one will say that the action is exact), then $\Cr {\Gamma}$ is exact with admissible short exact sequences of the form $0 \rightarrow (X, \psi) \xrightarrow{f} (X', \psi') \xrightarrow{g} (X'', \psi'') \rightarrow 0$ such that $0 \rightarrow X \xrightarrow{f} X' \xrightarrow{g} X'' \rightarrow 0$ is an admissible short exact sequence in $\Cr$.
 \end{enumerate}
\end{lem}

\begin{demo}
 \begin{enumerate}
  All these points are clear from the functoriality of the various constructions (kernel, splitting, \dots). \cqfd
 \end{enumerate}
\end{demo}

\begin{df}
 Let $n \in \N$. An \emph{$n$-associativity} is a functor from $\G^n \times \Cr$ to $\Cr$ built from the bifunctor $\tens$ and the object $\un \in \G$. For instance, $(- \tens -) \tens (\un \tens -)$ and $- \tens (- \tens -)$ are $2$-associativities. An \emph{associativity} is an $n$-associativity for some $n$.

 Let $H$ be a subgroup of ${\Gamma}$ and $(X, \psi) \in \Cr H$. Let $g_1, g_2, \dots, g_n, g'_1, g'_2, \dots, g'_m \in {\Gamma}$.
 Let $A_1$ be an $n$-associativity and $A_2$ an $m$-associativity. A \emph{$\psi$-structural isomorphism} from $A_1(\g_1, \dots, \g_n, X)$ to $A_2(\g'_1, \dots, \g'_m, X)$ is any isomorphism composed of structural isomorphisms of the $\G$-module structure, and of isomorphisms of the form $\Phi(\psi_h)$ or $\Phi(\psi_h^{-1})$ where $h \in H$ and $\Phi$ is a functor constructed from an associativity and objects of $\G$.
\end{df}

\begin{lem}[coherence]
 \label{lemcoh}
 Let $H$ be a subgroup of ${\Gamma}$. Let $(X, \psi) \in \Cr H$. Let $$g_1, g_2, \dots, g_n, g'_1, g'_2, \dots, g'_m \in {\Gamma}.$$ Let also $A_1$ be an $n$-associativity and $A_2$ be an $m$-associativity. Then all $\psi$-structural isomorphisms from $A_1(\g_1, \g_2, \dots, \g_n, X)$ to $A_2(\g'_1, \g'_2, \dots, \g'_m, X)$ are equal.
\end{lem}

\begin{demo}
 By invertibility of the $\psi$-structural isomorphisms, one can suppose that $A_1 = A_2$ and that one of the two $\psi$-structural isomorphisms is the identity. Let $f$ be a $\psi$-structural isomorphism from $A_1(\g_1, \g_2, \dots, \g_n, X)$ to itself. By using commutative diagrams built from diagrams of the form
 $$\xymatrix{
  X \ar[d]_\alpha & \h \tens X \ar[d]_\alpha \ar[l]_\psi \\
  \h \tens (\h^{-1} \tens X) \ar[d]_\psi & \h \tens (\h^{-1} \tens (\h \tens X)) \ar[d]_\alpha \ar[l]_\psi \\
  \h \tens X & \h \tens (\un \tens X) \ar@/_/[l]_\psi \ar@/^/[l]^\alpha 
 }$$
 with $h \in H$, one can suppose that $f$ contains only positive powers of $\psi$. By using commutative diagrams of the form
 $$\xymatrix{
  \Phi_1(X) \ar[r]^\alpha & \Phi_2(\h \tens X) \ar[r]^\psi & \Phi_2(X) \\
  \Phi_1(\h' \tens X) \ar[u]^\psi \ar[r]^\alpha & \Phi_2(\h \tens (\h' \tens X)) \ar[u]^\psi \ar[r]^\alpha & \Phi_2(\h\h' \tens X) \ar[u]^\psi
 }$$
 where $h, h' \in H$ and $\Phi_1$ and $\Phi_2$ are functors made from an associativity and objects of $\G$, one can suppose that $f$ is of the form $\alpha_1 \Phi(\psi_h) \alpha_2$ where $\alpha_1$ and $\alpha_2$ are structural morphisms of $\Cr$ and $\Phi$ is made from an associativity and objects of $\G$ and $h \in H$. Then, as $\alpha_1 \Phi(\psi_h) \alpha_2$ goes from $A_1(\g_1, \g_2, \dots, \g_n, X)$ into itself, $h$ is the neutral element of $H$ and as $\psi_e$ is structural, $f$ is a structural morphism of the $\G$-module structure which implies the result, by the MacLane coherence lemma (see \cite{Ma71}). \cqfd
\end{demo}

\begin{prop}
 \label{actssgrp}
 Let $H$ be a normal subgroup of ${\Gamma}$. The action of ${\Gamma}$ on $\Cr$ induces an action of $H$ on $\Cr$. Then:
 \begin{enumerate}
  \item The action of ${\Gamma}$ on $\Cr$ extends to an action of ${\Gamma}$ on $\Cr H$.
  \item For every $h \in H$, there is an isomorphism of functors from $\h \tens -$ to $\id_{\Cr H}$.
  \item The action of ${\Gamma}$ on $\Cr H$ induces an action of ${\Gamma}/H$ on $\Cr H$.
  \item There is an equivalence of categories between $(\Cr H) ({\Gamma}/H)$ and $\Cr {\Gamma}$.
 \end{enumerate}
\end{prop}

\begin{demo}
 \begin{enumerate}
  \item Let $(X, \psi) \in \Cr H$. If $g \in {\Gamma}$ and $h \in H$, define $(\g \tens \psi)_h$ such that the following diagram commutes:
  $$\xymatrix{
   \g \tens X & \g \tens (\g^{-1}\h\g \tens X) \ar[dl]^\alpha \ar[l]_\psi \\
   \h \tens (\g \tens X) \ar[u]^{(\g \tens \psi)_h} & 
  }$$
  For $g \in {\Gamma}$ and $h, h' \in H$, the following diagram is commutative from lemma \ref{lemcoh} because all arrows are $\psi$-structural:
  $$\xymatrix{
   \h \tens (\h' \tens (\g \tens X)) \ar[rr]^{\h \tens (\g \tens \psi)_{h'}} & & \h \tens (\g \tens X) \ar[d]_{(\g \tens \psi)_h} \\
   \h\h' \tens X \ar[u]^\alpha \ar[rr]^{(\g \tens \psi)_{hh'}} & & \g \tens X
  }$$
  which means that $\g \tens (X, \psi) = (\g \tens X, \g \tens \psi) \in \Cr H$ (it is easy to see that $(\g \tens \psi)_e$ is structural). The following diagram commutes for $g, g' \in {\Gamma}$ and $h \in H$:
  $$\xymatrix{
   \h \tens (\g \tens (\g' \tens X)) \ar[rr]^{(\g \tens (\g' \tens \psi))_h} \ar[d]_\alpha & & \g \tens (\g' \tens X) \ar[d]_\alpha \\
   \h \tens (\g \g' \tens X) \ar[rr]^{(\g\g' \tens \psi)_h} & & \g \g' \tens X.
  }$$
  Therefore the structural isomorphisms of $\Cr$ are also structural isomorphisms of $\Cr H$ which leads to the conclusion.
  \item For $h, h' \in H$, for the same reason, the following diagram commutes:
  $$\xymatrix{
   \h' \tens (\h \tens X) \ar[d]_{(\h \tens \psi)_{h'}} \ar[r]^\psi & \h' \tens X \ar[d]_\psi \\
   \h \tens X \ar[r]^\psi & X.
  }$$
  Therefore $\psi_h$ is an isomorphism from $\h \tens X$ to $X$ in the category $\Cr H$. If $(Y, \chi) \in \Cr H$ and $f: (X, \psi) \rightarrow (Y, \chi)$ is a morphism, the following diagram commutes:
  $$\xymatrix{
   \h \tens X \ar[r]^\psi \ar[d]_{\id_\h \tens f} & X \ar[d]_f \\
   \h \tens Y \ar[r]^\chi & Y.
  }$$
  If one denotes
  $$\phi_{h; Y, \chi} = \chi_h$$
  for every $(Y, \chi) \in \Cr H$, then $\phi_h$ is an isomorphism from the functor $\h \tens -$ to the functor $\id_{\Cr H}$.
  \item For $g \in {\Gamma}$, denote by $\bar g$ its class in ${\Gamma}/H$. Let ${\Gamma}_0 \subset {\Gamma}$ be a set of representatives of ${\Gamma}/H$ containing the neutral element. For every $g \in {\Gamma}_0$, let $\bar \g \tens - = \g \tens -$. Let $g, g' \in {\Gamma}_0$. There exists a unique decomposition $gg' = g''h$ with $g'' \in {\Gamma}_0$ and $h \in H$. The isomorphism of functors $\bar \alpha$, such that the diagram
  $$\xymatrix{
   \bar \g \bar \g' \tens - \ar[d]_{\bar \alpha} \ar@{=}[r] & \g'' \tens - & \g'' \tens (\h \tens -) \ar[l]_{\phi_h} \\
   \bar \g \tens (\bar \g' \tens -) \ar@{=}[r] & \g \tens (\g' \tens -) \ar[r]^\alpha & \g'' \h \tens - \ar[u]^\alpha
  }$$
  commutes, endows $\Cr H$ with an action of ${\Gamma}/H$. Indeed, all axioms of a $\G/\H$-module category are verified because every morphism considered in these axioms is $\psi$-structural and therefore the equalities are true.
  \item Let $(X, \psi', \psi'') \in (\Cr H) ({\Gamma}/H)$. For $g \in {\Gamma}$, there exists $g_0 \in {\Gamma}_0$ and $h \in H$ such that $g = g_0 h$. Let $\psi_g$ be such that the following diagram commutes:
  $$\xymatrix{
   \g \tens X \ar[r]^\alpha \ar[d]_{\psi_g} & \g_0 \tens (\h \tens X) \ar[r]^{\psi'} & \g_0 \tens X \ar@{=}[d] \\
   X & & \bar \g \tens X \ar[ll]_{\psi''}
  }$$
  Let now $g' \in {\Gamma}$, $g'_0 \in {\Gamma}_0$, $h' \in H$, $g''_0 \in {\Gamma}_0$ and $h'' \in H$ such that $g' = g'_0 h'$ and $gg' = g''_0 h''$. The following diagram commutes:
  $$\xymatrix{
   \g \tens (\g'_0 \tens (\h' \tens X)) \ar[r]^{\psi'} \ar[dd]_ \alpha & \g \tens (\g'_0 \tens X) \ar[d]_\alpha \ar[r]^{\psi''} & \g \tens X \ar[d]_\alpha \\
   & \g_0 \tens (\h \tens (\g'_0 \tens X)) \ar[r]^{\psi''} \ar[d]_{\g'_0 \tens \psi'}& \g_0 \tens (\h \tens X) \ar[d]_{\psi'} \\
   \g_0 \tens (\h \tens (\g'_0 \tens (\h' \tens X))) \ar[ur]^{\psi'} & \g_0 \tens (\g'_0 \tens X) \ar[r]^{\psi''} & \g_0 \tens X \ar[d]_{\psi''} \\
   \g''_0 \tens (\h'' \tens X) \ar[r]^{\psi'} \ar[u]^\alpha & \g''_0 \tens X \ar[r]^{\psi''} \ar[u]^{\bar \alpha} & X
  }$$
  (the two upper squares because $\alpha$ is functorial, the middle right square because $\psi''$ is a morphism from $\g'_0 \tens (X, \psi')$ to $(X, \psi')$, the lower right square by definition of $\psi''$ and the lower left pentagon because it is formed by $\psi'$-structural morphisms). Keeping only the border and composing with structural morphisms on the left, one gets the following commutative diagram:
  $$\xymatrix{
   \g \tens (\g' \tens X) \ar[r]^\psi \ar[d]_\alpha & \g \tens X \ar[d]_\psi \\
   \g \g' \tens X \ar[r]^\psi & X
  }$$
  hence $(X, \psi) \in \Cr {\Gamma}$. Let $\Phi(X, \psi', \psi'') = (X, \psi)$. Let $(Y, \chi', \chi'') \in (\Cr H) ({\Gamma}/H)$ be another object and $(Y, \chi) = \Phi(Y, \chi', \chi'')$. Let $f \in \Hom_{(\Cr H) ({\Gamma}/H)}(X, \psi', \psi''; Y, \chi', \chi'')$. The following diagram commutes: 
  $$\xymatrix{
   \g \tens X \ar[r]^\alpha \ar[d]_f \ar@/^1cm/[rrr]^\psi & \g_0 \tens (\h \tens X) \ar[d]_f \ar[r]^{\psi'} & g_0 \tens X \ar[d]_f \ar[r]^{\psi''} & X \ar[d]_f \\
   \g \tens Y \ar[r]^\alpha \ar@/_1cm/[rrr]_\chi & \g_0 \tens (\h \tens Y) \ar[r]^{\chi'} & \g_0 \tens Y \ar[r]^{\chi''} & Y
  }$$
  and, as a consequence, $f \in \Hom_{\Cr {\Gamma}}((X, \psi), (Y, \chi))$. By setting $\Phi(f) = f$, $\Phi$ is a functor. Now, if $(X, \psi) \in \Cr {\Gamma}$, for $h \in H$, let $\psi'_h = \psi_h$ and for $g \in {\Gamma}_0$, $\psi''_{\bar g_0} = \psi_{g_0}$. It is easy to check that $(X, \psi', \psi'') \in (\Cr H) ({\Gamma}/H)$ (each involved morphism is $\psi$-structural) and that $\Phi(X, \psi', \psi'') = (X, \psi)$. Finally, $\Phi$ is essentially surjective. Moreover, it is easy to see that $\Phi$ is fully faithful. Hence $\Phi$ is an equivalence of categories. \cqfd
 \end{enumerate}
\end{demo}

\subsection{A $\md k[\Gamma]$-module structure on $\Cr \Gamma$}

We denote by $k[\Gamma]$ the group algebra of $\Gamma$. This is a Hopf algebra (dual to $k(\Gamma)$), hence $\md k[\Gamma]$ is a monoidal category. An object of $\md k[\Gamma]$ will be denoted by $(V, r)$, where $V$ is a $k$-vector space, and $r : \Gamma \rightarrow \GL(V)$ a group homomorphism.

\begin{prop}
 \label{modkgmod}
 The category $\Cr {\Gamma}$ is a $\md k[{\Gamma}]$-module category in a natural way.
\end{prop}

\begin{demo}
 Let $\vect_0 k$ be the full subcategory of $\vect k$, whose objects are $k^n$ ($n \in \N$). Let $\Psi$ be the inclusion functor. It is easy to extend $\Psi$ to a monoidal equivalence of categories. Let $\Phi$ be a monoidal quasi-inverse of $\Psi$. One endows $\Cr$ with a structure of $\vect k$-module category by setting for every $V \in \vect k$ and $X \in \Cr$, $V \tens X = X^{\dim(V)}$. If $V, W \in \vect k$, $X, Y \in \Cr$, $f \in \Hom_k(V, W)$ and $g \in \Hom_\Cr(X, Y)$, one defines $f \tens g: V \tens X = X^{\dim(V)} \rightarrow W \tens Y = Y^{\dim(W)}$ by
 $$f \tens g = \left(\Phi(f)_{ij} g\right)_{1 \leq i \leq \dim(W), 1 \leq j \leq \dim(V)}.$$
  
 For every $V, W \in \vect k$ and $X \in \Cr$, $\alpha_{V, W, X}: (V \tens W) \tens X \rightarrow V \tens (W \tens  X)$ is defined by
  $$\alpha_{V, W, X; i, j, \ell} = \left\{
    \begin{array}{ll}
     \id_X & \quad \text{if} \quad i = j + \dim(V)(\ell-1), \\
     0 & \quad \text{else},.
    \end{array}
    \right.
  $$
  where $1 \leq i \leq \dim(V \tens W) = \dim(V) \dim(W)$, $1 \leq j \leq \dim(V)$ and $1 \leq \ell \leq \dim(W)$.
 
 It is now easy to check that $\Cr$ is $\vect_0 k$-module with these structural isomorphisms. As $\Phi$ is a monoidal equivalence, $\Cr$ is also $\vect k$-module.
 
 One deduces that $\Cr {\Gamma}$ is $\md k[{\Gamma}]$-module. Indeed, one remarks that for every $g \in {\Gamma}$, $V \in \vect k$ and $X \in \Cr$, one has $\g \tens (V \tens X) = \g \tens X^{\dim V} = (\g \tens X)^{\dim V} = V \tens (\g \tens X)$. If $(V, r) \in \md k[{\Gamma}]$ and $(X, \psi) \in \Cr {\Gamma}$, let $(V, r) \tens (X, \psi) = (V \tens X, r \tens \psi)$ where, for $g \in {\Gamma}$, $(r \tens \psi)_g: \g \tens (V \tens X) = V \tens (\g \tens X) \rightarrow V \tens X$ is defined by $(r \tens \psi)_g = r_g \tens \psi_g$. If one takes also $h \in {\Gamma}$, one gets the two commutative diagrams 
 $$\xymatrix{
  V \ar[r]^{r_h} & V \ar[d]_{r_g} & \g \tens (\h \tens X) \ar[rr]^{\id_\g \tens \psi_h} & & \g \tens X \ar[d]_{\psi_g} \\
  V \ar@{=}[u] \ar[r]^{r_{gh}} & V & \g\h \tens X \ar[u]^\alpha \ar[rr]^{\psi_{gh}} & & X
 }$$
 and applying the bifunctor $\tens: \vect k \times \Cr \rightarrow \Cr$ yields the commutative diagram
 $$\xymatrix{
  \g \tens (\h \tens (V \tens X)) \ar[rr]^{\id_\g \tens (r \tens \psi)_h} & & \g \tens (V \tens X) \ar[d]_{(r \tens \psi)_g} \\
  \g\h \tens (V \tens X) \ar[u]^\alpha \ar[rr]^{(r \tens \psi)_{gh}} & & V \tens X
 }$$
 hence $(V \tens X, r \tens \psi) \in \Cr {\Gamma}$. 
 
 Moreover, if $(V', r') \in \md k[{\Gamma}]$, $(Y, \chi) \in \Cr {\Gamma}$, $f \in \Hom_{\md k[{\Gamma}]}((V, r), (V', r'))$ and $f' \in \Hom_{\Cr {\Gamma}}((X, \psi), (Y, \chi))$, the two following diagrams commute for every $g \in {\Gamma}$:
 $$\xymatrix{
  V \ar[d]_{f} \ar[r]^{r_g} & V \ar[d]_{f} & \g \tens X \ar[d]_{\id_g \tens f'} \ar[r]^{\psi_g} & X \ar[d]_{f'} \\
  V' \ar[r]^{r'_g} & V' & \g \tens Y \ar[r]^{\chi_g} & Y
 }$$
 which shows by applying the bifunctor $\tens: \vect k \times \Cr \rightarrow \Cr$ that $f \tens f' \in   \Hom_{\Cr {\Gamma}}((V, r) \tens (X, \psi), (V', r') \tens (Y, \chi))$. This finishes the proof that $\Cr {\Gamma}$ is a $\md k[{\Gamma}]$-module category. \cqfd
\end{demo}

\begin{rem}
 The previous structure does not depend on the choice of $\Phi$ up to isomorphism.
\end{rem}

\subsection{A $k[{\Gamma}]$-linear structure on the equivariant category}
The aim of this section is to define new morphisms spaces on $\Cr {\Gamma}$ which are $k[{\Gamma}]$-modules. These new structures will be written in bold face. This gives a new category closely related to $\Cr {\Gamma}$. The main relationships between the two categories will be outlined.

\begin{nt}
 One denotes by $F$ the forgetful functor from $\Cr {\Gamma}$ to $\Cr$.
\end{nt}

Recall this classical lemma:

\begin{lem} \label{adjoncf}
 There is an isomorphism of trifunctors from $(\md k[{\Gamma}])^3$ to $\vect k$:
 $$\Hom_{\md k[{\Gamma}]}(?_1 \tens ?_2^*,?_3) \simeq \Hom_{\md k[{\Gamma}]}(?_1, ?_2 \tens ?_3).$$
 If $r$ is a $k[{\Gamma}]$-module, $r^*$ denotes its contragredient, or dual representation.
\end{lem}

Here is an easy lemma (for a detailed proof, see \cite[lemme 2.1.16]{De08-1}) :

\begin{lem}
 \label{homtens}
 With the previous notations, there is an isomorphism of quadrifunctors
 $$\Hom_\Cr(?_1 \tens -_1,?_2 \tens -_2) \simeq?_1^* \tens?_2 \tens \Hom_\Cr(-_1, -_2)$$
 where the $-$ are variables of $\Cr$ and the $?$ are variables of $\vect k$.
\end{lem}

 Let $(X, \psi), (Y, \chi) \in \Cr {\Gamma}$. As a $k$-vector space, let $$\uHom_{\Cr {\Gamma}}((X, \psi), (Y, \chi)) = \Hom_\Cr (X, Y).$$ If $g \in {\Gamma}$ and $f \in \uHom_{\Cr {\Gamma}}(X, Y)$, define $gf \in \uHom_{\Cr {\Gamma}}(X, Y)$ by the following commutative diagram:
 $$\xymatrix{
  X \ar[r]^{gf} & Y \\
  \g \tens X \ar[r]^{\id_\g \tens f} \ar[u]^{\psi_g} & \g \tens Y \ar[u]^{\chi_g}
 }$$
 Thus, one will prove in proposition \ref{rephom} that $\uHom_{\Cr {\Gamma}}(X, Y)$ acquires the structure of a $k[{\Gamma}]$-module.
 If ${\Gamma}: (X', \psi') \rightarrow (X, \psi)$, ${\Gamma}': (Y, \chi) \rightarrow (Y', \chi')$ are morphisms in $\Cr {\Gamma}$, define
 $$\uHom_{\Cr {\Gamma}}({\Gamma}, {\Gamma}') = \Hom_\Cr(F {\Gamma}, F {\Gamma}').$$

\begin{prop}
 \label{rephom}
 Defined in this way, $\uHom_{\Cr {\Gamma}}$ is a bifunctor from $\Cr {\Gamma} \times \Cr {\Gamma}$ to $\md k[{\Gamma}]$ contravariant in the first variable and covariant in the second one which satisfies:
 \begin{enumerate}
  \item for $X, Y, Z \in \Cr {\Gamma}$ the composition 
   $$\circ: \Hom_\Cr(FY, FZ) \tens \Hom_\Cr(FX, FY) \rightarrow \Hom_\Cr(FX, FZ)$$
    is a morphism of $k[{\Gamma}]$-modules 
   $$\circ: \uHom_{\Cr {\Gamma}}(Y, Z) \tens \uHom_{\Cr {\Gamma}}(X, Y) \rightarrow \uHom_{\Cr {\Gamma}}(X, Z)\text{;}$$
  \item \label{rephom2} there is an isomorphism of quadrifunctors
   $$\uHom_{\Cr {\Gamma}}(?_1 \tens -_1,?_2 \tens -_2) \simeq?_1^* \tens?_2 \tens \uHom_{\Cr {\Gamma}}(-_1, -_2)$$
   where the $?$ are variables in $\md k[{\Gamma}]$ and the $-$ are variables in $\Cr {\Gamma}$;
  \item there is an isomorphism of quadrifunctors
   $$\Hom_{\Cr {\Gamma}}(?_1 \tens -_1,?_2 \tens -_2) \simeq \Hom_{\md k[{\Gamma}]}(?_1 \tens?_2^*, \uHom_{\Cr {\Gamma}}(-_1, -_2))$$
   where the $?$ are variables in $\md k[{\Gamma}]$ and the $-$ are variables in $\Cr {\Gamma}$.
 \end{enumerate}
 In particular it endows $\Cr {\Gamma}$ with the structure of a $\md k[{\Gamma}]$-linear category.
\end{prop}

\begin{demo}
 For $g, h \in {\Gamma}$, the following diagram commutes:
 $$\xymatrix{
  & X \ar[rr]^{g(hf)} & & Y & \\
  & \g \tens X \ar[u]^\psi \ar[rr]^{\id_\g \tens hf} & &  \g \tens Y \ar[u]^\chi & \\
  \g\h \tens X \ar[uur]^\psi \ar@/_.5cm/[rrrr]_{\id_{\g\h} \tens f} & \g \tens (\h \tens X) \ar[l]_\alpha \ar[u]^\psi \ar[rr]^{\id_\g \tens (\id_\h \tens f)} & & \g \tens (\h \tens Y) \ar[u]^\chi \ar[r]^\alpha & \g\h \tens Y \ar[luu]_\chi
 }$$
 so $g(hf) = (gh)f$ and $\uHom_{\Cr {\Gamma}}(X, Y)$ is a representation of ${\Gamma}$.
 
 If ${\Gamma}: (X', \psi') \rightarrow (X, \psi)$, ${\Gamma}': (Y, \chi) \rightarrow (Y', \chi')$ are morphisms in $\Cr {\Gamma}$, and if $f \in \Hom_\Cr(X, Y)$, the following diagram commutes for every $g \in {\Gamma}$:
 $$\xymatrix{
  & X' \ar[rrr]^{g \Hom_\Cr({\Gamma}, {\Gamma}')(f)} \ar[ldd]_{\Gamma} & & & Y' \\
  & \g \tens X' \ar[u]_{\psi'} \ar[d]^{\id_\g \tens {\Gamma}} \ar[rrr]^{\id_\g \tens \Hom_\Cr ({\Gamma}, {\Gamma}')(f)} & & & \g \tens Y' \ar[u]^{\chi'} & \\
  X \ar@/_.5cm/[rrrrr]_{gf} & \g \tens X \ar[l]_\psi \ar[rrr]^{\id_\g \tens f} & & & \g \tens Y \ar[u]^{\id_\g \tens {\Gamma}'} \ar[r]^\chi & Y \ar[uul]_{{\Gamma}'}
 }$$
 and finally, $g \Hom_\Cr({\Gamma}, {\Gamma}')(f) = \Hom_\Cr({\Gamma}, {\Gamma}')(gf)$. 
 
 Hence $\Hom_\Cr({\Gamma}, {\Gamma}')$ is a morphism from $\uHom_{\Cr {\Gamma}} ((X, \psi), (Y, \chi))$ to $\uHom_{\Cr {\Gamma}} ((X', \psi'), (Y', \chi'))$, morphism which will be denoted $\uHom_{\Cr {\Gamma}}({\Gamma}, {\Gamma}')$ turning $\uHom_{\Cr {\Gamma}}$ into a bifunctor.

 Let us now prove the three additional properties:
 \begin{enumerate}
  \item is clear. 
  \item It is enough to show that the isomorphism of quadrifunctors ${\Gamma}$ defined in lemma \ref{homtens}
 remains an isomorphism. Let $g \in {\Gamma}$ and $(V, r), (V', r') \in \md k[{\Gamma}]$. By definition of a morphism of functors, the lower square of the following diagram commutes:
 $$\xymatrix{
  V^* \tens V' \tens \Hom_\Cr(X, Y) \ar[r]^{\Gamma} \ar[d]_{\id \tens \id \tens (\g \tens -)} & \Hom_\Cr(V \tens X, V' \tens Y) \ar[d]^{\g \tens -} \\
  V^* \tens V' \tens \Hom_\Cr(\g \tens X, \g \tens Y) \ar[r]^{\Gamma} \ar[d]_{\tr{r_g^{-1}} \tens r'_g \tens \Hom_\Cr(\psi_g^{-1}, \chi_g)} & \Hom_\Cr(V \tens (\g \tens X), V' \tens (\g \tens Y)) \ar[d]^{\Hom_\Cr(r_g^{-1} \tens \psi_g^{-1}, r'_g \tens \chi_g)} \\
  V^* \tens V' \tens \Hom_\Cr(X, Y) \ar[r]^{\Gamma} & \Hom_\Cr(V \tens X, V' \tens Y)
 }$$
 and the upper square commutes because $V \tens (\g \tens X) = \g \tens (V \tens X)$. This proves that ${\Gamma}_{V, V', X, Y}$ is a morphism of representations. 
 \item From lemma \ref{adjoncf} and (\ref{rephom2}), one gets
 \begin{align*}
  \Hom_{\md k[{\Gamma}]}(?_1 \tens?_2^*, \uHom_{\Cr {\Gamma}}(-_1, -_2)) &\simeq \Hom_{\md k[{\Gamma}]}(\un,?_1^* \tens?_2 \tens \uHom_{\Cr {\Gamma}}(-_1, -_2)) \\ &\simeq \Hom_{\md k[{\Gamma}]}(\un, \uHom_{\Cr {\Gamma}}(?_1 \tens -_1,?_2 \tens -_2))
 \end{align*}
 hence it is enough to see that $\Hom_{\md k[{\Gamma}]}(\un, \uHom_{\Cr {\Gamma}}(-_1, -_2)) \simeq \Hom_{\Cr {\Gamma}}(-_1, -_2)$. This is clear. Indeed, it is sufficient, if $f \in \Hom_{\md k[{\Gamma}]}(\un, \uHom_{\Cr {\Gamma}}((X, \psi), (Y, \chi)))$, to associate to it 
 $$f(1) \in \bigcap_{g \in {\Gamma}} \uHom_{\Cr {\Gamma}}((X, \psi), (Y, \chi))^g = \Hom_{\Cr {\Gamma}}((X, \psi), (Y, \chi))$$
 the last equality being the definition of $\Hom_{\Cr {\Gamma}}((X, \psi), (Y, \chi))$. \cqfd
 \end{enumerate}
\end{demo}

\begin{cor}
 There is an isomorphism of trifunctors
 $$\Hom_{\Cr {\Gamma}}(?^* \tens -, -) \simeq \Hom_{\Cr {\Gamma}} (-,? \tens -)$$
 where $?$ is a variable of $\md k[{\Gamma}]$ and the $-$ are variables of $\Cr {\Gamma}$. In particular, if $r \in \md k[{\Gamma}]$, the two functors from $\Cr {\Gamma}$ to itself $r \tens -$ and $r^* \tens -$ are adjoint.
\end{cor}

\begin{demo}
 Using proposition \ref{rephom}, one gets
 $$
  \Hom_{\Cr {\Gamma}}(?^* \tens -, -) \simeq \Hom_{\md k[{\Gamma}]}(?^*, \uHom_{\Cr {\Gamma}}(-, -)) 
  \simeq \Hom_{\Cr {\Gamma}}(-,? \tens -),
 $$
 which is the desired isomorphism. \cqfd
\end{demo}

\subsection{The functors $-[{\Gamma}]$ and $F$}

\begin{prop}
 For $X \in \Cr$, define
 $$\tilde X = \left(\bigoplus_{g \in {\Gamma}} \g \right) \tens X.$$
 For $g \in {\Gamma}$, denote by $\psi_g$ the unique structural isomorphism from $\g \tens \tilde X$ to $\tilde X$. Then $(\tilde X, \psi) \in \Cr {\Gamma}$.
\end{prop}

\begin{demo}
 For $g, h \in {\Gamma}$, the diagram
 $$\xymatrix{
  \g \tens (\h \tens \tilde X) \ar[rr]^{\id_\g \tens \psi_h} & & \g \tens \tilde X \ar[d]_{\psi_g} \\
  \g\h \tens \tilde X \ar[u]^\alpha \ar[rr]^{\psi_{gh}} & & \tilde X
 }$$
 commutes by unicity of the structural isomorphism from $\g \tens (\h \tens \tilde X)$ to $\tilde X$ (lemma \ref{lemcoh}). \cqfd
\end{demo}

\begin{df}
 The object $(\tilde X, \psi)$ will be denoted by $X[{\Gamma}]$.
\end{df}

\begin{lem}
 One can extend $-[{\Gamma}]$ to a functor from $\Cr$ to $\Cr {\Gamma}$ by setting
 $$f[{\Gamma}] = \left(\bigoplus_{g \in {\Gamma}} \id_\g \right) \tens f$$
 for every morphism $f$ of $\Cr$.
\end{lem}

\begin{demo}
 It is enough to see that if $f: X \rightarrow Y$ is a morphism in $\Cr$, then
 $$f[{\Gamma}] = \left(\bigoplus_{g \in {\Gamma}} \id_\g \right) \tens f$$
 is a morphism from $X[{\Gamma}]$ to $Y[{\Gamma}]$. Let $h \in {\Gamma}$. The diagram
 $$\xymatrix{
  \h \tens \left(\left(\bigoplus_{g \in {\Gamma}} \g \right) \tens X\right) \ar[rr]^{\id \tens (\id \tens f)} \ar[d]_\alpha & & \h \tens \left(\left(\bigoplus_{g \in {\Gamma}}\g \right) \tens Y\right) \ar[d]_\alpha \\
  \left(\bigoplus_{g \in {\Gamma}} \g \right) \tens X \ar[rr]^{\id \tens f} & & \left(\bigoplus_{g \in {\Gamma}}\g \right) \tens Y
 }$$
 commutes by functoriality of structural morphisms. \cqfd
\end{demo}

\begin{prop} \label{adjfc}
 \begin{enumerate} 
  \item There is an isomorphism of functors from $\Cr {\Gamma}$ to itself:
   $$(F-)[{\Gamma}] \simeq k[{\Gamma}] \tens -$$
   where $F: \Cr {\Gamma} \rightarrow \Cr$ is the forgetful functor and $k[{\Gamma}]$ denotes the regular representation of ${\Gamma}$.
  \item The functors $F: \Cr {\Gamma} \rightarrow \Cr$ and $-[{\Gamma}]: \Cr \rightarrow \Cr {\Gamma}$ are adjoint.
 \end{enumerate}
\end{prop}

\begin{demo}
 \begin{enumerate}
  \item Let $(X, \psi) \in \Cr$. Let $\chi$ be such that
   $$\left(\left(\bigoplus_{g \in {\Gamma}} \g \right) \tens X, \chi \right) = X[{\Gamma}].$$
   For $g, h \in {\Gamma}$, the following diagram commutes:
   $$\xymatrix{
    \h \tens (\g \tens X) \ar[r]^{\alpha} \ar[d]_{\h \tens \psi_g} & \h\g \tens X \ar[d]_{\psi_{hg}}\\
    \h \tens X \ar[r]^{\psi_h} & X.
   }$$
   Adding up these diagrams for $g \in {\Gamma}$, by setting $g' = hg$, the following diagram commutes:
   $$\xymatrix{
    \h \tens \left(\left( \bigoplus_g \g \right) \tens X\right) \ar[r]^{\chi_h} \ar@{=}[d] \ar@/_3cm/[ddd]_{\h \tens \bigoplus \psi_g} & \left( \bigoplus_{g'} \g' \right) \tens X \ar@{=}[d] \\
    \bigoplus_g (\h \tens (\g \tens X)) \ar[r]^{\bigoplus \alpha} \ar[d]_{\bigoplus \h \tens \psi_g} & \bigoplus_{g'} (\g' \tens X) \ar[d]_{\bigoplus \psi_{g'}}\\
    \bigoplus_g (\h \tens X) \ar[r]^{\bigoplus \psi_h} \ar@{=}[d] & \bigoplus_{g'} X \ar@{=}[d] \\    
    \h \tens \left( \bigoplus_g X\right) \ar[r]^{(k[{\Gamma}] \tens \psi)_h} &  \bigoplus_{g'} X
  }$$
  the equalities consisting of unique identifications through universal properties. As a consequence, $\bigoplus \psi_g$ is an isomorphism from $X[{\Gamma}]$ to $k[{\Gamma}] \tens (X, \psi)$. Moreover, if $f: (X, \psi) \rightarrow (Y, \psi')$ is a morphism, the following diagram commutes:
  $$\xymatrix{
   \left( \bigoplus_g \g \right) \tens X \ar[rr]^{(\bigoplus \g) \tens f} \ar@{=}[d] & & \left( \bigoplus_g \g \right) \tens Y \ar@{=}[d] \\
   \bigoplus_g (\g \tens X) \ar[rr]^{\bigoplus (\g \tens f)} \ar[d]_{\bigoplus \psi_g} & & \bigoplus_g (\g \tens Y) \ar[d]_{\bigoplus \psi'_g}\\
   \bigoplus_g X \ar[rr]^{k[{\Gamma}] \tens f = \bigoplus f} & &  \bigoplus_g Y  \\
  }$$
  the commutativity of the lower square coming from the fact that $f$ is a morphism in $\Cr {\Gamma}$. Finally, $\bigoplus \psi_g$ is a functorial isomorphism.
 \item Let $X \in \Cr$ and $(Y, \psi) \in \Cr {\Gamma}$. If $f \in \Hom_\Cr(X, Y)$, one defines
  $$\xi(f): \bigoplus_{g \in {\Gamma}} (\g \tens X) = \left(\bigoplus_{g \in {\Gamma}} \g \right) \tens X \rightarrow Y$$
  by
  $$(\xi(f))_g = \psi_g \circ (\g \tens f).$$
  
  For $g, h \in {\Gamma}$, the following diagram commutes:
  $$\xymatrix{
   \h \tens (\g \tens X) \ar[rr]^{\h \tens (\g \tens f)} \ar[d]_\alpha & & \h \tens (\g \tens Y) \ar[rr]^{\h \tens \psi_g} \ar[d]_\alpha & & \h \tens Y \ar[d]_{\psi_h} \\
   \h \g \tens X \ar[rr]^{\h\g \tens f} & & \h \g \tens Y \ar[rr]^{\psi_{hg}} & & Y.
  }$$
  By adding up these diagrams for $g\in {\Gamma}$, one obtains the following commutative diagram:
  $$\xymatrix{
   \h \tens \left(\left( \bigoplus \g \right) \tens X\right) \ar[rr]^{\h \tens \xi(f)} \ar[d]_{\chi_h} & & \h \tens Y \ar[d]_{\psi_h} \\
   \left( \bigoplus \g \right) \tens X \ar[rr]^{\xi(f)} & & Y
  }$$
  which shows that $\xi(f)$ is a morphism from $X[{\Gamma}]$ to $(Y, \psi)$. We now check that $\xi$ is functorial. Let $X' \in \Cr$ and $(Y', \psi') \in \Cr {\Gamma}$. Let $\eta: X' \rightarrow X$ be a morphism of $\Cr$ and $\theta: (Y, \psi) \rightarrow (Y', \psi')$ a morphism of $\Cr {\Gamma}$. Then, for $g \in {\Gamma}$,
  \begin{align*}
   (\theta \circ \xi(f) \circ \eta[{\Gamma}])_g &= \theta \circ \xi(f)_g \circ (\g \tens \eta) = \theta \circ \psi_g \circ (\g \tens f) \circ (\g \tens \eta) \\ &= \psi'_g \circ (\g \tens \theta) \circ (\g \tens f \eta) = \psi'_g \circ (\g \tens \theta f \eta) = \xi(\theta f \eta)_g
  \end{align*}
  hence, the following diagram commutes:
  $$\xymatrix{
   \Hom_\Cr(X, Y) \ar[rrr]^{\Hom_\Cr(\eta, F \theta)} \ar[d]_\xi & & & \Hom_\Cr(X', Y') \ar[d]_\xi \\
   \Hom_{\Cr {\Gamma}}(X[{\Gamma}], (Y, \psi)) \ar[rrr]^{\Hom_{\Cr {\Gamma}}(\eta[{\Gamma}], \theta)} & & & \Hom_{\Cr {\Gamma}}(X'[{\Gamma}], (Y', \psi'))
  }$$
  which leads to the functoriality of $\xi$. If $f' \in \Hom_{\Cr {\Gamma}}(X[{\Gamma}], (Y, \psi))$, let $\zeta(f') = f'_\un \lambda^{-1}$. One clearly has $\zeta(\xi(f)) = f$. Moreover,
  \begin{align*}
   \xi(\zeta(f'))_g &= \psi_g \circ (\g \tens (f'_\un \lambda^{-1})) = \psi_g \circ (\g \tens f'_\un) \circ (\g \tens \lambda^{-1}) \\ &= f'_g \circ \chi_g \restr{\un \tens X} \circ (\g \tens \lambda^{-1}) = f'_g \circ (\g \tens \lambda) \circ (\g \tens \lambda^{-1}) = f'_g
  \end{align*}
  thus $\xi$ and $\zeta$ are reciprocal morphisms. \cqfd
 \end{enumerate}
\end{demo}

\begin{cor}
 \label{strhom}
 Let $(X, \psi) \in \Cr {\Gamma}$. Then
 \begin{enumerate}
  \item The object $(X, \psi)$ is a direct summand of $X[{\Gamma}]$.
  \item If $(X, \psi)$ is indecomposable, ${\Gamma}$ acts transitively on the set of isoclasses of indecomposable summands of $X$.
 \end{enumerate}
\end{cor}

\begin{demo}
 \begin{enumerate}
  \item By proposition \ref{adjfc}, $X[{\Gamma}] \simeq k[{\Gamma}] \tens (X, \psi)$. This gives the result because $\un$ is a direct summand of $k[{\Gamma}]$. 
  \item Let $X = \bigoplus_{i=1}^\ell X_i$ in $\Cr$, the $X_i$ being indecomposable. Then $X[{\Gamma}] = \bigoplus_{i=1}^\ell X_i[{\Gamma}]$ and, as $\Cr {\Gamma}$ is Krull-Schmidt, $(X, \psi)$ is a direct summand of one of the $X_i[{\Gamma}]$. As a consequence, $\add(X) \subset \add F\left(X_i[{\Gamma}]\right) = \add \left(\bigoplus_{g \in {\Gamma}} (\g \tens X_i)\right)$. As $X_i$ is indecomposable and $\Cr$ is Krull-Schmidt, ${\Gamma}$ acts clearly transitively on this subcategory. \cqfd
 \end{enumerate}
\end{demo}

\begin{nt}
 Let $M \in \Cr {\Gamma}$. The number of indecomposable direct summands of $FM$ will be denoted by $\ell(M)$. The number of non isomorphic indecomposable direct summands of $FM$ will be denoted by $\# M$.
\end{nt}

\begin{rem}
 If $M$ is indecomposable and $\Gamma$ is cyclic, then $\# M = \ell(M)$. However, it is not the case in general, even if the group is commutative (see proposition \ref{relcycl} and the remark after it).
\end{rem}

\begin{lem}
 Let $X \in \Cr {\Gamma}$ be indecomposable. Then, $\# X$ divides $\ell(X)$. Their ratio is the number of copies of each indecomposable of $\add(FX)$ in $FX$.
\end{lem}

\begin{demo}
 Let $X_1, X_2$ be two indecomposable summands of $FX$. By corollary \ref{strhom}, as $X$ is indecomposable, there exists $g \in {\Gamma}$ such that $\g \tens X_1 \simeq X_2$. Moreover, $\g \tens X \simeq X$, therefore, the numbers of copies of $X_1$ and $X_2$ in $X$ are equal. \cqfd
\end{demo}

\begin{lem}
 If $\Cr$ is semisimple then $\Cr {\Gamma}$ is semisimple.
\end{lem}

\begin{demo}
 Let $X, Y$ two indecomposable objects of $\Cr {\Gamma}$ and $f \in \Hom_{\Cr {\Gamma}}(X, Y)$. As $\Cr$ is semisimple, $\Cr$ and $\Cr {\Gamma}$ are  abelian. Let $k: K \inj X$ the kernel of $f$. As $\Cr$ is semisimple, $Fk$ splits through an $\tilde h: FX \surj FK$. One gets $\tilde h \in \uHom_{\Cr {\Gamma}}(X, K)$. Let
 $$h = \frac{1}{\# {\Gamma}} \sum_{g \in {\Gamma}} g \cdot \tilde h$$
 with the result that $h$ is a morphism from $X$ to $K$ and that
 $$h k = \frac{1}{\# {\Gamma}} \sum_{g \in {\Gamma}} (g \cdot \tilde h) k = \frac{1}{\# {\Gamma}} \sum_{g \in {\Gamma}} (g \cdot \tilde h)(g \cdot k) = \frac{1}{\# {\Gamma}} \sum_{g \in {\Gamma}} g \cdot (\tilde h k) = \frac{1}{\# {\Gamma}}\sum_{g \in {\Gamma}} g \cdot \id_K = \id_K$$
 and finally, $k$ splits. As $\Cr {\Gamma}$ is Krull-Schmidt, $K$ is a direct summand of $X$. Finally, as $X$ is indecomposable, $\ker f = 0$ or $\ker f = X$. In the same way, $\coker f = 0$ or $\coker f = Y$ and, as a consequence, $f = 0$ or $f$ is invertible which shows that $\Cr {\Gamma}$ is semisimple. \cqfd
\end{demo}

\begin{df}
 \label{fonccomb}
 One denotes by $[\Cr]$ the semisimple $k$-category, whose simple objects are the isomorphism classes of indecomposable objects of $\Cr$. For $X = X_1 \oplus X_2 \oplus \dots \oplus X_n \in \Cr$ where $X_1, X_2, \dots, X_n$ are indecomposable, write
 $$[X] = [X_1] \oplus [X_2] \oplus \dots \oplus [X_n]$$
 where $[X_i]$ is the isomorphism class of $X_i$ (it is well defined because $\Cr$ is Krull-Schmidt). Moreover, if $X \in \Cr$ is indecomposable, one defines
 $$\End_{[\Cr]}([X]) = \End_\Cr(X) / \mathfrak{m}$$
 where $\mathfrak{m}$ is the maximal ideal of $\End_\Cr(X)$.
\end{df}

\begin{lem}
 Let $X_1, X_2, \dots, X_n \in \Cr$ be non isomorphic indecomposable objects. Then, for $i_1, i_2, \dots, i_n, j_1, j_2, \dots, j_n \in \N$,
 $$\dim_k(\Hom_{[\Cr]}([X_1^{i_1} \oplus X_2^{i_2} \oplus \dots \oplus X_n^{i_n}], [X_1^{j_1} \oplus X_2^{j_2} \oplus \dots \oplus X_n^{j_n}])) = \sum_{\ell = 1}^n i_\ell j_\ell c_\ell.$$
 Here, for every $\ell$, $c_\ell$ is the degree of the extension
 $$k \subset \End_\Cr(X_\ell) / \mathfrak{m}_\ell$$
 and $\mathfrak{m}_\ell$ denotes the maximal ideal of $\End_\Cr(X_\ell)$.
\end{lem}

\begin{demo}
 It is obvious. \cqfd
\end{demo}

\begin{lem}
 \label{lemcombeq}
 \begin{enumerate}
  \item The action of ${\Gamma}$ on $\Cr$ induces an action of ${\Gamma}$ on $[\Cr]$;
  \item if $k$ is algebraically closed and $[\Cr]$ has only one simple, then $[\Cr] {\Gamma} \simeq [\Cr {\Gamma}]$ as a $\md k[{\Gamma}]$-module category. 
 \end{enumerate}
\end{lem}

\begin{demo}
 First of all, up to an equivalence of categories, each isomorphism class of $\Cr$ can be supposed to contain exactly one object.
 \begin{enumerate}
  \item If $g \in {\Gamma}$ and $X \in \Cr$, let $\g \tens [X] = [\g \tens X]$. If $f: [X] \rightarrow [X]$ is a morphism, then $f$ comes from a morphism $f_0: X \rightarrow X$. Let $\g \tens f$ be the class of $\g \tens f_0$ modulo $\mathfrak m$. It is well defined as if $f_0$ is nilpotent, then $\g \tens f_0$ is also nilpotent. In the same way, one defines structural morphisms of $[\Cr]$ by projecting those of $\Cr$ modulo maximal ideals. It is now clear that $[\Cr]$ is a $\G$-module category.
  \item Let $X_0 \in \Cr$ be the only indecomposable object up to isomorphism. As $[\Cr] {\Gamma}$ and $[\Cr {\Gamma}]$ are semisimple and $k$ is algebraically closed, it is enough to see that the simple objects of both categories are in bijection. If $[(X, \psi)] \in [\Cr {\Gamma}]$, one can associate to it $([X], \psi') \in [\Cr] {\Gamma}$ by reducing matrix coefficients of $\psi$ modulo the maximal ideal of $\End_\Cr(X_0)$. Conversely, if $([X], \psi') \in [\Cr] {\Gamma}$, one car associate to it $[(X, \psi)]$ where the matrix coefficients of $\psi$ are those of $\psi'$ multiplied by $\id_{X_0}$. It is clearly a bijection. \cqfd
 \end{enumerate}
\end{demo}

\begin{lem}
 \label{nbiso}
 Suppose that $k$ is algebraically closed. Let $X \in \Cr$ be an indecomposable object and $Y \in \add(X[{\Gamma}])$ an indecomposable object of $\Cr \Gamma$. Then $X[{\Gamma}]$ has $\ell(Y)/\#(Y)$ indecomposable direct summands isomorphic to $Y$.
\end{lem}

\begin{demo}
 It will be proved in three steps. 
 \begin{enumerate}
  \item \label{casa} $\# Y = 1$. Then, up to restriction to $\add(FY) = \add(X)$, one can use lemma \ref{lemcombeq} and
  \begin{align*}
   \dim \Hom_{[\Cr {\Gamma}]}([Y], [k[{\Gamma}] \tens Y]) &= \dim \Hom_{[\Cr] {\Gamma}}([Y], k[{\Gamma}] \tens [Y]) \\ &= \dim \Hom_{[\Cr]}(F[Y], F[Y]) = \ell(Y)^2.
  \end{align*}
  Moreover, $X[{\Gamma}]^{\ell(Y)} = k[{\Gamma}] \tens Y$ which implies the result in this case.
  \item $\# Y = \# {\Gamma}$. In this case, $X[{\Gamma}]$ is indecomposable and therefore $Y \simeq X[{\Gamma}]$ which clearly implies the result.
  \item General case. Let $\tilde Y$ be the set of indecomposable summands of $FY$ up to isomorphism. The action of ${\Gamma}$ induces a morphism ${\Gamma} \rightarrow \mathfrak{S}_{\tilde Y}$, whose kernel $H$ is normal in ${\Gamma}$. So, one can apply proposition \ref{actssgrp}. As $\Cr {\Gamma} \simeq \Cr H ({\Gamma}/H)$, one gets a partial forgetful functor from $\Cr {\Gamma}$ to $\Cr H$. Let $\tilde Y \in \Cr H$ be the image of $Y$ by this functor. Let $Y'$ be an indecomposable direct summand of $\tilde Y$. As $\# Y' = 1$, case (\ref{casa}) applies to $Y'$. Moreover, $X[{\Gamma}] \simeq X[H][{\Gamma}/H]$ and $\# Y = \# {\Gamma}/H$. Hence, $Y \simeq Y'[{\Gamma}/H]$ appears in $X[{\Gamma}]$ the number of times $Y'$ appears in $X[H]$, that is $\ell(Y') = \ell(Y)/\#({\Gamma}/H) = \ell(Y)/\# Y$. \cqfd
 \end{enumerate}
\end{demo}

\subsection{Approximations}

\begin{df}
 One defines $\Add(\Cr)$ to be the class of all full sub-$k$-categories of $\Cr$ which are stable under isomorphisms and direct summands. If $E$ is a collection of objects of $\Cr$, one denotes by $\add(E)$ the smallest category of $\Add(\Cr)$ containing $E$. The category $\Tr \in \Add(\Cr)$ is said to be finitely generated if $\Tr$ is of the form $\add(M)$ for some object $M \in \Cr$. The subclass of $\Add(\Cr)$ consisting in all finitely generated categories will be denoted by $\adds(\Cr)$. If $F$ is a functor from $\Cr$ to a $k$-category $\Cr'$ and $\Tr \in \Add(\Cr)$, $F(\Tr)$ will denote $\add(\{F(X) \,|\, X \in \Tr\})$. If $\Cr$ is an $\Mr$-module category for some monoidal category $\Mr$, $\Add(\Cr)^\Mr$ (resp. $\adds(\Cr)^\Mr$) will denote the class of elements of $\Add(\Cr)$ (resp. $\adds(\Cr)$) which are sub-$\Mr$-module categories of $\Cr$. 
\end{df}

\begin{df}
 Let $\Tr \in \Add(\Cr)$ and $M \in \Cr$. A \emph{left (resp. right) $\Tr$-approximation of $M$} is an object $N \in \Tr$ and a morphism $f \in \Hom_\Cr (M, N)$ (resp. $\in \Hom_\Cr (N, M)$) such that for every $N' \in \Tr$ and $f': M \rightarrow N'$ (resp. $f': N' \rightarrow M$), $f'$ factors through $f$.
\end{df}

\begin{df}
 A morphism $f$ of the category $\Cr$ is said to be \emph{left (resp. right) minimal} if every morphism $g$ such that $f = g \circ f$ (resp. $f = f \circ g$) is an isomorphism.
\end{df}

The following lemmas are folklore. For detailed proofs, the reader is referred to \cite{AuReSm95} or \cite{De08-1}.

\begin{lem}
 \label{sdminimal}
 Let $X, Y, X', Y' \in \Cr$ and $f \in \Hom_\Cr (X, Y)$, $f' \in \Hom_\Cr (X', Y')$. 
 \begin{enumerate}
  \item The morphism $f$ is right minimal if and only if there is no decomposition  $X \simeq X_0 \oplus X_1$ such that $f \restr{X_0} = 0$ and $X_0 \neq 0$.
  \item The morphism $f$ is left minimal if and only if there is no decomposition $Y = Y_0 \oplus Y_1$ such that the corestriction of $f$ to $Y_0$ vanishes and $Y_0 \neq 0$.
  \item The morphisms $f$ and $f'$ are both left (resp. right) minimal if and only if $f \oplus f'$ is left (resp. right) minimal.
 \end{enumerate}
\end{lem}

\begin{lem}
 \label{sdapprox}
 Let $M, N, M', N' \in \Cr$, $\Tr \in \Add(\Cr)$, $f: M \rightarrow N$ and $f': M' \rightarrow N'$. Then $f$ and $f'$ are both left (resp. right) $\Tr$-approximation if and only if $f \oplus f'$ is a left (resp. right) $\Tr$-approximation.
\end{lem}

\begin{lem}
 If $\Tr \in \adds(\Cr)$ and $M \in \Cr$, then there exists a minimal left (resp. right) $\Tr$-approximation of $M$ which is unique up to (non unique) isomorphism. Moreover, any left (resp. right) $\Tr$-approximation of $M$ is the direct sum of the minimal left (resp. right) $\Tr$-approximation of $M$ and a morphism of the form $0 \rightarrow N$ (resp. $N \rightarrow 0$).
\end{lem}

\begin{lem}
 \label{admapp}
 Let $M \in \Cr$ and $\Tr \in \Add(\Cr)$. If $\Tr$ contains the injective envelope (resp. projective cover) of $M$, then every left (resp. right) $\Tr$-approximation of $M$ is an admissible monomorphism (resp. epimorphism).
\end{lem}

\begin{lem}
 \label{Fminapprox}
 Let $\T \in \Add(\Cr {\Gamma})^{\md k[{\Gamma}]}$, $X, Y \in \Cr {\Gamma}$ and $f \in \Hom_{\Cr {\Gamma}}(X, Y)$.
 \begin{enumerate}
  \item $f$ is a left (resp. right) $\T$-approximation if and only if $Ff$ is a left (resp. right) $F\T$-approximation.
  \item $f$ is a minimal left (resp. right) $\T$-approximation if and only if $Ff$ is a minimal left (resp. right) $F\T$-approximation.
 \end{enumerate}
\end{lem}

\begin{demo}
 By duality, it is enough to prove the statement for left approximations.
 \begin{enumerate}
  \item First of all, as $\T$ is $\md k[{\Gamma}]$-stable, if $FY \in F\T$, then $Y \in \T$ (and, by definition, if $Y \in \T$, $FY \in F \T$). As $F$ and $-[{\Gamma}]$ are adjoint, the following diagram is commutative for every $T \in \T$:
  $$\xymatrix{
   \Hom_\Cr(FY, FT) \ar[rrr]^{\Hom_\Cr(Ff, FT)} \ar[d]_\sim & & & \Hom_\Cr(FX, FT) \ar[d]_\sim \\
   \Hom_{\Cr {\Gamma}}(Y, (FT)[{\Gamma}]) \ar[rrr]^{\Hom_{\Cr {\Gamma}}(f, (FT)[{\Gamma}])} & & & \Hom_{\Cr {\Gamma}}(X, (FT)[{\Gamma}]). \\
  }$$
  The first line is surjective for every $T$ if and only if $Ff$ is a left $F\T$-approximation. As $\T$ is $\md k[{\Gamma}]$-stable, $(FT)[{\Gamma}] = k[{\Gamma}] \tens T \in \T$ and therefore the second line is surjective for every $T$ if and only if $f$ is a left $\T$-approximation (because $T$ is a direct summand of $(FT)[{\Gamma}]$).
  \item If $Ff$ is left minimal then $f$ is clearly left minimal. Conversely, suppose that $f$ is a minimal left $\T$-approximation. Let $\tilde f: FX \rightarrow \tilde Y$ be a minimal left $F\T$-approximation of $FX$. Let $f' = \tilde f[{\Gamma}]: (FX)[{\Gamma}] \rightarrow \tilde Y[{\Gamma}]$. Then, one gets
  $$Ff' = \bigoplus_{g \in {\Gamma}} \g \tens \tilde f.$$
  For every $g \in {\Gamma}$, $\g \tens \tilde f$ is a minimal left $\T$-approximation because $F\T$ is ${\Gamma}$-stable, and, as a consequence, by using lemmas \ref{sdapprox} and \ref{sdminimal}, $Ff'$ is also a minimal left $F\T$-approximation. Then, $f'$ is a minimal left $\T$-approximation. By unicity of a minimal left $\T$-approximation, and using lemmas \ref{sdapprox} and \ref{sdminimal}, $f' \simeq f \oplus g$ where $g$ is a minimal left $\T$-approximation of $\bigoplus_{g \in {\Gamma} \setminus \{e\}} \g \tens X$. Finally, as $Ff' = Ff \oplus Fg$ is left minimal, $Ff$ is also left minimal by lemma \ref{sdminimal}. \cqfd
 \end{enumerate}
\end{demo}

\subsection{Action on an exact category}

The action of ${\Gamma}$ on $\Cr$ is now supposed to be exact. Hence, $\Cr {\Gamma}$ is exact. It is easy to see that the functor $-[{\Gamma}]$ from $\Cr$ to $\Cr {\Gamma}$ is exact. 

\begin{lem}
 If $X \in \Cr$ is injective (resp. projective), then $X[{\Gamma}]$ is injective (resp. projective).
\end{lem}

\begin{demo}
 Suppose that $X$ is injective. Let $0 \rightarrow X[{\Gamma}] \xrightarrow{f} Y \xrightarrow{g} Z \rightarrow 0$ be an admissible short exact sequence in $\Cr {\Gamma}$. By definition, $0 \rightarrow F(X[{\Gamma}]) \xrightarrow{f} FY \xrightarrow{g} FZ \rightarrow 0$ is an admissible short exact sequence in $\Cr$. Applying $\Hom_{\Cr {\Gamma}}(-,X[{\Gamma}])$ and $\Hom_\Cr(-,X)$ gives the long exact sequences
 $$0 \rightarrow \Hom_{\Cr {\Gamma}}(Z, X[{\Gamma}]) \rightarrow \Hom_{\Cr {\Gamma}}(Y, X[{\Gamma}]) \rightarrow \Hom_{\Cr {\Gamma}}(X[{\Gamma}], X[{\Gamma}]) \rightarrow \dots$$
 $$0 \rightarrow \Hom_\Cr (FZ, X) \rightarrow \Hom_\Cr (FY, X) \rightarrow \Hom_\Cr (F(X[{\Gamma}]), X) \rightarrow \Ext_\Cr^1 (FZ, X) = 0.$$

 The isomorphism of bifunctors $\Hom_{\Cr {\Gamma}}(-, -[{\Gamma}]) \simeq \Hom_\Cr(F-, -)$ permits to conclude that $\Hom_{\Cr {\Gamma}}(f, X[{\Gamma}]): \Hom_{\Cr {\Gamma}}(Y, X[{\Gamma}]) \rightarrow \Hom_{\Cr {\Gamma}}(X[{\Gamma}], X[{\Gamma}])$ is surjective. Therefore the admissible short exact sequence $0 \rightarrow X[{\Gamma}] \xrightarrow{f} Y \xrightarrow{g} Z \rightarrow 0$ splits. Hence, $X[{\Gamma}]$ is injective. The proof is similar for the projective case. \cqfd
\end{demo}

\begin{cor}
 \label{injcr}
 If $\Cr$ has enough injective (resp. projective) objects, then $\Cr {\Gamma}$ has also enough injective (resp. projective) objects. Moreover, for $X \in \Cr {\Gamma}$, there exists an injective resolution $I_\bullet$ (resp. a projective resolution $P_\bullet$) of $X$ such that $FI_\bullet$ (resp. $FP_\bullet$) is an injective resolution (resp. a projective resolution) of $FX$.
\end{cor}

\begin{demo}
 Suppose that $\Cr$ has enough injective objects. Let $X \in \Cr {\Gamma}$. There exists an admissible short exact sequence $0 \rightarrow FX \rightarrow I \rightarrow Y \rightarrow 0$ in $\Cr$ where $I$ is injective. As the action of ${\Gamma}$ is exact, it gives an admissible short exact sequence $0 \rightarrow (FX)[{\Gamma}] \rightarrow I[{\Gamma}] \rightarrow Y[{\Gamma}] \rightarrow 0$. Moreover, $X$ is a direct summand of $(FX)[{\Gamma}]$ and, as the composition of two admissible monomorphisms is an admissible monomorphism, one gets an admissible monomorphism $X \rightarrow I[{\Gamma}]$. For the second part, it is enough to apply inductively the first part because $I[{\Gamma}]$ and $F(I[{\Gamma}])$ are both injective in $\Cr {\Gamma}$ and $\Cr$, respectively. \cqfd
\end{demo}

One supposes now that $\Cr$ has enough injectives or enough projectives.

\begin{df}
 For $X \in \Cr {\Gamma}$ and $n \in \N$, $\uExt^n_{\Cr {\Gamma}}(X, -)$ will denote the right derived functor of $\uHom_{\Cr {\Gamma}}(X, -)$ if $\Cr$ has enough injective objects and $\uExt^n_{\Cr {\Gamma}}(-, X)$ will denote the left derived functor of $\uHom_{\Cr {\Gamma}}(-, X)$ if $\Cr$ has enough projective objects. If $\Cr$ has both enough injective and projective objects, the two definitions coincide as usual.
\end{df}

Hence, one gets on $\Cr {\Gamma}$ a structure of $\md k[{\Gamma}]$-exact category. All usual homological results remain true in this context. The following proposition links these properties with the usual $k$-exact structure and summarizes some properties of $\uExt$.

\begin{prop}
 \label{repext}
 For every $n \in \N$, there are functorial isomorphisms (the $?$ are variables in $\md k[{\Gamma}]$ and the $-$ are variables in $\Cr {\Gamma}$):
 \begin{enumerate}
  \item $\uExt^n_{\Cr {\Gamma}}(?_1 \tens -_1,?_2 \tens -_2) \simeq?_1^* \tens?_2 \tens \uExt^n_{\Cr {\Gamma}}(-_1, -_2)$;
  \item $\Ext^n_{\Cr {\Gamma}}(?_1 \tens -_1,?_2 \tens -_2) \simeq \Hom_{\md k[{\Gamma}]}(?_1 \tens?_2^*, \uExt^n_{\Cr {\Gamma}}(-_1, -_2))$;
  \item $F \uExt^n_{\Cr {\Gamma}} (-, -) \simeq \Ext^n_\Cr (F-, F-)$;
  \item $\Ext^n_{\Cr {\Gamma}}(?^* \tens -, -) \simeq \Ext^n_{\Cr {\Gamma}}(-,? \tens -)$;
  \item $\Ext^n_{\Cr {\Gamma}}(?[{\Gamma}], -) \simeq \Ext^n_\Cr(?, F-)$;
  \item $\Ext^n_{\Cr {\Gamma}}(-,?[{\Gamma}]) \simeq \Ext^n_\Cr(F-,?)$ 
 \end{enumerate}
 where $F$ denotes the forgetful functors from $\Cr {\Gamma}$ to $\Cr$ and from $\md k[{\Gamma}]$ to $\vect k$. Moreover, these isomorphisms commute with long exact sequences obtained from admissible short exact sequences in $\Cr {\Gamma}$.
\end{prop}

\begin{demo}
 These are easy consequences of lemma \ref{injcr} together with proposition \ref{rephom} and standard homological constructions. \cqfd
\end{demo}

\begin{cor}
 \label{carinjproj}
 If $X \in \Cr {\Gamma}$, $X$ is injective (resp. projective) if and only if $FX$ is injective (resp. projective).
\end{cor}

\begin{demo}
 If $X$ is injective, for every $Y \in \Cr$, 
 $$\Ext^n_\Cr(Y, FX) \simeq \Ext^n_{\Cr {\Gamma}}(Y[{\Gamma}], X) = 0$$
 so that $FX$ is injective. If $FX$ is injective, for all $Y \in \Cr {\Gamma}$, as $Y$ is a direct summand of $(FY)[{\Gamma}]$,
 $$\Ext^n_{\Cr {\Gamma}}(Y, X) \subset \Ext^n_{\Cr {\Gamma}}((FY)[{\Gamma}], X) \simeq \Ext^n_\Cr(FY, FX) = 0$$
 hence $X$ is injective. The proof is the same for the projective case. \cqfd
\end{demo}

\begin{lem}
 For every representation $r \in \md k[{\Gamma}]$, the functor $r \tens -$ from $\Cr {\Gamma}$ to $\Cr {\Gamma}$ is exact.
\end{lem}

\begin{demo}
 Let $0 \rightarrow X \rightarrow Y \rightarrow Z \rightarrow 0$ be an admissible short exact sequence of $\Cr {\Gamma}$. By definition, $0 \rightarrow FX \rightarrow FY \rightarrow FZ \rightarrow 0$ is an admissible short exact sequence of $\Cr$. As $0 \rightarrow F(r \tens X) \rightarrow F(r \tens Y) \rightarrow F(r \tens Z) \rightarrow 0$ is isomorphic to $0 \rightarrow (FX)^{\dim r} \rightarrow (FY)^{\dim r} \rightarrow (FZ)^{\dim r} \rightarrow 0$, this is an admissible short exact sequence of $\Cr$ and therefore, by definition, $0 \rightarrow r \tens X \rightarrow r \tens Y \rightarrow r \tens Z \rightarrow 0$ is an admissible short exact sequence of $\Cr {\Gamma}$. \cqfd
\end{demo}

\subsection{Action on a Frobenius stably $2$-Calabi-Yau category}

As before, the group ${\Gamma}$ is supposed to be finite, of cardinality non divisible by the characteristic of $k$ and the $k$-category $\Cr$ is exact, $\Hom$-finite, Krull-Schmidt. Recall that $\Cr$ is called \emph{Frobenius} if it has enough projectives and enough injectives and if the projective objects and the injective objects are the same. Recall that $\Cr$ is called \emph{(stably) $2$-Calabi-Yau} if there is a functorial isomorphism $c: \Ext^1_\Cr(-_1, -_2) \simeq \Ext^1_\Cr(-_2, -_1)^*$. In the following, $\Cr$ will be supposed to be Frobenius and stably $2$-Calabi-Yau. The category $\Cr {\Gamma}$ is Frobenius by corollary \ref{carinjproj}. One fixes an isomorphism of bifunctors $c: \Ext^1_\Cr(-_1, -_2) \simeq \Ext^1_\Cr(-_2, -_1)^*$.

\begin{df}
 \label{act2CY}
 The action of ${\Gamma}$ on $\Cr$ is said to be \emph{$2$-Calabi-Yau} (for $c$) if it is exact and for every $g \in {\Gamma}$, the following diagram commutes:
 $$\xymatrix{
  \Ext^1_\Cr(-_1, -_2) \ar[d]^c \ar[rr]^{\g \tens -} & & \Ext^1_\Cr(\g \tens -_1, \g \tens -_2) \ar[d]^c \\
  \Ext^1_\Cr(-_2, -_1)^* & & \Ext^1_\Cr(\g \tens -_2, \g \tens -_1)^* \ar[ll]^{(\g \tens -)^*}
 }$$
\end{df}

From now on, the action of ${\Gamma}$ on $\Cr$ is assumed to be $2$-Calabi-Yau (for $c$).

\begin{prop}
 \label{CG2CY}
 \begin{enumerate}
  \item \label{CG2CY1} The functorial isomorphism of vector spaces $c$ is also a functorial isomorphism of $k[{\Gamma}]$-modules:
   $$\uExt_{\Cr {\Gamma}}^1(-_1, -_2) \simeq \uExt^1_{\Cr {\Gamma}}(-_2, -_1)^*$$
   (recall that for any $X, Y \in \Cr {\Gamma}$, the underlying vector space of the $\md k[{\Gamma}]$-module $\uExt_{\Cr {\Gamma}}^1(X, Y)$ is $\Ext_{\Cr}^1(FX, FY)$).
  \item The category $\Cr {\Gamma}$ is $2$-Calabi-Yau.
 \end{enumerate}
\end{prop}

\begin{demo}
 Recall that there is an isomorphism of functors from $\md k[{\Gamma}]$ into itself:
  $$\Hom_{\md k[{\Gamma}]}(\un, -) \simeq \Hom_{\md k[{\Gamma}]}(-, \un)^*.$$
 \begin{enumerate}
  \item The only thing to prove is that, for any $(X, \psi), (Y, \chi) \in \Cr {\Gamma}$, $c_{X, Y}$ is in fact a morphism of representations of ${\Gamma}$ from $\uExt_{\Cr {\Gamma}}^1((X, \psi), (Y, \chi))$ to $\uExt^1_{\Cr {\Gamma}}((Y, \chi), (X, \psi))^*$. For $g \in {\Gamma}$, it is enough to show that the following diagram commutes:
  $$\xymatrix{
   \Ext_\Cr^1(X, Y) \ar[rr]^{c_{X, Y}} \ar[d]_{\g \tens -} & & \Ext_\Cr^1(Y, X)^* \ar[d]^{\left((\g \tens -)^*\right)^{-1}} \\
   \Ext_\Cr^1(\g \tens X, \g \tens Y) \ar[rr]^{c_{\g \tens X, \g \tens Y}} \ar[d]_{\Ext^1_\Cr(\psi_g^{-1}, \chi_g)} & & \Ext_\Cr^1(\g \tens Y, \g \tens X)^* \ar[d]|{\vphantom{()}=}_{\Ext^1_\Cr(\chi_g, \psi_g^{-1})^*}^{\left(\Ext^1_\Cr(\chi_g^{-1}, \psi_g)^*\right)^{-1}} \\
   \Ext_\Cr^1(X, Y) \ar[rr]^{c_{X, Y}} & & \Ext_\Cr^1(Y, X)^* 
  }$$
  (the left side comes from the action of $g$ on $\Ext_\Cr^1(X, Y)$ and the right side comes from the inverse of the adjoint of the action of $g$); the upper square commutes because the action of ${\Gamma}$ is $2$-Calabi-Yau and the lower square commutes because the isomorphism $c$ is functorial. 
  \item Denote by $\mathbf{c}$ the isomorphism of functors of (\ref{CG2CY1}). Let
   $$\tilde c = \Hom_{\md k[{\Gamma}]}(\un, \mathbf{c}).$$
   Then $\tilde c$ is an isomorphism of functors 
    $$\Hom_{\md k[{\Gamma}]}(\un, \uExt^1_{\Cr {\Gamma}}(-_1, -_2)) \longrightarrow \Hom_{\md k[{\Gamma}]}(\un, \uExt^1_{\Cr {\Gamma}}(-_2, -_1)^*).$$
   The reciprocal isomorphism is
   $$\tilde c^{-1} = \Hom_{\md k[{\Gamma}]}(\un, \mathbf{c}^{-1}).$$
   Moreover, lemma \ref{repext} leads to
   $$\Hom_{\md k[{\Gamma}]}(\un, \uExt^1_{\Cr {\Gamma}}(-_1, -_2)) \simeq \Ext^1_{\Cr {\Gamma}}(-_1, -_2)$$
   and 
   \begin{align*}
    \Hom_{\md k[{\Gamma}]}(\un, \uExt^1_{\Cr {\Gamma}}(-_2, -_1)^*) &\simeq \Hom_{\md k[{\Gamma}]}(\uExt^1_{\Cr {\Gamma}}(-_2, -_1), \un) \\
    &\simeq \Hom_{\md k[{\Gamma}]}(\un, \uExt^1_{\Cr {\Gamma}}(-_2, -_1))^* \\
    &\simeq \Ext^1_{\Cr {\Gamma}}(-_2, -_1)^*
   \end{align*}
  which finishes the proof. \cqfd  
 \end{enumerate}
\end{demo}

\subsection{Computation of $(kQ){\Gamma}$ and $(\Lambda_Q){\Gamma}$}

\label{acc}
The aim of this section is to summarize some properties proved in \cite{De-1} useful to compute equivariant categories for categories of modules. The assumptions on $k$ and ${\Gamma}$ are the same as before. If $\Lambda$ is a $k$-algebra and if ${\Gamma}$ acts on $\Lambda$, the action being denoted exponentially, the skew group algebra of $\Lambda$ under the action of ${\Gamma}$ is by definition the $k$-algebra whose underlying $k$-vector space is $k[{\Gamma}] \tens_k \Lambda$ and whose multiplication is linearly generated by $(g \tens a) (g' \tens a') = gg' \tens a^{g'^{-1}} a'$ for all $g, g' \in {\Gamma}$ and $a, a' \in \Lambda$ (see \cite{ReRi85}). It will be denoted by $\Lambda {\Gamma}$. Identifying $k[{\Gamma}]$ and $\Lambda$ with subalgebras of $\Lambda {\Gamma}$, an alternative definition is
$$\Lambda {\Gamma} = \langle \Lambda, k[{\Gamma}] \,|\, \forall (g, a) \in {\Gamma} \times \Lambda, g a g^{-1} = a^g \rangle_{k\text{-alg}}$$

The following links skew group algebras with equivariant categories.

\begin{prop}
 \label{eqskeweq}
 The action of ${\Gamma}$ on $\Lambda$ induces an action of ${\Gamma}$ on $\md \Lambda$. Moreover, there is a canonical equivalence of categories between  $\md (\Lambda {\Gamma})$ and $\md (\Lambda) {\Gamma}$.
\end{prop}

\begin{demo}
 If $g \in {\Gamma}$ and $(V, r) \in \md \Lambda$, one denotes by $\g \tens (V, r)$ the representation $(V, \g \tens r)$ of $\Lambda$ where, if $a \in \Lambda$, $(\g \tens r)_a = r_{g^{-1} a}$. If $f \in \Hom_{\md \Lambda}((V, r), (V', r'))$, one defines $\id_\g \tens f = f$. Extending this definition by linearity on the whole category $\G$, $\md \Lambda$ is a $\G$-module category. The structural isomorphisms are identities. 
 
 Let $(V, r, \psi) \in \md (\Lambda) {\Gamma}$. For $g \in {\Gamma}$ and $a \in \Lambda$, one defines
 $$r^\psi_{g \tens a} = \psi_g r_a.$$
 Since for every $g, g' \in {\Gamma}$ and $a, a' \in \Lambda$ the diagram
 $$\xymatrix{
  V \ar[rrr]^{r_{a'}} \ar[drrr]_{r_{g'^{-1}(a) a'}} & & & V \ar[rrr]^{\psi_{g'}} \ar[d]|{r_{g'^{-1}(a)}=(\g' \tens r)_a} & & & V \ar[d]^{r_a} \\
  & & & V \ar[drrr]_{\psi_{gg'}} \ar[rrr]^{\psi_{g'}} & & & V \ar[d]^{\psi_g} \\
  & & & & & & V
 }$$
 is commutative, one gets that $r^\psi_{g \tens a} r^\psi_{g' \tens a'} = r^\psi_{(g \tens a)(g' \tens a')}$. Therefore $(V, r^\psi)$ is a representation of $\Lambda {\Gamma}$. 
 
 If $f \in \Hom_{\md(\Lambda) {\Gamma}}((V, r, \psi), (V', r', \psi'))$,  it is clear that 
 $$f \in \Hom_{\md(\Lambda {\Gamma})}((V, r^\psi), (V', r'^{\psi'})).$$ 
 Hence $(V, r, \psi) \mapsto (V, r^\psi)$ is a functor from $\md(\Lambda) {\Gamma}$ to $\md(\Lambda {\Gamma})$. 
 
 Let $(V, r^0) \in \md(\Lambda {\Gamma})$. If $a \in \Lambda$, let $r_a = r^0_{1 \tens a}$ which gives that $(V, r) \in \md(\Lambda)$. If $g \in {\Gamma}$, let $\psi_g = r^0_{g \tens 1}$. For $g \in {\Gamma}$ and $a \in \Lambda$, the following diagram commutes:
 $$\xymatrix{
  V \ar[rr]^{\psi_g} \ar[d]_{r_{g^{-1}(a)}} \ar[drr]^{r^0_{g \tens g^{-1}(a)}} & & V \ar[d]^{r_a} \\
  V \ar[rr]^{\psi_g} & & V.
 }$$
 Hence $\psi_g$ is an isomorphism from $\g \tens (V, r)$ to $(V, r)$. It is clear that $(V, r, \psi) \in \md(\Lambda){\Gamma}$. 
 
 If $f \in \Hom_{\md(\Lambda {\Gamma})}((V, r^0), (V', r'^0))$, one gets immediately
 $$f \in \Hom_{\md(\Lambda) {\Gamma}}((V, r, \psi), (V', r', \psi')).$$
 The two constructed functors are mutually inverse. \cqfd
\end{demo}

Let now $Q = (Q_0, Q_1)$ be a quiver. Consider an action of ${\Gamma}$ on the path algebra $kQ$ permuting the set of primitive idempotents $\{e_i \,|\, i \in Q_0\}$. We now define a new quiver $Q_{\Gamma}$. 

Let $\tilde Q_0$ be a set of representatives of the classes of $Q_0$ under the action of ${\Gamma}$. For $i \in Q_0$, let ${\Gamma}_i$ denote the subgroup of ${\Gamma}$ stabilizing $e_i$, let $i_\circ \in \tilde Q_0$ be the representative of the class of $i$ and let $\kappa_i \in {\Gamma}$ be such that $\kappa_i i_\circ = i$.

For $(i, j) \in \tilde Q_0^2$, ${\Gamma}$ acts on $O_i \times O_j$ where $O_i$ and $O_j$ are the orbits of $i$ and $j$ under the action of ${\Gamma}$. A set of representatives of the classes of this action will be denoted by $F_{ij}$. 

For $i, j \in Q_0$, define $A_{ij} = e_j \left(\rad(kQ)/\rad(kQ)^2\right) e_i$ where $\rad(kQ)$ is the Jacobson radical of $kQ$. We regard $A_{ij}$ as a left $k[{\Gamma}_i \inter {\Gamma}_j]$-module by restricting the action of ${\Gamma}$.

The quiver $Q_{\Gamma}$ has vertex set 
 $$Q_{\Gamma,0} = \bigcup_{i \in \tilde Q_0} \{i\} \times \irr({\Gamma}_i)$$
 where $\irr({\Gamma}_i)$ is a set of representatives of isomorphism classes of irreducible representations of ${\Gamma}_i$. The set of arrows of $Q_{\Gamma}$ from $(i, \rho)$ to $(j, \sigma)$ is a basis of
 $$\bigoplus_{(i', j') \in F_{ij}} \Hom_{\md k[{\Gamma}_{i'} \inter {\Gamma}_{j'}]}((\kappa_{i'} \cdot \rho) \restr{{\Gamma}_{i'} \inter {\Gamma}_{j'}} \tens A_{i'j'}, (\kappa_{j'} \cdot \sigma) \restr{{\Gamma}_{i'} \inter {\Gamma}_{j'}})^*$$
 where the representation $\kappa_{i'} \cdot \rho$ of ${\Gamma}_{i'}$ is the same as $\rho$ as a vector space, and $(\kappa_{i'} \cdot \rho)_{g} = \rho_{\kappa_{i'}^{-1} g \kappa_{i'}}$ for $g \in {\Gamma}_{i'} = \kappa_{i'} {\Gamma}_i \kappa_{i'}^{-1}$. 

\begin{thm}[\citeb{theorem 1}{De-1}] \label{th1}
 There is an equivalence of categories
 $$\md k\left(Q_{\Gamma}\right) \simeq \md \left(kQ\right){\Gamma}.$$
\end{thm}
 Theorem \ref{th1} was also proved by Reiten and Riedtmann in \cite[\S 2]{ReRi85} for cyclic groups. The following theorem deals with the case of preprojective algebras $\Lambda_Q$.

\begin{thm}[\citeb{theorem 2}{De-1}] \label{th2}
  If ${\Gamma}$ acts on $k \bar Q$, where $\bar Q$ is the double quiver of $Q$, by permuting the primitive idempotents $e_i$, and if for all $g \in {\Gamma}$, $r^g = r$ where $r$ is the preprojective relation of this quiver, then $\left(\bar Q\right)_{\Gamma}$ is of the form $\bar Q'$ for some quiver $Q'$ and $(\Lambda_Q){\Gamma}$ is Morita equivalent to $\Lambda_{Q'}$.
\end{thm}

One can always extend an action on $kQ$ to an action on $k \bar Q$ and this yields:
 
\begin{cor}[\citeb{corollary 1}{De-1}] \label{c1}
 An action of ${\Gamma}$ on a path algebra $k Q$ permuting the primitive idempotents induces naturally an action of ${\Gamma}$ on $k \bar Q$ and $\left(\bar Q\right)_{\Gamma}$ is isomorphic to the double quiver of $Q_{\Gamma}$. Moreover, there is an equivalence of categories
  $$\md \Lambda_{Q_{\Gamma}} \simeq \md \Lambda_Q {\Gamma}.$$
\end{cor}

\section{Categorification of skew-symmetrizable cluster algebras}

 In this part, $\Cr$ is supposed to be exact, Frobenius, $\Hom$-finite, stably $2$-Calabi-Yau and Krull-Schmidt. As before, ${\Gamma}$ is a finite group whose cardinality is not divisible by the characteristic of $k$. The group ${\Gamma}$ is supposed to act on $\Cr$, the action being exact and $2$-Calabi-Yau (see definition \ref{act2CY}).

The results of this section generalize works by Gei\ss, Leclerc and Schröer (in particular \cite{GeLeSc06-1} and \cite{GeLeSc}) in the context of preprojective algebras, and works of Dehy, Fu, Keller, Palu and others in the context of exact categories (see \cite{DeKe08}, \cite{FuKe}, \cite{Pa08}).

\subsection{Mutation of maximal ${\Gamma}$-stable rigid subcategories}

\begin{df}
 A category $\T \in \Add(\Cr)$ (resp. $\in \Add(\Cr {\Gamma})$) is said to be \emph{rigid} if there are no non trivial extensions between its objects. If moreover every rigid category $\T'$ containing $\T$ is equal to $\T$, then we say that $\T$ is \emph{maximal rigid}.
\end{df}

The following definition was introduced in \cite{Iy07}:
\begin{df}
 A category $\T \in \Add(\Cr)$ (resp. $\in \Add(\Cr {\Gamma})$) is said to be \emph{cluster-tilting} or \emph{maximal $1$-orthogonal} if for any $X \in \Cr$ (resp $\in \Cr {\Gamma}$), the following are equivalent:
 \begin{itemize}
  \item $\forall Y \in \T, \Ext^1(X, Y) = 0$;
  \item $X \in \T$.
 \end{itemize}
\end{df}

Clearly, any cluster-tilting category is maximal rigid.

\begin{df}
 If $\T \in \Add(\Cr)^{\Gamma}$ is rigid, $\T$ is said to be \emph{maximal ${\Gamma}$-stable rigid} if for every rigid $\T' \in \Add(\Cr)^{\Gamma}$  containing $\T$, $\T' = \T$. Similarly, if $\T \in \Add(\Cr {\Gamma})^{\md k[{\Gamma}]}$ is rigid, $\T$ is said to be \emph{maximal $\md k[{\Gamma}]$-stable rigid} if for every rigid $\T' \in \Add(\Cr {\Gamma})^{\md k[{\Gamma}]}$ containing $\T$, $\T' = \T$.
\end{df}

\begin{rem}
 Being maximal ${\Gamma}$-stable (resp. $\md k[{\Gamma}]$-stable) rigid is weaker than being maximal rigid and ${\Gamma}$-stable (resp. $\md k[{\Gamma}]$-stable). For instance, the quiver
 $$\xymatrix{
  1 \ar@/^/[r]^\alpha & 2 \ar@/^/[l]^{\alpha^*}
 }$$
 with relations $\alpha \alpha^* = \alpha^*\alpha = 0$, on which $\Z/2\Z$ acts by exchanging the two vertices, has no maximal rigid object which is $\Z/2\Z$-stable. The only maximal ${\Gamma}$-stable rigid object is the direct sum of the projective ones.
\end{rem}

\begin{prop}
 \label{bijamasbasc}
 The functors $-[{\Gamma}]$ and $F$ induce reciprocal bijections between:
 \begin{enumerate}
  \item $\Add(\Cr)^{\Gamma}$ and $\Add(\Cr {\Gamma})^{\md k[{\Gamma}]}$;
  \item the set of rigid $\T \in \Add(\Cr)^{\Gamma}$ and the set of rigid $\T \in \Add(\Cr {\Gamma})^{\md k[{\Gamma}]}$;
  \item the set of maximal ${\Gamma}$-stable rigid $\T \in \Add(\Cr)^{\Gamma}$ and the set of maximal $\md k[{\Gamma}]$-stable rigid $\T \in \Add(\Cr {\Gamma})^{\md k[{\Gamma}]}$;
  \item the set of cluster-tilting $\T \in \Add(\Cr)^{\Gamma}$ and the set of cluster-tilting $\T \in \Add(\Cr {\Gamma})^{\md k[{\Gamma}]}$.
 \end{enumerate}
 Moreover, all these bijections restrict to bijections between the corresponding finitely generated classes.
\end{prop}

\begin{demo}
 \begin{enumerate}
  \item If $\Dr \in \Add(\Cr)^{\Gamma}$, it is clear that $\Dr[{\Gamma}] \in \Add(\Cr {\Gamma})^{\md k[{\Gamma}]}$. Similarly, if $\Dr' \in \Add(\Cr {\Gamma})^{\md k[{\Gamma}]}$, it is easy to see that $F \Dr' \in \Add(\Cr)^{\Gamma}$. Suppose now that $X \in F(\Dr[{\Gamma}])$. This means that $X$ is a direct summand of $F(Y[{\Gamma}])$ for some $Y \in \Dr$. But $F(Y[{\Gamma}])$ is the direct sum of the $\g \tens Y$ for $g \in {\Gamma}$, hence $F(Y[{\Gamma}]) \in \Dr$ and finally, $X \in \Dr$. If $X \in \Dr$, $X$ is a direct summand of  $F(X[{\Gamma}])$ so that $X$ is in $F(\Dr[{\Gamma}])$. One concludes that $F(\Dr[{\Gamma}]) = \Dr$.   
   
   Suppose that $X \in F(\Dr')[{\Gamma}]$. It means that $X$ is a direct summand of some $F(Y)[{\Gamma}] \simeq k[{\Gamma}] \tens Y$ where $Y \in \Dr'$ and, as $\Dr'$ is $\md k[{\Gamma}]$-stable, $F(Y)[{\Gamma}]$ and $X$ are in $\Dr'$. On the other hand, if $X \in \Dr'$, as $X$ is a direct summand of $F(X)[{\Gamma}]$, $X$ is in $F(\Dr')[{\Gamma}]$. Finally $F(\Dr')[{\Gamma}] = \Dr'$. Therefore, $F$ and $-[{\Gamma}]$ induce reciprocal bijections.
  \item Suppose that $\T \in \Add(\Cr)^{\Gamma}$ is rigid. Let $X \in \T[{\Gamma}]$. By definition, there exists $Y \in \T$ and $X' \in \T[{\Gamma}]$ such that $Y[{\Gamma}] \simeq X \oplus X'$. Thus, one gets
  \begin{align*}
  \Ext^1_{\Cr {\Gamma}}(X, X) &\subset \Ext^1_{\Cr {\Gamma}}(Y[{\Gamma}], Y[{\Gamma}]) \simeq \Ext^1_\Cr(Y, F(Y[{\Gamma}])) \\ &\subset \Ext^1_\Cr(F(Y[{\Gamma}]), F(Y[{\Gamma}]))= 0
  \end{align*}
  because, as $Y \in \T$ and $\T$ is ${\Gamma}$-stable, $F(Y[{\Gamma}]) \in \T$. As a consequence, $X$ is rigid and therefore $\T[{\Gamma}]$ is rigid.
  
  Suppose now that $\T \in \Add(\Cr {\Gamma})^{\md k[{\Gamma}]}$ is rigid. Let $X \in F\T$. By definition, there exists $Y \in \T$ and $X' \in F\T$ such that $FY = X \oplus X'$. One gets
  \begin{align*}
  \Ext^1_\Cr(X, X) &\subset \Ext^1_\Cr(FY, FY) \simeq \Ext^1_{\Cr {\Gamma}}(Y, (FY)[{\Gamma}]) \\ &\subset \Ext^1_{\Cr {\Gamma}}((FY)[{\Gamma}], (FY)[{\Gamma}]) = 0
  \end{align*}
  because, as $Y \in \T$ and $\T$ is $\md k[{\Gamma}]$-stable, $(FY)[{\Gamma}] \simeq k[{\Gamma}] \tens Y \in \T$. As a consequence, $F\T$ is rigid.
  \item This is clear because the two bijections are increasing (with respect to inclusion).
  \item Suppose that $\T \in \Add(\Cr)^{\Gamma}$ is cluster-tilting. Let $X \in \Cr {\Gamma}$ such that for every $Y \in \T[{\Gamma}]$, $\Ext^1_{\Cr {\Gamma}}(X, Y) = 0$. In particular, for every $Z \in \T$, $\Ext^1_\Cr(FX, Z) \simeq \Ext^1_{\Cr {\Gamma}}(X, Z[{\Gamma}]) = 0$. Therefore, as $\T$ is cluster-tilting, $FX \in \T$, $(FX)[{\Gamma}] \in \T[{\Gamma}]$, and $X \in \T[{\Gamma}]$ because $X$ is a direct summand of $(FX)[{\Gamma}]$. Finally, as $\T[{\Gamma}]$ is rigid, $\T[{\Gamma}]$ is cluster-tilting.
  
  Suppose that $\T \in \Add(\Cr {\Gamma})^{\md k[{\Gamma}]}$ is cluster-tilting. Let $X \in \Cr$ be such that for every $Y \in F\T$, $\Ext^1_\Cr(X, Y) = 0$. In particular, for all $Z \in \T$, $\Ext^1_{\Cr {\Gamma}}(X[{\Gamma}], Z) \simeq \Ext^1_\Cr(X, FZ) = 0$. Therefore $X[{\Gamma}] \in \T$,  $F(X[{\Gamma}]) \in F\T$ and, $X \in F\T$ as $X$ is a direct summand of $F(X[{\Gamma}])$. Finally, as $F\T$ is rigid, $F\T$ is cluster-tilting. \cqfd  
 \end{enumerate}
\end{demo}

\begin{lem}
 \label{cokapp}
 Let $\T \in \Add(\Cr {\Gamma})^{\md k[{\Gamma}]}$ be rigid, and let $X \in \Cr {\Gamma}$ be such that $k[{\Gamma}] \tens X$ is rigid. If
 $$0 \rightarrow X \xrightarrow{f} T \xrightarrow{g} Y \rightarrow 0$$
 is an admissible short exact sequence and $f$ is a left $\T$-approximation, then the category $$\add(\T, k[{\Gamma}] \tens Y)$$ is rigid.
\end{lem}

\begin{demo}
 If $T' \in \T$, applying $\Hom_{\Cr {\Gamma}}(-, k[{\Gamma}] \tens T')$ to the admissible short exact sequence yields the long exact sequence
 \begin{align*}
  0 &\rightarrow \Hom_{\Cr {\Gamma}}(Y, k[{\Gamma}] \tens T') \rightarrow \Hom_{\Cr {\Gamma}}(T, k[{\Gamma}] \tens T') \xrightarrow{\Hom_{\Cr {\Gamma}}(f, k[{\Gamma}] \tens T')} \Hom_{\Cr {\Gamma}}(X, k[{\Gamma}] \tens T') \\
    &\rightarrow \Ext^1_{\Cr {\Gamma}}(Y, k[{\Gamma}] \tens T') \rightarrow \Ext^1_{\Cr {\Gamma}}(T, k[{\Gamma}] \tens T') = 0.
 \end{align*}
 As $f$ is a left $\T$-approximation, $\Hom_{\Cr {\Gamma}}(f, k[{\Gamma}] \tens T')$ is surjective, and as a consequence, $\Ext^1_{\Cr {\Gamma}}(Y, k[{\Gamma}] \tens T') = \Ext^1_{\Cr {\Gamma}}(k[{\Gamma}] \tens Y, T') = 0$. Moreover, as $k[{\Gamma}] \tens -$ is an exact functor, $0 \rightarrow k[{\Gamma}] \tens X \rightarrow k[{\Gamma}] \tens T \xrightarrow{k[{\Gamma}] \tens g} k[{\Gamma}] \tens Y \rightarrow 0$ is an admissible short exact sequence. Therefore, applying the functor $\Hom_{\Cr {\Gamma}}(X, -)$ gives rise to the long exact sequence
 \begin{align*}
  0 &\rightarrow \Hom_{\Cr {\Gamma}}(X, k[{\Gamma}] \tens X) \rightarrow \Hom_{\Cr {\Gamma}}(X, k[{\Gamma}] \tens T) \xrightarrow{\Hom_{\Cr {\Gamma}}(X, k[{\Gamma}] \tens g)} \Hom_{\Cr {\Gamma}}(X, k[{\Gamma}] \tens Y) \\
  &\rightarrow \Ext^1_{\Cr {\Gamma}}(X, k[{\Gamma}] \tens X) = 0
 \end{align*}
 because $k[{\Gamma}] \tens X$ is rigid. Furthermore, applying $\Hom_{\Cr {\Gamma}}(-, k[{\Gamma}] \tens Y)$ to the first admissible short exact sequence yields the long exact sequence
 \begin{align*}
  0 &\rightarrow \Hom_{\Cr {\Gamma}}(Y, k[{\Gamma}] \tens Y) \rightarrow \Hom_{\Cr {\Gamma}}(T, k[{\Gamma}] \tens Y) \xrightarrow{\Hom_{\Cr {\Gamma}}(f, k[{\Gamma}] \tens Y)} \Hom_{\Cr {\Gamma}}(X, k[{\Gamma}] \tens Y) \\
   &\rightarrow \Ext^1_{\Cr {\Gamma}}(Y, k[{\Gamma}] \tens Y) \rightarrow \Ext^1_{\Cr {\Gamma}}(T, k[{\Gamma}] \tens Y) = 0
 \end{align*}
 Let now $h \in \Hom_{\Cr {\Gamma}}(X, k[{\Gamma}] \tens Y)$. By the previous argument, $h$ factorizes through $k[{\Gamma}] \tens g$. Let $h' \in \Hom_{\Cr {\Gamma}}(X, k[{\Gamma}] \tens T)$ be such that $h = (k[{\Gamma}] \tens g)h'$. As $k[{\Gamma}] \tens T \in \T$ and $f$ is a left $\T$-approximation, there exists $t \in \Hom_{\Cr {\Gamma}}(T, k[{\Gamma}] \tens T)$ such that $h' = tf$. Hence, $h = (k[{\Gamma}] \tens g)tf$ and $\Hom_{\Cr {\Gamma}}(f, k[{\Gamma}] \tens Y)$ is surjective. Therefore, $\Ext^1_{\Cr {\Gamma}}(Y, k[{\Gamma}] \tens Y) = 0$ and $T \oplus k[{\Gamma}] \tens Y$ is rigid. Finally, $\add(\T, k[{\Gamma}] \tens Y)$ is rigid. \cqfd
\end{demo}

\begin{prop}
 \label{propdfmut}
 Let $\T \in \adds(\Cr {\Gamma})^{\md k[{\Gamma}]}$ rigid, and let $X \in \Cr {\Gamma}$ be indecomposable such that $X \notin \T$. Suppose that $\T$ contains all projective objects of $\Cr {\Gamma}$, and that $\add(\T, k[{\Gamma}] \tens X)$ is rigid. Then, there exist two admissible short exact sequences which are unique up to isomorphism 
 $$0 \rightarrow X \xrightarrow{f} T \xrightarrow{g} Y \rightarrow 0 \quad \text{and} \quad 0 \rightarrow Y' \xrightarrow{f'} T' \xrightarrow{g'} X \rightarrow 0$$
 such that
 \begin{enumerate}
  \item $f$ and $f'$ are minimal left $\T$-approximations; \label{propdfmut1}
  \item $g$ and $g'$ are minimal right $\T$-approximations; \label{propdfmut2}
  \item $\add(\T, k[{\Gamma}] \tens Y)$ and $\add(\T, k[{\Gamma}] \tens Y')$ are rigid; \label{propdfmut3}
  \item $Y \notin \T$ and $Y' \notin \T$; \label{propdfmut4}
  \item $Y$ and $Y'$ are indecomposable; \label{propdfmut5}
  \item $\add(k[{\Gamma}] \tens X) \inter \add(k[{\Gamma}] \tens Y) = 0$ and $\add(k[{\Gamma}] \tens X) \inter \add(k[{\Gamma}] \tens Y') = 0$. \label{propdfmut6}
 \end{enumerate}
\end{prop}

\begin{demo}
 By symmetry, it is enough to prove the results for the first admissible short exact sequence. Using lemma \ref{admapp}, a minimal left $\T$-approximation of $X$ is an admissible monomorphism, which implies the existence and the unicity of the admissible short exact sequence. Then (\ref{propdfmut1}) is satisfied by definition. As $X \notin \T$, the admissible short exact sequence does not split. Hence, $Y \notin \add(\T, k[{\Gamma}] \tens X)$ which is rigid. This proves (\ref{propdfmut4}) and (\ref{propdfmut6}). Moreover, lemma \ref{cokapp} shows that $\add(\T, k[{\Gamma}] \tens Y)$ is rigid. Hence, (\ref{propdfmut3}) is proved. For $T' \in \T$, applying $\Hom_{\Cr {\Gamma}}(T', -)$ to the admissible short exact sequence yields the long exact sequence
 \begin{align*}
  0 &\rightarrow \Hom_{\Cr {\Gamma}}(T', X) \rightarrow \Hom_{\Cr {\Gamma}}(T', Y) \xrightarrow{\Hom_{\Cr {\Gamma}}(T', g)} \Hom_{\Cr {\Gamma}}(T', Z) \\
    &\rightarrow \Ext^1_{\Cr {\Gamma}}(T', X) = 0
 \end{align*}
 and, therefore, $\Hom_{\Cr {\Gamma}}(T', g)$ is surjective, which implies that $g$ is a right $\T$-approximation.
 
 Let $T_0$ be a direct summand of $T$ on which $g$ vanishes. Let $\pi$ be a projection on $T_0$. As $g \pi = 0$ and $f$ is a kernel of $g$, there exists $f' \in \Hom_{\Cr {\Gamma}}(T, X)$ such that $\pi = ff'$. Therefore, $ff'(\id_T - \pi) = 0$ and, as $f$ is an admissible monomorphism, $f'(\id_T - \pi) = 0$. Thus $(f'f)^2 = f' \pi f = f'f$. As every idempotent splits and $X$ is indecomposable, $f'f = 0$ or $f'f = \id_X$. As the admissible short exact sequence does not split, $f'f = 0$. Finally, $\pi = ff' = (ff')^2 = 0$ so that $T_0 = 0$ and $g$ is minimal. This proves (\ref{propdfmut2}).
 
 Suppose now that $Y = Y_1 \oplus Y_2$. Let $g_1: T_1 \rightarrow Y_1$ and $g_2: T_2 \rightarrow Y_2$ be minimal right  $\T$-approximations. Then $g_1 \oplus g_2$ is a minimal right $\T$-approximation using lemmas \ref{sdminimal} and \ref{sdapprox}. By unicity of minimal approximations, $g \simeq g_1 \oplus g_2$. Therefore, as $g$ is an admissible epimorphism, $g_1$ and $g_2$ are admissible epimorphisms. Hence they have kernels $f_1: X_1 \rightarrow T_1$ and $f_2: X_2 \rightarrow T_2$ and by unicity of the kernel, $f \simeq f_1 \oplus f_2$. As a consequence, as $X \simeq X_1 \oplus X_2$ is indecomposable, $X_1 = 0$ or $X_2 = 0$ and as $f$ is minimal, $T_1 = 0$ or $T_2 = 0$. So, $Y_1 = 0$ or $Y_2 = 0$. Finally, $Y$ is indecomposable. This proves (\ref{propdfmut5}). \cqfd
\end{demo}

One now defines left and right mutations. Corollary \ref{mutvois} below claims that under some mild assumptions, these two notions coincide.

\begin{df}
 Retaining the notation of proposition \ref{propdfmut}, one writes
  $$\mu^r_X(\add(\T, k[{\Gamma}] \tens X)) = \add(\T, k[{\Gamma}] \tens Y) \quad \mu^l_X(\add(\T, k[{\Gamma}] \tens X)) = \add(\T, k[{\Gamma}] \tens Y').$$
 The map $\mu^r_X$ is called the \emph{right $X$-mutation} and $\mu^l_X$ is called the \emph{left $X$-mutation}. The ordered pair $(X, Y)$ (resp. $(X, Y')$) is called a \emph{right (resp. left) exchange pair associated with $\T$}.
\end{df}

\begin{rem}
 This is not ambiguous. Indeed, if $\T' \in \Add(\Cr {\Gamma})^{\md k[{\Gamma}]}$ and $X \in \T'$ is indecomposable then there exists a unique $\T \in \Add(\T')^{\md k[{\Gamma}]}$ such that $\add(k[{\Gamma}] \tens X) \inter \T = 0$ and $\T' = \add(\T, k[{\Gamma}] \tens X)$: it is the full subcategory of $\T'$ consisting of objects which have no common factors with $k[{\Gamma}] \tens X$.
\end{rem}

\begin{lem}
 \label{invmdG}
 Let $(X, Y)$ be a left (resp. right) exchange pair. Then,
 \begin{enumerate}
  \item $\# X = \# Y$ and $\ell(X) = \ell(Y)$;
  \item if $(X', Y')$ is another left (resp. right) exchange pair and if $X' \in \add(k[{\Gamma}] \tens X)$, then $Y' \in \add(k[{\Gamma}] \tens Y)$.
 \end{enumerate}
\end{lem}

\begin{demo}
 Let $0 \rightarrow X \rightarrow T \rightarrow Y \rightarrow 0$ be an admissible short exact sequence corresponding to the exchange pair. Using lemmas \ref{sdapprox}, \ref{sdminimal}, \ref{Fminapprox} and the fact that the class of admissible monomorphisms (resp. epimorphisms) is stable under direct summands, one proves that $0 \rightarrow FX \rightarrow FT \rightarrow FY \rightarrow 0$ is a direct sum of admissible short exact sequences of the same form as in proposition \ref{propdfmut}. Thus, $\ell(X) = \ell(Y)$. Moreover, it is clear that for two isomorphic indecomposable direct summands of $FX$, the two corresponding indecomposable direct summands of $FY$ are also isomorphic, and conversely. As a consequence, $\# X = \# Y$. Let $0 \rightarrow X' \rightarrow T' \rightarrow Y' \rightarrow 0$ be an admissible short exact sequence corresponding to the exchange pair $(X', Y')$. If $X_0$ is an indecomposable direct summand of $FX$, then, as $\add(k[{\Gamma}] \tens X) = \add(X_0[{\Gamma}]) = \add(k[{\Gamma}] \tens X')$, $X_0$ is also an indecomposable direct summand of $FX'$. Hence, the admissible short exact sequence $0 \rightarrow X_0 \xrightarrow{f} T_0 \rightarrow Y_0 \rightarrow 0$ where $f$ is a minimal left $\add(FT')$-approximation appears both as a direct summand of $0 \rightarrow FX \rightarrow FT \rightarrow FY \rightarrow 0$ and as a direct summand of  $0 \rightarrow FX' \rightarrow FT' \rightarrow FY' \rightarrow 0$. Finally $FY$ and $FY'$ have $Y_0$ as a common direct summand, which implies that $Y' \in \add(Y_0[{\Gamma}]) = \add(k[{\Gamma}] \tens Y)$ because $Y$ and $Y'$ are indecomposable. \cqfd
\end{demo}

\begin{lem} \label{lemv}
 Let $X, Y \in \Cr \Gamma$ be indecomposable, with $\ell(X) = \ell(Y)$ and $\# X = \# Y$. The following are equivalent:
 \begin{enumerate} 
  \item \label{lemv1} For every indecomposable object $X' \in \add(k[{\Gamma}] \tens X)$, 
  $$\dim \Ext^1_{\Cr {\Gamma}}(X', Y) = \left\{ \begin{array}{ll} 1 & \quad \text{if } X' \simeq X \\ 0 & \quad \text{else.} \end{array} \right.$$
  \item \label{lemv2} For every indecomposable object $Y' \in \add(k[{\Gamma}] \tens Y)$, 
  $$\dim \Ext^1_{\Cr {\Gamma}}(X, Y') = \left\{ \begin{array}{ll} 1 & \quad \text{if } Y' \simeq Y \\ 0 & \quad \text{else.} \end{array} \right.$$
 \end{enumerate}
\end{lem}

\begin{demo}
 Let $\tilde X$ (resp. $\tilde Y$) be the set of isomorphism classes of indecomposable direct summands of $k[{\Gamma}] \tens X$ (resp. $k[{\Gamma}] \tens Y$). Using lemma \ref{nbiso},
 \begin{align}
  \dim \Ext^1_\Cr(FX, FY) &= \dim \Ext^1_{\Cr {\Gamma}}(k[{\Gamma}] \tens X, Y) = \sum_{X' \in \tilde X} \frac{\ell(X) \ell(X')}{\# X'}\dim \Ext^1_{\Cr {\Gamma}}(X', Y) \label{eq1} \\
  &= \dim \Ext^1_{\Cr {\Gamma}}(X, k[{\Gamma}] \tens Y) = \sum_{Y' \in \tilde Y} \frac{\ell(Y) \ell(Y')}{\# Y'}\dim \Ext^1_{\Cr {\Gamma}}(X, Y') \label{eq2}
 \end{align}
 If (\ref{lemv1}) holds, then (\ref{eq1}) yields 
 $$\dim \Ext^1_\Cr(FX, FY) = \ell(X)^2/\# X.$$ As $\dim \Ext^1_{\Cr {\Gamma}}(X, Y) = 1$, the corresponding term in (\ref{eq2}) is equal to $\ell(Y)^2/\# Y = \ell(X)^2/\# X$, and therefore all other terms vanish. The converse is proved similarly. \cqfd
\end{demo}

\begin{df}
 With the notation of lemma \ref{lemv}, if the two equivalent assumptions are satisfied, $X$ and $Y$ are called \emph{neighbours}.
\end{df}

\begin{lem}
 \label{neicr}
 If $X, Y \in \Cr {\Gamma}$ are indecomposable objects satisfying $\ell(X) = \ell(Y)$ and $\# X = \# Y$, the following are equivalent:
  \begin{enumerate}
  \item For every indecomposable object $X_0 \in \add(FX)$, there exists an indecomposable object $Y_0 \in \add(FY)$ such that for every indecomposable $Y_0' \in \add(FY)$, 
   $$\dim \Ext^1_\Cr(X_0, Y_0') = \left\{\begin{array}{ll} 1 & \quad \text{if } Y_0' \simeq Y_0 \\ 0 & \quad \text{else.} \end{array} \right.$$ \label{neicr1}
  \item For every indecomposable object $Y_0 \in \add(FY)$, there exists an indecomposable object $X_0 \in \add(FX)$ such that for every indecomposable $X_0' \in \add(FX)$, 
  $$\dim \Ext^1_\Cr(X_0', Y_0) = \left\{\begin{array}{ll} 1 & \quad \text{if } X_0' \simeq X_0 \\ 0 & \quad \text{else.} \end{array} \right.$$ \label{neicr2}
  \item $\dim \Ext^1_\Cr(FX, FY) = \ell(X)^2/\# X$. \label{neicr3}
 \end{enumerate}
 Moreover, if $X$ and $Y$ are neighbours then these three conditions are satisfied.
\end{lem}

\begin{demo}
 By symmetry, it is enough to prove that (\ref{neicr1}) and (\ref{neicr3}) are equivalent. Let $X_0 \in \add(FX)$. Let $FX = \bigoplus_{k=1}^n X_i$ and $FY = \bigoplus_{k=1}^n Y_i$ where the $X_i$ and the $Y_i$ are indecomposable in $\Cr$. Let $\tilde X$ be a set of representatives of isomorphism classes of indecomposable summands of $FX$ and $\tilde Y$ a set of representatives of isomorphism classes of indecomposable summands of $FY$.
 \begin{align*}
  \dim \Ext^1_\Cr(FX, FY) &= \sum_{i=1}^n \sum_{j=1}^n \dim \Ext^1_\Cr(X_i, Y_j) = \left(\frac{\ell(X)}{\# X}\right)^2 \sum_{X_0' \in \tilde X} \sum_{Y_0' \in \tilde Y} \dim \Ext^1_\Cr(X_0', Y_0') \\ &= \frac{\ell(X)^2}{\# X} \sum_{Y_0' \in \tilde Y} \dim \Ext^1_\Cr(X_0, Y_0')
 \end{align*}
 which yields the equivalence because the $\dim \Ext^1_\Cr(X_0, Y_0')$ are non-negative integers.
 
 The proof of lemma \ref{invmdG} implies that if $X$ and $Y$ are neighbours, then (\ref{neicr3}) is satisfied. \cqfd
\end{demo}

\begin{prop}
 \label{neighapp} Let $X, Y \in \Cr {\Gamma}$ be neighbours such that $k[{\Gamma}] \tens X$ and $k[{\Gamma}] \tens Y$ are rigid. Let $$0 \rightarrow X \xrightarrow{f} M \xrightarrow{g} Y \rightarrow 0$$ be a non-split admissible short exact sequence (which is unique up to isomorphism because $X$ and $Y$ are neighbours). 
 
 Then $\add(k[{\Gamma}] \tens (M \oplus X))$ and $\add(k[{\Gamma}] \tens (M \oplus Y))$ are rigid and $X, Y \notin \add(k[{\Gamma}] \tens M)$. 
 
 Moreover, if there is $\T \in \Add(\Cr {\Gamma})^{\md k[{\Gamma}]}$ such that $\add(\T, k[{\Gamma}] \tens X)$ and $\add(\T, k[{\Gamma}] \tens Y)$ are maximal $\md k[{\Gamma}]$-stable rigid, then $f$ is a minimal left $\T$-approximation and $g$ is a minimal right $\T$-approximation.
\end{prop}

\begin{demo}
 In order to show that $\Ext^1_{\Cr {\Gamma}} (k[{\Gamma}] \tens M, k[{\Gamma}] \tens X) = 0$, it is enough to show that $\Ext^1_{\Cr {\Gamma}}(M, X') = 0$ for every indecomposable object $X' \in \add(k[{\Gamma}] \tens X)$. Let $X' \in \add(k[{\Gamma}] \tens X)$ be indecomposable. Applying $\Hom_{\Cr  {\Gamma}}(-,X')$ to the admissible short exact sequence yields the following long exact sequence:
 \begin{align*}
  0 &\rightarrow \Hom_{\Cr {\Gamma}}(Y, X') \rightarrow \Hom_{\Cr {\Gamma}}(M, X') \xrightarrow{\Hom_{\Cr {\Gamma}}(f, X')} \Hom_{\Cr {\Gamma}}(X, X') \\
    &\xrightarrow{\delta} \Ext^1_{\Cr {\Gamma}}(Y, X') \rightarrow \Ext^1_{\Cr {\Gamma}}(M, X') \rightarrow \Ext^1_{\Cr {\Gamma}}(X, X') = 0
 \end{align*}
 If $X' \simeq X$ then $\dim \Ext^1_{\Cr {\Gamma}}(Y, X') = 1$. Thus, it is enough to prove that $\delta \neq 0$, which is equivalent to $\Hom_{\Cr {\Gamma}}(f, X')$ is not surjective. If it was surjective, there would exist $g \in \Hom_{\Cr  {\Gamma}}(M, X')$ such that $gf$ is an isomorphim. As a consequence, $[(gf)^{-1}g] f = \id_X$ and the admissible short exact sequence would split, which is not the case. 
 
 If $X'$ and $X$ are not isomorphic, as $X$ and $Y$ are neighbours, $\Ext^1_{\Cr {\Gamma}}(Y, X') = 0$ and the result is clear. Similarly, $\Ext^1_{\Cr {\Gamma}} (k[{\Gamma}] \tens M, k[{\Gamma}] \tens Y) = 0$.
 
 In order to show that $\Ext^1_{\Cr {\Gamma}}(k[{\Gamma}] \tens M, k[{\Gamma}] \tens M) = 0$, it is enough to prove that 
 $$\Ext^1_{\Cr {\Gamma}}(k[{\Gamma}] \tens M, M) = 0.$$
 Applying the functor $\Hom_{\Cr {\Gamma}}(k[{\Gamma}] \tens M, -)$ to the admissible short exact sequence induces the following long exact sequence:
 $$\Ext^1_{\Cr {\Gamma}}(k[{\Gamma}] \tens M, X) = 0\rightarrow \Ext^1_{\Cr {\Gamma}}(k[{\Gamma}] \tens M, M) \rightarrow \Ext^1_{\Cr {\Gamma}}(k[{\Gamma}] \tens M, Y) = 0$$
 which yields the result. 
 
 If $X$ was in $\add(k[{\Gamma}] \tens M)$, there would exist an indecomposable object $X' \in \add(k[{\Gamma}] \tens X)$ such that $X' \in \add(M)$. Let $M = X' \oplus M'$. If $X' \simeq X$, then $\dim \Ext^1_{\Cr {\Gamma}}(Y, X) = 1$ and $\Ext^1_{\Cr {\Gamma}}(Y, M) = 0$ are contradictory. Suppose that $X'$ and $X$ are not isomorphic. Then there exists an indecomposable object $Y' \in \add(k[{\Gamma}] \tens Y)$ such that $\Ext^1_{\Cr {\Gamma}}(X', Y') \neq 0$ and, by definition of neighbours, $Y'$ and $Y$ are not isomorphic. Applying the functor $\Hom_{\Cr {\Gamma}}(Y', -)$ to the admissible short exact sequence yields the long exact sequence
 $$\Ext^1_{\Cr {\Gamma}}(Y', X) = 0 \rightarrow \Ext^1_{\Cr {\Gamma}}(Y', X' \oplus M') \rightarrow \Ext^1_{\Cr {\Gamma}}(Y', Y) = 0$$
 the first equality coming from the fact that $Y'$ and $Y$ are not isomorphic and that $X$ and $Y$ are neighbours. As a consequence, the central term vanishes, which contradicts the hypothesis.
 
 Suppose now that $\T$ exists. For $T \in \T$, applying $\Hom_{\Cr {\Gamma}}(T,-)$ yields the following long exact sequence
 $$\Ext^1_{\Cr {\Gamma}}(T, X) = 0 \rightarrow \Ext^1_{\Cr {\Gamma}}(T, M) \rightarrow \Ext^1_{\Cr {\Gamma}}(T, Y) = 0$$
 and therefore, $\Ext^1_{\Cr {\Gamma}}(T, M) = 0$. Hence $\add(\T, k[{\Gamma}] \tens (M \oplus X)))$ is $\md k[{\Gamma}]$-stable rigid. As a consequence, as $\add(\T, k[{\Gamma}] \tens X)$ is maximal $\md k[{\Gamma}]$-stable rigid, $M \in \add(\T, k[{\Gamma}] \tens X)$. As $\add(k[{\Gamma}] \tens X) \inter \add(k[{\Gamma}] \tens M) = 0$, one gets $M \in \T$. Applying $\Hom_{\Cr {\Gamma}}(-, T)$ gives rise to the following long exact sequence:
 \begin{align*}
  0 &\rightarrow \Hom_{\Cr {\Gamma}}(Y, T) \rightarrow \Hom_{\Cr {\Gamma}}(M, T) \xrightarrow{\Hom_{\Cr {\Gamma}}(f, T)} \Hom_{\Cr {\Gamma}}(X, T) \\&\rightarrow \Ext^1_{\Cr {\Gamma}}(Y, T) = 0
 \end{align*}
 and, as a consequence, $\Hom_{\Cr {\Gamma}}(f, T)$ is surjective for every $T \in \T$, and the morphism $f$ is a left $\T$-approximation. If $f$ were not minimal, there would exist a decomposition $M \simeq M_0 \oplus M_1$ with $M_0 \neq 0$ such that $f = f_0 \oplus f_1$ where $f_1 = \pi_{M_1} f$ and $f_0 = 0$. Thus, $f_0$ and $f_1$ would be admissible monomorphisms. As $Y$ is indecomposable, one of the cokernels of $f_0$ and $f_1$ vanishes. As $\coker f_0 = M_0 \neq 0$, $\coker f_1 = 0$, and, as a consequence, the admissible short exact sequence splits which is a contradiction. In the same way, $g$ is a minimal right $\T$-approximation. \cqfd
\end{demo}

\begin{cor}
 \label{mutvois}
 Let $\T \in \Add(\Cr {\Gamma})^{\md k[{\Gamma}]}$ and $X \in \Cr {\Gamma}$ be an indecomposable object such that $X \notin \T$ and $\add(\T, k[{\Gamma}] \tens X)$ is maximal $\md k[{\Gamma}]$-stable rigid. Then, the following are equivalent:
 \begin{enumerate}
  \item There exists an indecomposable object $Y \in \Cr {\Gamma}$ such that $\mu^l_X(\add(\T, k[{\Gamma}] \tens X)) = \add(\T, k[{\Gamma}] \tens Y)$ and $X$ and $Y$ are neighbours. \label{mutvois1}
  \item There exists an indecomposable object $Y' \in \Cr {\Gamma}$ such that $\mu^r_X(\add(\T, k[{\Gamma}] \tens X)) = \add(\T, k[{\Gamma}] \tens Y')$ and $X$ and $Y'$ are neighbours. \label{mutvois2}
 \end{enumerate}
 In this case, $Y \simeq Y'$ and if one denotes
 $$\mu_X(\add(\T, k[{\Gamma}] \tens X)) = \mu^l_X(\add(\T, k[{\Gamma}] \tens X)) = \mu^r_X(\add(\T, k[{\Gamma}] \tens X)),$$
 then
 $$\mu_Y(\mu_X(\add(\T, k[{\Gamma}] \tens X)) = \add(\T, k[{\Gamma}] \tens X).$$
\end{cor}

\begin{demo}
 If (\ref{mutvois1}) is true, proposition \ref{neighapp} induces the admissible short exact sequence $0 \rightarrow X \rightarrow T' \rightarrow Y \rightarrow 0$ satisfying the conditions of proposition \ref{propdfmut}, which proves (\ref{mutvois2}) and the fact that $Y \simeq Y'$. Hence $\mu_Y(\mu_X(\add(\T, k[{\Gamma}] \tens X)) = \mu^l_Y(\mu^r_X(\add(\T, k[{\Gamma}] \tens X)) = \add(\T, k[{\Gamma}] \tens X)$. By a similar argument, (\ref{mutvois2}) implies (\ref{mutvois1}). \cqfd
\end{demo}

\begin{df}
 In the situation of the corollary, we write
  $$\mu_X(\add(\T, k[{\Gamma}] \tens X)) = \add(\T, k[{\Gamma}] \tens Y),$$
  and we say that $\{X, Y\}$ is \emph{an exchange pair associated with $\T$}.
\end{df}

\subsection{Rigid quasi-approximations}

\begin{df}
 Let $X \in \Cr {\Gamma}$ (resp. $\in \Cr$). An epimorphism $f: X \twoheadrightarrow Y$ will be called a \emph{left rigid quasi-approximation of $X$} if the following conditions are satisfied:
 \begin{itemize}
  \item $k[{\Gamma}] \tens Y$ (resp. $\bigoplus_{g \in {\Gamma}} \g \tens Y$) is rigid;
  \item $Y$ has an injective envelope, without direct summand in $$\add(k[{\Gamma}] \tens X) \quad \left(\text{resp. }\add\left(\bigoplus_{g \in {\Gamma}} \g \tens X\right)\right);$$
  \item If $k[{\Gamma}] \tens Z$ (resp. $\bigoplus_{g \in {\Gamma}} \g \tens Z$) is rigid, then every morphism from $X$ to $Z$, without invertible matrix coefficient, factorizes through $f$.
 \end{itemize}
\end{df}

\begin{rem}
 The definition of a right rigid quasi-approximation is obtained from the previous one by reversing the arrows. All the following results can be adapted to this case.
\end{rem}

\begin{rem}
 As for the case of approximations, minimal quasi-approximations are unique up to (non unique) isomorphism.
\end{rem}

The following lemma gives an easy way to get quasi-approximations:

\begin{lem}
 \label{qapsc}
 Let $\Dr$ be an abelian category endowed with an action of ${\Gamma}$. Let $\Cr$ be an exact full subcategory of $\Dr$. Let $P$ be a projective indecomposable object of $\Cr$. If
 \begin{itemize}
  \item every monomorphism of $\Cr$ to a projective object is admissible,
  \item every monomorphism of $\Cr$ from $P$ to an indecomposable object is admissible,
  \item $P$ has a simple socle $S$ in $\Dr$,
  \item $\bigoplus_{g \in {\Gamma}} \g \tens S$ is rigid in $\Dr$,
  \item the cokernel of $S \inj P$ is in $\Cr$,
 \end{itemize}
 then this cokernel is a left rigid quasi-approximation of $P$ in $\Cr$.
\end{lem}

\begin{demo}
 Consider the following short exact sequence in $\Dr$:
 $$0 \rightarrow S \rightarrow P \rightarrow Q \rightarrow 0.$$
 
 Applying the functor $\Hom_\Dr\left(\bigoplus_{g \in {\Gamma}} \g \tens S, -\right)$ yields the following long exact sequence:
 \begin{align*}
  0 &\rightarrow \Hom_\Dr\left(\bigoplus_{g \in {\Gamma}} \g \tens S, S\right) \rightarrow \Hom_\Dr\left(\bigoplus_{g \in {\Gamma}} \g \tens S, P\right) \rightarrow \Hom_\Dr\left(\bigoplus_{g \in {\Gamma}} \g \tens S, Q\right) \\
  &\rightarrow \Ext^1_\Dr\left(\bigoplus_{g \in {\Gamma}} \g \tens S, S\right) = 0.
 \end{align*}
 As the socle of $P$ is simple, 
 $$\dim_k \Hom_\Dr\left(\bigoplus_{g \in {\Gamma}} \g \tens S, P\right) = \# \{g \in {\Gamma} \,|\, \g \tens S \simeq S\} = \dim_k \Hom_\Dr\left(\bigoplus_{g \in {\Gamma}} \g \tens S, S\right)$$
 and therefore $\Hom_\Dr\left(\bigoplus_{g \in {\Gamma}} \g \tens S, Q\right) = 0$.
 
 Applying the functor $\Hom_\Dr\left(-, \bigoplus_{g \in {\Gamma}} \g \tens Q\right)$ yields the long exact sequence
 \begin{align*}
  0 &\rightarrow \Hom_\Dr\left(Q, \bigoplus_{g \in {\Gamma}} \g \tens Q\right) \rightarrow \Hom_\Dr\left(P, \bigoplus_{g \in {\Gamma}} \g \tens Q\right) \rightarrow \Hom_\Dr\left(S, \bigoplus_{g \in {\Gamma}} \g \tens Q\right) \\
  &\rightarrow \Ext^1_\Dr\left(Q, \bigoplus_{g \in {\Gamma}} \g \tens Q\right) \rightarrow \Ext^1_\Dr\left(P, \bigoplus_{g \in {\Gamma}} \g \tens Q\right) = 0 
 \end{align*}
 and, as $\Hom_\Dr\left(S, \bigoplus_{g \in {\Gamma}} \g \tens Q\right) = 0$, one gets  
 $$\Ext^1_\Dr\left(Q, \bigoplus_{g \in {\Gamma}} \g \tens Q\right) = 0.$$ 
 
 Let now $S'$ be the (not necessarily simple) socle of $Q$. As 
 $$\Hom_\Dr\left(\bigoplus_{g \in {\Gamma}} \g \tens S, Q\right) = 0,$$
 $S'$ does not contain any direct summand isomorphic to some $\g \tens S$. Hence, if one denotes $\alpha: S' \rightarrow I$ the injective envelope of $S'$, $I$ has no direct summand of the form $\g \tens P$. Moreover, as $I$ is injective, one can complete the following commutative diagram:
 $$\xymatrix{
  S' \ar[r]^\alpha \ar[d]_\iota & I \\
  Q \ar[ur]_\beta &
 }$$
 Here, $\beta$ is an admissible monomorphism since $I$ is injective. Hence, the second condition of the definition of a left quasi-approximation is proved. 
 
 Let now $Z \in \Cr$ be indecomposable and $f: P \rightarrow Z$. Suppose that $f$ is not a monomorphism. Let $K \rightarrow P$ be the kernel of $f$ in $\Dr$ and let $S'$ be the socle of $K$. By composing the two morphisms, one has a non vanishing morphism $g: S' \rightarrow P$ such that $fg = 0$. As the socle of $P$ is simple, $S' = S$ and finally, $f$ vanishes on $S$ and therefore factorizes through $P \rightarrow Q$ which is the cokernel of $S \rightarrow P$. If $f$ is a monomorphism, by hypothesis, it is admissible in $\Cr$ so it splits and finally $f = \id_P$ as $Z$ is indecomposable. The third condition is proved. \cqfd
\end{demo}

\begin{lem}
 If every indecomposable projective object of $\Cr$ has a left rigid quasi-approximation, then every indecomposable projective object of $\Cr {\Gamma}$ has a left rigid quasi-approximation.
\end{lem}

\begin{demo}
 Let $(P, \psi) \in \Cr  {\Gamma}$ be an indecomposable projective object. By lemma \ref{carinjproj}, $P \in \Cr$ is projective. Let $f: P \twoheadrightarrow X$ be the direct sum of the minimal left rigid quasi-approximations of its indecomposable direct summands. It is also a minimal left rigid quasi-approximation since $$\add X = \add\left(\bigoplus_{g \in {\Gamma}} \g \tens X_0\right)$$ where $P_0 \twoheadrightarrow X_0$ is one of the minimal left rigid quasi-approximations of the indecomposable direct summands of $P$. Hence $\bigoplus_{g \in {\Gamma}} \g \tens X$ is rigid. Therefore, for all $g \in {\Gamma}$, there exists a unique morphism $\chi_g$ which makes the following diagram commutative:
 $$\xymatrix{
  \g \tens P \ar[r]^{\g \tens f} \ar[d]_{\psi_g} & \g \tens X \ar[d]_{\chi_g} \\
  P \ar[r]^{f} & X
 }$$
 (The existence comes from the definition of a rigid quasi-approximation, the unicity comes from the fact that $f$ is an epimorphism). Clearly, $(X, \chi) \in \Cr {\Gamma}$. So it is easy to see that $f: (P, \psi) \twoheadrightarrow (X, \chi)$ is a left rigid quasi-approximation. \cqfd
\end{demo}

\subsection{Endomorphisms}

All projective indecomposable objects of $\Cr$ will be supposed here to have left rigid quasi-approximations. All results remain valid if they are supposed to have right rigid quasi-approximations.

Let $\T \in \adds(\Cr {\Gamma})^{\md k[{\Gamma}]}$ be maximal $\md k[{\Gamma}]$-stable rigid. Let $T \in \Cr {\Gamma}$ and $\tilde T \in FT$ be basic such that $\T = \add(T)$ and $F \T = \add(\tilde T)$ (one can find such $T$ and $\tilde T$ since $\T$ is finitely generated). Write $E = \End_{\Cr {\Gamma}}(T)$ and $\tilde E = \End_\Cr(\tilde T)$. If $X \in \T$ is indecomposable, denote by $S_X$ the corresponding simple representation of $E$, that is, the head of the projective $E$-module $\Hom_\Cr(X, T)$. Likewise, if $X \in F \T$ is indecomposable, denote by $\tilde S_X$ the corresponding simple representation of $\tilde E$.

\begin{df}
 Let $\Dr \in \Add(\Cr {\Gamma})^{\md k[{\Gamma}]}$ and let $X \in \Dr$ be indecomposable. A \emph{$\md k[{\Gamma}]$-loop of $\Dr$ at $X$} is an irreducible morphism $X \rightarrow X'$ of $\Dr$ where $X' \in \add(k[{\Gamma}] \tens X)$ is indecomposable. A \emph{$\md k[{\Gamma}]$-$2$-cycle of $\Dr$ at $X$} is a couple of irreducible morphisms $X \rightarrow Y$ and $Y \rightarrow X'$ of $\Dr$ where $X' \in \add(k[{\Gamma}] \tens X)$ is indecomposable.
\end{df}

\begin{df}
 Let $\Dr \in \Add(\Cr)^{\Gamma}$ and $X \in \Dr$ be indecomposable. A \emph{${\Gamma}$-loop of $\Dr$ at $X$} is an irreducible morphism $X \rightarrow \g \tens X$ of $\Dr$ where $g \in {\Gamma}$. A \emph{${\Gamma}$-$2$-cycle of $\Dr$ at $X$} is a couple of irreducible morphisms $X \rightarrow Y$ and $Y \rightarrow \g \tens X$ of $\Dr$ where $g \in {\Gamma}$.
\end{df}

\begin{lem}
 \label{loop} 
 Let $\Dr \in \Add(\Cr {\Gamma})^{\md k[{\Gamma}]}$. Let $X \in \Dr$ be indecomposable and $X'$ be a direct summand of $FX$. Then $\Dr$ has no $\md k[{\Gamma}]$-loops (resp. $\md k[{\Gamma}]$-$2$-cycles) at $X$ if and only if $F \Dr$ has no ${\Gamma}$-loops (resp. ${\Gamma}$-$2$-cycles) at $X'$.
\end{lem}

\begin{demo}
 The proof is the same for loops and $2$-cycles. Hence, it will be done only for loops. Suppose that $\Dr$ has no $\md k[{\Gamma}]$-loops at $X$. Let $f \in \End_{\Cr}(F X)$. Then $f[{\Gamma}] \in \End_{\Cr {\Gamma}}((F X)[{\Gamma}]) \simeq \End_{\Cr {\Gamma}}(k[{\Gamma}] \tens X)$. As $\Dr$ has no $\md k[{\Gamma}]$-loops at $X$, $f[{\Gamma}]$ factorizes through $Y \in \Dr$ such that $\add(Y) \inter \add(k[{\Gamma}] \tens X) = 0$. Hence $F(f[{\Gamma}])$ factorizes through $F Y$. As $\add(k[{\Gamma}] \tens X)$ is $\md k[{\Gamma}]$-stable,  $\add(F Y) \inter \add(F (k[{\Gamma}] \tens X)) = 0$. As $f$ is a direct summand of $F(f[{\Gamma}])$, $f$ factorizes also through $F Y$ and, as a consequence, $F \Dr$ has no ${\Gamma}$-loops at $X'$. 
 
 Conversely, suppose that $F \Dr$ has no ${\Gamma}$-loops at $X'$. Let $f \in \End_{\Cr {\Gamma}}(X)$. As $F f \in \End_{\Cr}(FX)$, $Ff$  factorizes through $Y \in \Dr$ such that $\add(Y) \inter \add \left(\bigoplus_{g \in {\Gamma}} \g \tens X' \right) = 0$. Thus, $(Ff)[{\Gamma}]$ factorizes through $Y[{\Gamma}]$ and $f$, as a direct summand of $(Ff)[{\Gamma}]$, factorizes also through $Y[{\Gamma}]$ and finally, as $$\add \left(\bigoplus_{g \in {\Gamma}} \g \tens X' \right) = \add F(X'[{\Gamma}]) = \add F X,$$ one gets $\add(Y[{\Gamma}]) \inter \add(X) = 0$. Eventually, $\Dr$ has no $\md k[{\Gamma}]$-loops at $X$. \cqfd
\end{demo}

\begin{lem}
 \label{psbclvois}
 Let $\T_0 \in \Add(\T)^{\md k[{\Gamma}]}$ and let $(X, Y)$ be a left (resp. right) exchange pair associated to $\T_0$ such that $\add(\T_0, k[{\Gamma}] \tens X) = \T$. The following are equivalent:
 \begin{enumerate}
  \item $\T$ has no $\md k[{\Gamma}]$-loops at $X$; \label{psbclvois1}
  \item For all indecomposable $X' \in \add(k[{\Gamma}] \tens X)$, every non invertible morphism from $X$ to $X'$ factorizes through $\T_0$. \label{psbclvois2}
  \item $X$ and $Y$ are neighbours. \label{psbclvois3}
 \end{enumerate}
\end{lem}

\begin{demo}
 For the proof, the exchange pair will be considered to be a left exchange pair.
 
 The equivalence of (\ref{psbclvois1}) and (\ref{psbclvois2}) is clear. Let
 $$0 \rightarrow X \xrightarrow{f} T' \rightarrow Y \rightarrow 0$$
 be the admissible short exact sequence corresponding to the exchange pair $(X, Y)$. Let $X' \in \add(k[{\Gamma}] \tens X)$. Applying $\Hom_{\Cr {\Gamma}}(-, X')$ leads to the long exact sequence:
 \begin{align*}
  0 &\rightarrow \Hom_{\Cr {\Gamma}}(Y, X') \rightarrow \Hom_{\Cr {\Gamma}}(T', X') \xrightarrow{\Hom_{\Cr {\Gamma}}(f, X')} \Hom_{\Cr {\Gamma}}(X, X') \\
  &\rightarrow \Ext^1_{\Cr {\Gamma}}(Y, X') \rightarrow \Ext^1_{\Cr {\Gamma}}(T', X') = 0.
 \end{align*}
 As $f$ is a left $\T_0$-approximation, for any element of $\Hom_{\Cr {\Gamma}}(X, X')$, factorizing through $\T_0$ is equivalent to being in the image of $\Hom_{\Cr {\Gamma}}(f, X')$. There are two cases:
 \begin{itemize}
  \item if $X' \simeq X$ then every non invertible morphism from $X$ to $X'$ factorizes through $\T_0$ if and only if $\Ext^1_{\Cr {\Gamma}}(Y, X') \simeq k$ because, as $k$ is algebraically closed, $\End_{\Cr {\Gamma}}(X) / \mathfrak{m} \simeq k$ where $\mathfrak{m}$ is the maximal ideal of $\End_{\Cr {\Gamma}}(X)$;
  \item if $X'$ and $X$ are not isomorphic then every morphism from $X$ to $X'$ factorizes through $\T_0$ if and only if $\Ext^1_{\Cr {\Gamma}}(Y, X') = 0$.
 \end{itemize}
 Combining these two cases, the equivalence of (\ref{psbclvois2}) and (\ref{psbclvois3}) is proved. \cqfd
\end{demo}

\begin{prop}
 \label{gldim}
 Suppose that $\T$ has no $\md k[{\Gamma}]$-loops. Then
 $$\gldim(E) = \gldim(\tilde E) \left\{\begin{array}{ll} \leq 2 & \quad \text{if } T \text{ is projective} \\ = 3 & \quad \text{else.} \end{array} \right.$$
\end{prop}

\begin{demo}
 The proof is the same for $E$ and for $\tilde E$ since $\tilde \T$ has no ${\Gamma}$-loops and since $T$ is projective if and only if $\tilde T$ is projective. Hence, it is enough to do it for $E$. Let $X \in \T$ be indecomposable and $\T_0 \in \T^{\md k[{\Gamma}]}$ be such that $\T = \add(\T_0, k[{\Gamma}] \tens X)$ and $X \notin \T_0$. Suppose that $X$ is not projective. There exist two exchange admissible short exact sequences
 $$0 \rightarrow X \xrightarrow{f} T' \rightarrow Y \rightarrow 0 \quad \text{and} \quad 0 \rightarrow Y \rightarrow T'' \rightarrow X \rightarrow 0$$ since, using the previous lemma, $X$ and $Y$ are neighbours and using corollary \ref{mutvois}, one gets $\mu_Y(\mu_X(\T)) = \T$. Applying $\Hom_{\Cr {\Gamma}}(-, T)$ to these sequences yields the following long exact sequences:
 \begin{align*}
  0 &\rightarrow \Hom_{\Cr {\Gamma}}(Y, T) \rightarrow \Hom_{\Cr {\Gamma}}(T', T) \rightarrow \Hom_{\Cr {\Gamma}}(X, T) \\  
  &\rightarrow \Ext^1_{\Cr {\Gamma}}(Y, T) \simeq \Ext^1_{\Cr {\Gamma}}(Y, X') \rightarrow \Ext^1_{\Cr {\Gamma}}(T', T) = 0
 \end{align*}
 where $X'$ is the largest direct summand of $T$ contained in $\add(k[{\Gamma}] \tens X)$ and
 \begin{align*}
  0 &\rightarrow \Hom_{\Cr {\Gamma}}(X, T) \rightarrow \Hom_{\Cr {\Gamma}}(T'', T) \rightarrow \Hom_{\Cr {\Gamma}}(Y, T) \rightarrow \Ext^1_{\Cr {\Gamma}}(X, T) = 0
 \end{align*}
 As $T$ is basic and $X'$ is $\md k[{\Gamma}]$-stable, $X'$ contains exactly one object of each isomorphism class of $\add(k[{\Gamma}] \tens X)$. As $X$ and $Y$ are neighbours, $\dim \Ext^1_{\Cr {\Gamma}}(Y, X') = 1$. Thus, $\Ext^1_{\Cr {\Gamma}}(Y, X')  \simeq S_X$ as an $E$-module. Therefore, combining these two long exact sequences yields the following long exact sequence of $E$-modules:
 $$0 \rightarrow \Hom_{\Cr {\Gamma}}(X, T) \rightarrow \Hom_{\Cr {\Gamma}}(T'', T) \rightarrow \Hom_{\Cr {\Gamma}}(T', T) \rightarrow \Hom_{\Cr {\Gamma}}(X, T) \rightarrow S_X \rightarrow 0$$
 which is a projective resolution of $S_X$. As a consequence, $\prdim(S_X) \leq 3$. As $X \notin \add(T'')$, $\Hom_E(\Hom_{\Cr {\Gamma}}(T'', T), S_X) = 0$ and therefore $$\Ext^3(S_X, S_X) \simeq \Hom_E(\Hom_{\Cr {\Gamma}}(X, T), S_X)$$ has dimension $1$. Finally, $\prdim(S_X) = 3$.
 
 Suppose now that $X$ is projective. Show that $\prdim(S_X) \leq 2$. Let $\pi: X \twoheadrightarrow Y$ be a left rigid  quasi-approximation. By definition, the injective envelope of $Y$ does not intersect $\add(k[{\Gamma}] \tens X)$ and, using lemma \ref{admapp}, there is an admissible short exact sequence
 $$0 \rightarrow Y \xrightarrow{f} T' \rightarrow Z \rightarrow 0$$
 where $f$ is a left $\add(\T_0)$-approximation. Using lemma \ref{cokapp} and as $X$ is projective, $\add(\T, k[{\Gamma}] \tens Z)$ is $\md k[{\Gamma}]$-stable rigid, and, since $\T$ is maximal $\md k[{\Gamma}]$-stable rigid, $Z \in \T$. Applying $\Hom_{\Cr {\Gamma}}(-,T)$ to this admissible short exact sequence yields the long exact sequence
 $$
  0 \rightarrow \Hom_{\Cr {\Gamma}}(Z,T) \rightarrow \Hom_{\Cr {\Gamma}}(T',T) \rightarrow \Hom_{\Cr {\Gamma}}(Y,T) \rightarrow \Ext^1_{\Cr {\Gamma}}(Z,T) = 0
 $$
 Moreover, $\Hom_{\Cr {\Gamma}}(\pi, T): \Hom_{\Cr {\Gamma}}(Y, T) \rightarrow \Hom_{\Cr {\Gamma}}(X,T)$ is injective because $\pi$ is an epimorphism, and its cokernel has dimension $1$: this cokernel is $S_X$. One deduces the following long exact sequence:
 $$
  0 \rightarrow \Hom_{\Cr {\Gamma}}(Z,T) \rightarrow \Hom_{\Cr {\Gamma}}(T',T) \rightarrow \Hom_{\Cr {\Gamma}}(X,T) \rightarrow S_X \rightarrow 0
 $$
 which is a projective resolution of $S_X$. \cqfd
\end{demo}

Recall this theorem of Happel:

\begin{thm}[\citeb{section 1.4}{Ha87}]
 \label{bascde}
 If $A$ is a $k$-algebra and $X$ is a tilting $A$-module, then $\md A$ and $\md \End_A(X)^{\op}$ are derived equivalent.
\end{thm}

The following proposition explains the relationship between any maximal $\md k[{\Gamma}]$-stable rigid category $\T'$ and the initial one $\T$.

\begin{prop}
 \label{deq}
 Suppose that $\T' \in \Add(\Cr {\Gamma})^{\md k[{\Gamma}]}$ is maximal $\md k[{\Gamma}]$-stable rigid. Then $\T'$ is finitely generated. Let $T' \in \T'$ be basic such that $\T' = \add(T')$. Let $E' = \End_{\Cr {\Gamma}}(T')$. Let $M = \Hom_{\Cr {\Gamma}}(T', T)$. Then $M$ is a tilting module on $E$ and
 $$\End_{E}(M) \simeq E'^{\op}.$$
 In particular, there is a derived equivalence between $E$ and $E'$ and $\T$ and $\T'$ contain the same number of indecomposable objects up to isomorphism.
 
 The same holds for $\tilde T'$ and $\tilde E'$.
\end{prop}

\begin{demo}
 Using lemma \ref{admapp}, one gets the admissible short exact sequence
 \begin{equation} \label{ses0}
  0 \rightarrow T'_2 \rightarrow T'_1 \xrightarrow{g} T \rightarrow 0
 \end{equation}
 where $g$ is a minimal right $\T'$-approximation. Using lemma \ref{cokapp}, $T'_2 \in \T'$. Let $T' \in \add(\T')$ be basic such that $T'_1, T'_2 \in \add(T')$ (in fact $\add(T') = \T'$ will be proved later).
 
 By the same argument, there is an admissible short exact sequence
 \begin{equation} \label{ses1}
 0 \rightarrow T' \xrightarrow{f} T_1 \rightarrow T_2 \rightarrow 0
 \end{equation}
 where $f$ is a minimal left $\T$-approximation and $T_2 \in \T$. Applying $\Hom_{\Cr {\Gamma}}(-, T)$ to (\ref{ses1}) yields the following long exact sequence:
 \begin{equation} \label{ses2}
 0 \rightarrow \Hom_{\Cr {\Gamma}}(T_2, T) \rightarrow \Hom_{\Cr {\Gamma}}(T_1, T) \rightarrow \Hom_{\Cr {\Gamma}}(T', T) = M \rightarrow \Ext^1_{\Cr {\Gamma}}(T_2, T) = 0
 \end{equation}
 and, as a consequence, $\prdim_E M \leq 1$. Now, applying the functor $\Hom_{\Cr {\Gamma}}(T', -)$ to (\ref{ses1}) gives the following long exact sequence:
  $$
  0 \rightarrow \Hom_{\Cr {\Gamma}}(T', T') \rightarrow \Hom_{\Cr {\Gamma}}(T', T_1) \rightarrow \Hom_{\Cr {\Gamma}}(T', T_2) \rightarrow \Ext^1_{\Cr {\Gamma}}(T', T') = 0.
  $$
  
  Applying the functor $\Hom_E(-, M)$ to (\ref{ses2}) induces the long exact sequence
  \begin{align*} 0 &\rightarrow \Hom_E(\Hom_{\Cr {\Gamma}}(T', T), M) \rightarrow \Hom_E(\Hom_{\Cr {\Gamma}}(T_1, T), M) \rightarrow \Hom_E(\Hom_{\Cr {\Gamma}}(T_2, T), M) \\
  &\rightarrow \Ext^1_E(\Hom_{\Cr {\Gamma}}(T', T), M) \rightarrow \Ext^1_E(\Hom_{\Cr {\Gamma}}(T_1, T), M) = 0 \end{align*}
  the last equality coming from the fact that $\Hom_{\Cr {\Gamma}}(T_1, T)$ is a projective $E$-module.
  
  Let us show that the morphism of functors from $\T$ to $\md E'^{\op}$
  \begin{align*}
   \Phi:\Hom_{\Cr {\Gamma}}(T', -) &\rightarrow \Hom_E(\Hom_{\Cr {\Gamma}}(-, T), \Hom_{\Cr {\Gamma}}(T', T)) \\
    \phi &\mapsto \Hom_{\Cr {\Gamma}}(\phi, T)
  \end{align*}
  is an isomorphism. By additivity and since $\T = \add(T)$, it is enough to look at $\Phi_T$. Let $\phi \in \Hom_{\Cr {\Gamma}}(T', T)$. Then $\Phi_T(\phi)(\id_T) = \phi$. Finally, $\Phi_T$ is injective. Moreover, if $\phi' \in \Hom_E(\Hom_{\Cr {\Gamma}}(T, T), \Hom_{\Cr {\Gamma}}(T', T))$ then $\Phi_T(\phi'(\id_T)) = \phi'$ so that $\Phi_T$ is surjective.
  
   By comparing the two previous long exact sequences, one gets the following isomorphisms of $E'^{op}$-modules.
  $$\End_E(M) \simeq \Hom_E(\Hom_{\Cr {\Gamma}}(T', T), M) \simeq \End_{\Cr {\Gamma}}(T') = E'^{\op}$$
   and
  $$\Ext^1_E(M, M) = \Ext^1_E(\Hom_{\Cr {\Gamma}}(T', T), M) = 0.$$
  
  Applying $\Hom_{\Cr {\Gamma}}(-, T)$ to (\ref{ses0}) yields the following long exact sequence:
  $$0 \rightarrow \Hom_{\Cr {\Gamma}}(T, T) \rightarrow \Hom_{\Cr {\Gamma}}(T'_1, T) \rightarrow \Hom_{\Cr {\Gamma}}(T'_2, T) \rightarrow \Ext^1_{\Cr {\Gamma}}(T, T) = 0.$$
  Hence, $M$ is a tilting $E$-module. From theorem \ref{bascde}, one deduces that $E$ and $E'$ are derived equivalent and that $\T$ and $\add(T')$ contain the same number of isomorphism classes of indecomposable objects. It is now obvious that $\add(T') = \T'$ (if $X \in \T' \setminus \add(T')$, the proof can be done replacing $T'$ by $T' \oplus X$ and we would obtain that $\T$, $\add(T' \oplus X)$ and $\add(T')$ contain the same number of isomorphism classes of indecomposable objects. This would be a contradiction). \cqfd 
\end{demo}

Recall this theorem of Igusa:

\begin{thm}[\citeb{3.2, b}{Ig90}]
 \label{thig}
 Let $A$ be a $k$-algebra of finite dimension and finite global dimension. Let $\phi$ be an automorphism of $A$ such that there exists a family of primitive idempotents of $A$ on which $\phi$ acts as a permutation. Then, the Gabriel quiver of $A$ has no arrows between any two vertices of the same orbit of $\phi$.
\end{thm}

Recall also the following theorem of Iyama, in a particular case:

\begin{thm}[\citeb{5.1 (3)}{Iy07-1}]
 \label{amasbasc}
 Assume that $\tilde \T' \in \Add(\Cr)$ is rigid and contains the projective-injective objects of $\Cr$. Let $\tilde T' \in \tilde \T'$ be basic such that $\tilde \T' = \add(\tilde T')$, $\tilde E' = \End_\Cr(\tilde T')$. Then $\gldim(\tilde E') \leq 3$ if and only if $\tilde \T'$ is cluster-tilting.
\end{thm}

We shall also need the following proposition of Bongartz:

\begin{prop}[\citeb{p. 463}{Bo83}]
 \label{bo}
 Let $Q$ be a quiver and $I$ be an admissible ideal of $kQ$ such that $kQ / I$ is finite dimensional. Denote by $J$ the Jacobson radical of $kQ$ (that is the ideal of $kQ$ generated by arrows). Let $i, j \in Q_0$. Then
 $$\dim e_j (I/(IJ + JI)) e_i = \dim \Ext^2_{\md kQ/I}(S_i, S_j)$$
 where $S_i$ and $S_j$ are the simple representations of $kQ$ supported on vertices $i$ and $j$.
\end{prop}

Finally, recall a particular case of a theorem by Lenzing:

\begin{thm}[\citeb{satz 5}{Le69}]
 \label{hl}
 If $A$ is a $k$-algebra of finite dimension and finite global dimension, then every nilpotent element of $A$ is in the additive subgroup $[A, A]$ generated by commutators.
\end{thm}

We can now state and prove the main result of this section.

\begin{thm}
 \label{recapend} 
 Suppose that there exists a category $\T \in \adds(\Cr {\Gamma})^{\md k[{\Gamma}]}$ which is maximal $\md k[{\Gamma}]$-stable rigid without $\md k[{\Gamma}]$-loops. 
 
 Let $\T' \in \adds(\Cr {\Gamma})^{\md k[{\Gamma}]}$ be maximal $\md k[{\Gamma}]$-stable rigid, $\tilde \T' = F\T'$. Let $T' \in \T'$,  $\tilde T' \in \tilde \T'$ be basic such that $\T' = \add(T')$ and $\tilde \T' = \add(\tilde T')$. Let $E' = \End_{\Cr {\Gamma}}(T')$ and  $\tilde E' = \End_\Cr(\tilde T')$. Then:
 \begin{enumerate}
  \item $\T'$ has no $\md k[{\Gamma}]$-loops; $\tilde \T'$ has no ${\Gamma}$-loops; \label{recapend1} 
  \item $\displaystyle \gldim(E') = \gldim(\tilde E') \left\{ \begin{array}{ll} = 3 & \quad \text{if } T' \text{ is not projective} \\ \leq 2 & \quad \text{else;} \end{array} \right.$ \label{recapend2} 
  \item $\T'$ and $\tilde \T'$ are cluster-tilting; \label{recapend3} 
  \item for every simple $E'$-modules $S$ and $S'$ such that $\add(k[{\Gamma}] \tens S) = \add(k[{\Gamma}] \tens S')$, one has $\Ext^1_{E'}(S, S') = \Ext^2_{E'}(S, S') = 0$; for every simple $\tilde E'$-modules $S$ and $S'$ such that $\add\left(\bigoplus_{g \in {\Gamma}} \g \tens S \right) = \add\left(\bigoplus_{g \in {\Gamma}} \g \tens S' \right)$, one has $\Ext^1_{\tilde E'}(S, S') = \Ext^2_{\tilde E'}(S, S') = 0$; \label{recapend4} 
  \item $\T'$ has no $\md k[{\Gamma}]$-$2$-cycles; $\tilde \T'$ has no ${\Gamma}$-$2$-cycles. \label{recapend5} 
 \end{enumerate}
\end{thm}

\begin{demo}
 \begin{enumerate}
  \item Using lemma \ref{loop}, it is enough to prove that $\tilde \T'$ has no ${\Gamma}$-loops. For $g \in {\Gamma}$, we show that $\tilde \T'$ has no $\langle g \rangle$-loops. For that, we show that $\tilde E'$ satisfies the hypothesis of theorem \ref{thig}. Using proposition \ref{deq}, $\tilde E'$ is of finite global dimension since $\tilde E$ is. From proposition \ref{relcycl}, there exists $X \in \Cr \langle g \rangle$ such that $FX \simeq \tilde T'$. It induces an action of $\langle g \rangle$ on $\tilde E' \simeq \uHom_{\Cr \langle g \rangle}(X, X)$. Using propositions \ref{relcycl} and \ref{actssgrp} implies that $FX$ can be split up into a direct summand of indecomposable objects such that $\langle g \rangle$ acts on it by permuting these objects. Then $\langle g \rangle$ acts on the family of primitive idempotents corresponding to these objects by permutation and therefore theorem \ref{thig} applies.
  \item This follows from (\ref{recapend1}) and proposition \ref{gldim}.
  \item This follows from (\ref{recapend2}) and theorem \ref{amasbasc}.
  \item As $\T'$ (resp. $\tilde \T'$) has no $\md k[{\Gamma}]$-loops (resp. ${\Gamma}$-loops), 
  $$\Ext^1_{E'}(S, S') = 0 \quad (\text{resp. }\Ext^1_{\tilde E'}(S, S') = 0).$$
  Concerning $\Ext^2$, if $X \in \T'$ is not projective, $S_X$ has the following projective resolution, given in proposition \ref{gldim}: 
  $$0 \rightarrow \Hom_{\Cr {\Gamma}}(X, T) \rightarrow \Hom_{\Cr {\Gamma}}(T'', T) \rightarrow \Hom_{\Cr {\Gamma}}(T', T) \rightarrow \Hom_{\Cr  {\Gamma}}(X, T) \rightarrow S_X \rightarrow 0.$$
  As $\add(T'') \inter \add(k[{\Gamma}] \tens X) = 0$, $\Hom_{E'}(\Hom_{\Cr {\Gamma}}(T'', T), S_{X'}) = 0$ and therefore 
  $$\Ext^2_{E'}(S_X, S_{X'}) = 0$$
  for every indecomposable object $X' \in \add(k[{\Gamma}] \tens X)$. 
  
  If $X$ is projective, the projective resolution is
  $$ 0 \rightarrow \Hom_{\Cr {\Gamma}}(Z,T) \rightarrow \Hom_{\Cr {\Gamma}}(T',T) \rightarrow \Hom_{\Cr {\Gamma}}(X,T) \rightarrow S_X \rightarrow 0.$$
  As $\add(k[{\Gamma}] \tens X) \inter \add(T') = 0$, $\add(k[{\Gamma}] \tens X) \inter \add(Z) = 0$ and by the same argument as before, $\Ext^2_{E'}(S_X, S_{X'}) = 0$ for any indecomposable object $X' \in \add(k[{\Gamma}] \tens X)$. The same argument works for $\tilde E'$.
  \item It is enough to show this for $\T'$, thanks to lemma \ref{loop}. Suppose that $\T'$ has a $\md k[{\Gamma}]$-$2$-cycle. As $\End_{\Cr {\Gamma}}(T') \simeq \Hom_{\md k[{\Gamma}]}(\un, \uHom_{\Cr {\Gamma}}(T', T'))$, there exist two arrows $\tilde a$ and $\tilde b$ in the Gabriel quiver of $\tilde T'$ such that $\tilde a \tilde b$ is a ${\Gamma}$-$2$-cycle of $\tilde \T'$ and $ab$ is a $\md k[{\Gamma}]$-$2$-cycle of $\T'$ with $a = \sum_{g \in {\Gamma}} g \cdot \tilde a$ and $b = \sum_{g \in {\Gamma}} g \cdot \tilde b$. As $a b$ is nilpotent, theorem \ref{hl} induces the following identity in $E'$:
   $$a b = \sum_{i = 1}^n \lambda_i [u_i, v_i]$$
  where, for each $i$, $u_i = \sum_{g \in {\Gamma}} g \cdot \tilde u_i$ and $v_i = \sum_{g \in {\Gamma}} g \cdot \tilde v_i$ where $\tilde u_i$ and $\tilde v_i$ are paths of the Gabriel quiver of $\tilde T'$. One can suppose that $\tilde u_1 = \tilde a$ and $\tilde v_1 = \tilde b$ and, without loss of generality that $\tilde a \tilde b$ and $\tilde b \tilde a$ do not appear as terms of $[\tilde u_i, \tilde v_i]$ for $i \geq 2$.
  
  If $\lambda_1 \neq 1$, 
  $$\tilde a \tilde b = e_{t(\tilde a)} \tilde a \tilde b e_{s(\tilde b)}= \frac{1}{1-\lambda_1} \left(-\lambda_1 e_{t(\tilde a)} \tilde b \tilde a e_{s(\tilde b)} + \sum_{i=2}^n e_{t(\tilde a)} \lambda_i [u_i, v_i] e_{s(\tilde b)} \right)$$
  is a non trivial identity in $\End_\Cr(FT')$ (that is, an identity which is false in the path algebra of the quiver). It contains $\tilde a \tilde b$ with a non zero coefficient in the left hand side.
  
  If $\lambda_1 = 1$, as $\tilde a \tilde b$ is a ${\Gamma}$-$2$-cycle, there exists $h \in {\Gamma}$ such that $t(\tilde a) = h \cdot s(\tilde b)$ and
  $$e_{h \cdot t(\tilde b)} (h \cdot \tilde b) \tilde a e_{s(\tilde a)}= -\sum_{i=2}^n \lambda_i e_{h \cdot t(\tilde b)} [u_i, v_i] e_{s(\tilde a)}$$
  is a non trivial identity in $\End_\Cr(FT')$. It contains $(h \cdot \tilde b)\tilde a$ with a non zero coefficient in the left hand side.
  
  In these two cases, one gets a non trivial relation, and as a consequence, using proposition \ref{bo}, $\Ext^2_\Cr(S_{s(\tilde b)}, S_{t(\tilde a)}) \neq 0$ or $\Ext^2_\Cr(S_{s(\tilde a)}, S_{t(h \cdot \tilde b)}) \neq 0$ which contradicts (\ref{recapend4}). \cqfd
 \end{enumerate}
\end{demo}

\begin{df}
 \label{eqkG}
 Let $X, Y \in \Cr {\Gamma}$. We write $X \sim Y$ if $\add(k[{\Gamma}] \tens X) = \add(k[{\Gamma}] \tens Y)$ and $X$ and $Y$ are said to be \emph{equivalent modulo $\md k[{\Gamma}]$}.
\end{df}

The following theorem summarizes the results concerning mutation:

\begin{thm}
 \label{recapmutcat}
 Suppose that there exists a category $\T \in \Add(\Cr {\Gamma})$ maximal $\md k[{\Gamma}]$-stable rigid which has no $\md k[{\Gamma}]$-loops. Let $\T' \in \Add(\Cr {\Gamma})$ be maximal $\md k[{\Gamma}]$-stable rigid. Then $\T'$ is cluster-tilting (hence maximal rigid). 
  
 Let $X \in \T'$ be an indecomposable object of $\Cr {\Gamma}$, and $\T'_0 \in \Add(\Cr {\Gamma})^{\md k[{\Gamma}]}$ satisfying $X \notin \T'_0$ and $\add(\T'_0, k[{\Gamma}] \tens X) = \T'$. If $X$ is projective, every $Y \in \Cr {\Gamma}$ which is indecomposable such that $\add(\T'_0, k[{\Gamma}] \tens Y)$ is maximal rigid is equivalent to $X$ modulo $\md k[{\Gamma}]$. If $X$ is not projective, there exists a unique $Y \in \Cr \Gamma$ such that $\{X, Y\}$ is an exchange pair associated with $\T'_0$. Moreover, in this case, $X$ and $Y$ are neighbours.
 
 If $X' \in \T'$ is equivalent to $X$ modulo $\md k[{\Gamma}]$ and if $\{X, Y\}$, $\{X', Y'\}$ are two exchange pairs associated with $\T'_0$, then $Y$ and $Y'$ are equivalent modulo $\md k[{\Gamma}]$. If $\bar X$ and $\bar Y$ denote the equivalence classes of $X$ and $Y$ modulo $\md k[{\Gamma}]$, one will denote
 $$\mu_{\bar X}(\T') = \add(\T'_0, k[{\Gamma}] \tens Y).$$
 
 Hence one has $\mu_{\bar X}(\mu_{\bar Y}(\T')) = \T'$.
\end{thm}

\begin{demo}
 This follows from proposition \ref{propdfmut}, lemmas \ref{psbclvois} and \ref{invmdG}, corollary \ref{mutvois} and theorem \ref{recapend}. Note that if $X$ is not projective, the existence and unicity of $Y$ such that $(X, Y)$ is a left exchange pair associated with $\T'_0$ is clear by proposition \ref{propdfmut}. Using lemma \ref{psbclvois}, one deduces that $X$ and $Y$ are neighbours. Therefore, corollary \ref{mutvois} implies that $\{X, Y\}$ is an (unordered) exchange pair. \cqfd
\end{demo}

\subsection{Exchange matrices}

\label{echange}
As in the previous section, indecomposable projective objects of $\Cr$ are supposed to have left rigid quasi-approximations. As before, all results remain valid if they have right rigid quasi-approximations.

As before, $\T \in \adds(\Cr {\Gamma})^{\md k[{\Gamma}]}$ is maximal rigid $\md k[{\Gamma}]$-stable (and finitely generated). One supposes moreover that $\T$ has no $\md k[{\Gamma}]$-loops. Let $\tilde \T = F \T$. Thanks to the previous section, $\T$ and $\tilde \T$ are cluster-tilting. Let $T \in \Cr {\Gamma}$ and $\tilde T \in \Cr$ be basic such that $\T = \add(T)$ and $\tilde \T = \add(\tilde T)$.

One denotes by $Q$ the Gabriel quiver of $\End_{\Cr {\Gamma}}(T)$ and by $\tilde Q$ the quiver of $\End_\Cr(\tilde T)$. Denote by $Q_0 / \md k[{\Gamma}]$ the set of equivalence classes modulo $\md k[{\Gamma}]$ as in definition \ref{eqkG}. There is a canonical bijection between the sets $Q_0 / \md k[{\Gamma}]$ and $\tilde Q_0 / {\Gamma}$. If $\bar X \in \tilde Q_0 / {\Gamma}$, $\bar X^\circ \in Q_0 / \md k[{\Gamma}]$ will denote the image of $\bar X$ by this bijection. Denote by $P \subset Q_0$ and $\tilde P \subset \tilde Q_0$ the sets of vertices corresponding to projective objects.

\begin{df} \label{dfmate}
 If $\bar X \in (\tilde Q_0 \setminus \tilde P)/{\Gamma}$ and $\bar Z \in \tilde Q_0 / {\Gamma}$, put
 $$b_{\bar Z \bar X} = \frac{\#\{q \in \tilde Q_1 \,|\, s(q) \in \bar Z, t(q) \in \bar X\} - \#\{q \in \tilde Q_1 \,|\, s(q) \in \bar X, t(q) \in \bar Z\}}{\# \bar X}$$
 Denote by $B(\T)$ the matrix having these entries. It will be called the \emph{exchange matrix of $\T$}.
\end{df}

\begin{rems}
 \begin{itemize}
  \item As $\tilde \T$ has no ${\Gamma}$-$2$-cycles, 
  $$\#\{q \in \tilde Q_1 \,|\, s(q) \in \bar X, t(q) \in \bar Z\} = 0 \quad \text{or} \quad \#\{q \in \tilde Q_1 \,|\, s(q) \in \bar Z, t(q) \in \bar X\} = 0.$$ 
  \item As $\tilde \T$ has no ${\Gamma}$-loops, $b_{\bar X \bar X} = 0$. 
  \item It is easy to see that
   $$b_{\bar Z \bar X} = \#\{q \in \tilde Q_1 \,|\, s(q) \in \bar Z, t(q) = X\} - \#\{q \in \tilde Q_1 \,|\, s(q) = X, t(q) \in \bar Z\}$$
   for every $X \in \bar X$, hence $B(\T)$ has integer coefficients.
  \item The exchange matrix is clearly skew-symmetrizable (right multiplication by the diagonal matrix $(\# \bar X)_{\bar X \in \left(\tilde Q_0 \setminus \tilde P\right)/{\Gamma}}$ gives a skew-symmetric matrix).
 \end{itemize} 
\end{rems}

\begin{rem}
 The matrix of definition \ref{dfmate} coincides with the exchange matrix of \cite{FuKe} and \cite{GeLeSc06-1}.
\end{rem}

Fix now the three following mutations and matrices:
\begin{itemize}
 \item $\mu$ on $\Cr {\Gamma}$ and $B(\T)$ are defined as before;
 \item $\tilde \mu$ is the mutation on $\Cr \simeq \Cr\{e\}$ defined as before replacing $\Gamma$ by the trivial group $\{e\}$ and $\tilde B(\T)$ is the corresponding exchange matrix;
 \item $\mu^\circ$ is the mutation on $\Cr {\Gamma} \simeq (\Cr \Gamma) \{e\}$ defined as before replacing $\Cr$ by $\Cr \Gamma$ and $\Gamma$ by $\{e\}$ and $B^\circ(\T)$ is the corresponding exchange matrix.
\end{itemize}

\begin{rem}
 The definitions of $\tilde \mu$ and $\mu^\circ$ coincide with those of \cite{FuKe} and \cite{GeLeSc06-1}.
\end{rem}

\begin{prop}
 \label{relmut}
 Let $\bar X \in (\tilde Q_0 \setminus \tilde P)/{\Gamma}$ and $\bar Z \in \tilde Q_0 / {\Gamma}$. Then
 \begin{enumerate}
  \item $\displaystyle F \mu_{\bar X}(\T) = \left(\prod_{X' \in \bar X} \tilde \mu_{X'}\right)(\tilde \T)$ where the $\tilde \mu_{X'}$ commute. \label{relmut1}
  \item For every $X \in \bar X$, $\displaystyle b_{\bar Z \bar X} = \sum_{Z \in \bar Z} \tilde b_{ZX}$. \label{relmut2}
  \item $\displaystyle \mu_{\bar X}(\T) = \left(\prod_{X' \in \bar X^\circ} \mu^\circ_{X'}\right)(\T)$ where the $\mu^\circ_{X'}$ commute. \label{relmut3}
  \item For every $Z \in \bar Z^\circ$, $\displaystyle b_{\bar Z \bar X} = \frac{\# Z}{\ell(Z)} \sum_{X \in \bar X^\circ} \frac{\ell(X)}{\# X} b^\circ_{ZX}$. \label{relmut4}
 \end{enumerate}
\end{prop}

\begin{demo}
 \begin{enumerate}
  \item Let $0 \rightarrow X \rightarrow T \rightarrow Y \rightarrow 0$ be the admissible short exact sequence of $\Cr {\Gamma}$ corresponding to the mutation $\mu_X$ in $\T$. Then $X$ and $Y$ are neighbours and, as a consequence, one can write $\bar X = \{X_i \,|\, i \in \llbracket 1, \# X \rrbracket \}$ and $\bar Y = \{Y_i\,|\, i \in \llbracket 1, \# Y \rrbracket\}$ in such a way that for every $i, j \in \llbracket 1, \ell(X) \rrbracket = \llbracket 1, \ell(Y) \rrbracket$, $\dim \Ext^1_\Cr(X_i, Y_i) = \delta_{ij}$ by using lemma \ref{neicr}. Thus, for $i \in  \llbracket 1, \ell(X) \rrbracket$, there exists a non split admissible short exact sequence $0 \rightarrow X_i \rightarrow T_i \rightarrow Y_i \rightarrow 0$ in $\Cr$. As $\tilde \T$ has no ${\Gamma}$-loops, none of the $X_j$ and $Y_j$ is in $\add(T_i)$ and finally, the result is clear. 
  \item This is an easy consequence of the definition.
  \item The proof is the same as for (\ref{relmut1}).
  \item Let $X \in \bar X$ and $\T_0 \in \Add(\Cr {\Gamma})^{\md k[{\Gamma}]}$ be such that $\add(X[{\Gamma}]) \inter F\T_0 = 0$ and $\T = \add(\T_0,  X[{\Gamma}])$. Let $0 \rightarrow X \rightarrow T \rightarrow Y \rightarrow 0$ be an admissible short exact sequence in $\Cr$  corresponding to the mutation $\tilde \mu_X$. As $\tilde \T$ has no ${\Gamma}$-loops, $T \in \add(F\T_0)$. As a consequence, lemma \ref{Fminapprox} gives
  $$0 \rightarrow X[{\Gamma}] \rightarrow T[{\Gamma}] \rightarrow Y[{\Gamma}] \rightarrow 0 \simeq \bigoplus_{X' \in \bar X^\circ} \left[0 \rightarrow X' \rightarrow T'_{X'} \rightarrow Y'_{X'} \rightarrow 0\right]^{\ell(X') / \# X'}$$
  where, for every $X' \in \bar X^\circ$, $0 \rightarrow X' \rightarrow T'_{X'} \rightarrow Y'_{X'} \rightarrow 0$ is an admissible short exact sequence corresponding to $\mu^\circ_{X'}$ at $\T$, and the exponents $\ell(X') / \# X'$ come from lemma \ref{nbiso}. Let now $Z \in \bar Z^\circ$. By definition of $B^\circ(T)$, the number of copies of $Z$ in the middle term of the right hand side is
  $$\sum_{X' \in \bar X^\circ} \frac{\ell(X')}{\# X'} \max(0, -b^\circ_{ZX'}) = \max\left(0, -\sum_{X' \in \bar X^\circ} \frac{\ell(X')}{\# X'} b^\circ_{ZX'}\right)$$
  the equality coming from the fact that all $b^\circ_{ZX}$ have the same sign, because $\T$ has no $\md k[{\Gamma}]$-$2$-cycles.
  
  Moreover, the number of copies of $Z$ in the middle term of the left hand side is $\ell(Z)/\# Z$ times the number of copies of an element of $\bar Z$ in $T$, that is
  $$\frac{\ell(Z)}{\# Z} \sum_{Z' \in \bar Z} \max(0, -\tilde b_{Z'X}) = \max \left(0, -\frac{\ell(Z)}{\# Z} \sum_{Z' \in \bar Z} \tilde b_{Z'X}\right) = \max \left(0, -\frac{\ell(Z)}{\# Z} b_{\bar Z \bar X}\right).$$
  
  One deduces the identity:
  $$\max\left(0, -\sum_{X' \in \bar X^\circ} \frac{\ell(X')}{\# X'} b^\circ_{ZX'}\right) = \max \left(0, -\frac{\ell(Z)}{\# Z} b_{\bar Z \bar X}\right)$$
  and using the same argument for the admissible short exact sequence corresponding to the mutation from $Y$ to $X$, one deduces that
  $$\max\left(0, \sum_{X' \in \bar X^\circ} \frac{\ell(X')}{\# X'} b^\circ_{ZX'}\right) = \max \left(0, \frac{\ell(Z)}{\# Z} b_{\bar Z \bar X}\right).$$
  
 These two equalities yield the result. \cqfd
 \end{enumerate}
\end{demo}

\begin{rem}
 The formula (\ref{relmut2}) was given in Dynkin cases by Dupont in \cite{Du08} (see also erratum \cite{Du}). On the other hand, (\ref{relmut4}) is a generalization of the formula given by Yang in \cite{Ya09} for cluster algebras of finite type (for a cyclic group).
\end{rem}

\begin{thm}\label{cfmutfz}
 Let $X \in \T$ be an indecomposable non projective object in $\Cr {\Gamma}$. Then
 $$B\left(\mu_X(\T)\right) = \mu_{\bar X^\circ} \left(B(\T)\right)$$
 where $\bar X$ is the orbit of $X$ in $Q_0/ \md k[{\Gamma}]$ and $\mu_{\bar X}$ the mutation of matrices defined by Fomin and Zelevinsky.
\end{thm}

\begin{demo}
 The result is known for $\tilde B(T)$ and $\tilde \mu$ (the proof is similar as the one in \cite[\S 14]{GeLeSc06-1} for example). Let $\bar X = \{X_1, X_2, \dots X_n\}$. For $i \in \llbracket 1, n \rrbracket$, denote by ${^i \tilde b}_{ZY}$ the coefficients of the matrix
 $$\tilde B \left(\tilde \mu_{X_i} \tilde \mu_{X_{i-1}} \dots \tilde \mu_{X_1}(\tilde \T)\right).$$
 Then, by an easy induction, by using the fact that $\tilde \T$ has no ${\Gamma}$-loops nor ${\Gamma}$-$2$-cycles, one has
 \begin{itemize}
  \item $\displaystyle ^i \tilde b_{ZY} = - \tilde b_{ZY}$ if $Z = X_j$ or $Y = X_j$ where $1 \leq j \leq i$;
  \item $\displaystyle ^i \tilde b_{ZY} = \tilde b_{ZY}$ if $Z = X_j$ or $Y = X_j$ where $i < j \leq n$; 
  \item $\displaystyle ^i \tilde b_{ZY} = \tilde b_{ZY} + \sum_{j = 1}^i \frac{|\tilde b_{ZX_j}| \tilde b_{X_j Y} + \tilde b_{ZX_j}|\tilde b_{X_j Y}|}{2}$ else.
 \end{itemize}
 
 Denoting by $^* b$ the coefficients of $B\left(\mu_X(\T)\right)$, by using proposition \ref{relmut}, for $\bar Y \in (\tilde Q_0 \setminus \tilde P)/{\Gamma}$ and $Z \in \tilde Q_0 / {\Gamma}$, for $Y \in \bar Y$,
 \begin{itemize}
 \item $\displaystyle ^* b_{\bar Z \bar Y} = \sum_{Z \in \bar Z} {^n \tilde b}_{Z Y} = \sum_{Z \in \bar Z} (-\tilde b_{ZY}) = -b_{\bar Z \bar Y}$ if $\bar Z = \bar X$ or $\bar Y = \bar X$;
 \item $\displaystyle ^* b_{\bar Z \bar Y} = \sum_{Z \in \bar Z} {^n \tilde b}_{Z Y} = \sum_{Z \in \bar Z} \left[\tilde b_{ZY} + \sum_{j = 1}^n \frac{|\tilde b_{ZX_j}| \tilde b_{X_j Y} + \tilde b_{ZX_j}|\tilde b_{X_j Y}|}{2}\right]$ else. Then, using the fact that all $\tilde b_{Z X_j}$ (resp. all $\tilde b_{X_j Y}$) are of the same sign (as $\tilde \T$ has no ${\Gamma}$-$2$-cycles),
   \begin{align*}
    \sum_{Z \in \bar Z} \sum_{j = 1}^n \frac{|\tilde b_{ZX_j}| \tilde b_{X_j Y} + \tilde b_{ZX_j}|\tilde b_{X_j Y}|}{2} &= \sum_{j = 1}^n \frac{\left|\sum_{Z \in \bar Z} \tilde b_{ZX_j}\right| \tilde b_{X_j Y} + \left(\sum_{Z \in \bar Z} \tilde b_{ZX_j}\right)|\tilde b_{X_j Y}|}{2} \\ &= \sum_{j = 1}^n \frac{|b_{\bar Z \bar X}| \tilde b_{X_j Y} + b_{\bar Z \bar X} |\tilde b_{X_j Y}|}{2} \\ &= \frac{|b_{\bar Z \bar X}| \left( \sum_{j=1}^n \tilde b_{X_j Y}\right) + b_{\bar Z \bar X} \left|\sum_{j=1}^n \tilde b_{X_j Y}\right|}{2} \\ &= \frac{|b_{\bar Z \bar X}| b_{\bar X \bar Y} + b_{\bar Z \bar X} |b_{\bar X \bar Y}|}{2}
   \end{align*}
 \end{itemize}
 and the result is proved. \cqfd
\end{demo}

\subsection{Cluster characters}

\label{caramas}
The ideas of this section generalize results of \cite{FuKe}. 

 We retain notation and hypothesis of previous sections. One supposes moreover that there exists $\T \in \Add(\Cr)^{\Gamma}$ maximal ${\Gamma}$-stable rigid without ${\Gamma}$-loops. Thus, $\T$ is cluster-tilting using theorem \ref{recapmutcat}. One denotes by $T_1, T_2, \dots, T_n$ the indecomposable objects of $\T$ up to isomorphism, the $T_i$ for $i \in \llbracket r+1, n \rrbracket$ being the projective objects. The action of ${\Gamma}$ on $\T$ induces an action on $\llbracket 1, n \rrbracket$. If $i \in \llbracket 1, n \rrbracket$, $\ig$ denotes its equivalence class modulo ${\Gamma}$. One denotes by $\bar I$ the set of these equivalence classes. Let $T = T_1 \oplus T_2 \oplus \dots \oplus T_n$ and $E = \End_\Cr(T)$. For $i \in \llbracket 1, n \rrbracket$, $S_i$ denotes the simple $E$-module corresponding to $T_i$ (the head of $\Hom_\Cr(T_i, T)$). 

 Thus, all hypothesis to apply results of \cite{FuKe} hold.

\begin{nt}[\cite{FuKe}]
 For $L, N \in \md E$, let
 $$\langle L, N \rangle_\tau = \dim_k \Hom_B(L, N) - \dim_k \Ext^1_B(L, N)\;;$$
 $$\langle L, N \rangle_3 = \sum_{i = 0}^3 (-1)^i \dim_k \Ext^i_B(L, N).$$
\end{nt}

\begin{rem}
 One keeps the notation $\langle , \rangle_3$ introduced in \cite{FuKe}, but, here, as the global dimension of $E$ is less than $3$, this form is the Euler form (see also \cite[remark 2.4]{FuKe}).
\end{rem}

\begin{nt}
 If $L \in \md E$, $\udim L$ denotes its image in $K_0(E)$, that is its dimension vector relatively to the $S_i$.
\end{nt}

\begin{prop}[\citeb{proposition 2.1}{FuKe}]
 If $L, N \in \md E$, then $\langle L, N \rangle_3$ depends only on $\udim L$ and $N$.
\end{prop}

\begin{nt}
 Using previous proposition, if $L, N \in \md E$, one can set
 $$\langle \udim L, N \rangle_3 = \langle L, N \rangle_3.$$
\end{nt}

Define $\pi$ to be the following canonical projection:
\begin{align*}
  \pi: \Q\left[x_i^{\pm 1}\right]_{i \in \llbracket 1, n \rrbracket} &\rightarrow \Q\left[x_{\ig}^{\pm 1}\right]_{\ig \in \bar I} \\
  x_i &\mapsto x_{\ig}.
 \end{align*}

\begin{df}
 For $X \in \Cr$, define the Laurent polynomial $P_X$ of $\Q\left[x_{\ig}^{\pm 1}\right]_{\ig \in \bar I}$ by $P_X = \pi(X'_X)$ where $X'_X$ is the Laurent polynomial of $\Q\left[x_i^{\pm 1}\right]_{i \in I}$ defined by Fu and Keller in \cite{FuKe}. In other words, 
 \begin{equation} P_X = \left( \prod_{i = 1}^n x_{\ig}^{\langle \Hom_\Cr(T, X), S_i \rangle_\tau} \right) \sum_{e \in \N^{\llbracket 1, n \rrbracket}} \left(\chi(\Gr_e(\Ext^1_\Cr(T, X))) \prod_{i = 1}^n x_{\ig}^{-\langle e, S_i \rangle_3}\right),\label{defpx}\end{equation}
 where $\Gr_e(\Ext^1_\Cr(T, X))$ is the variety of $E$-submodules $B$ of $\Ext^1_\Cr(T, X)$ such that $\udim B = e$ and $\chi$ is the Euler characteristic with respect to \'etale cohomology with proper support.
\end{df}

\begin{lem} \label{indeppx}
 The Laurent polynomial $P_X$ depends only on the class of $X$ modulo ${\Gamma}$.
\end{lem}

\begin{demo}
 As $T$ is ${\Gamma}$-invariant, for every $g \in {\Gamma}$, 
 $$\Hom_\Cr(T, \g \tens X) = \Hom_\Cr(\g \tens T, \g \tens X) = \g^{-1} \tens \Hom_\Cr(T, X)$$
 where $\md E$ is canonically endowed with the action of ${\Gamma}$ induced by the action of ${\Gamma}$ on $\Cr$. As a consequence,
 $$\langle \Hom_\Cr(T, \g \tens X), S_i \rangle_\tau = \langle \g^{-1} \tens \Hom_\Cr(T, X), S_i \rangle_\tau = \langle  \Hom_\Cr(T, X), \g \tens S_i \rangle_\tau$$
 which leads to the conclusion concerning the first factor of the right-hand side of (\ref{defpx}). In the same way, one gets
 $$\Gr_e(\Ext^1_\Cr(T, \g \tens X)) \simeq \Gr_{g\cdot e}(\Ext^1_\Cr(T, X))$$
 and
 $$\langle e, S_i \rangle_3 = \langle g \cdot e, \g \tens S_i \rangle_3$$
 which yields the conclusion concerning the second factor of the right-hand side of (\ref{defpx}). \cqfd
\end{demo}

 According to lemma \ref{indeppx}, it makes sens to denote
 $$P_{\bar X} = P_X$$
 where $\bar X$ is the class of $X$ modulo ${\Gamma}$.

Here is the analogous of theorem \cite[theorem 2.2]{FuKe}:

\begin{thm}
 \label{fmut1}
 \begin{enumerate}
  \item For $\ig \in \bar I$, $P_{\bar T_i} = x_{\ig}$.
  \item If $X, Y \in \Cr$, $P_{\bar X \oplus \bar Y} = P_{\bar X} P_{\bar Y}$.
  \item If $X, Y \in \Cr$ and $\dim \Ext^1_\Cr(X, Y) = 1$, and if one fixes two non split admissible short exact sequences
  $$0 \rightarrow X \rightarrow Z \rightarrow Y \rightarrow 0 \quad \text{and} \quad 0 \rightarrow Y \rightarrow Z' \rightarrow X \rightarrow 0$$
  then $P_{\bar X} P_{\bar Y} = P_{\bar Z} + P_{\bar Z'}$.
 \end{enumerate}
\end{thm}

\begin{demo}
 This follows from \cite[theorem 2.2]{FuKe} by applying the ring morphism $\pi$. \cqfd
\end{demo}

\begin{cor} \label{pcfmutfz}
 The $P_{\bar X}$ satisfy the mutation formulas of Fomin and Zelevinsky encoded by the exchange matrices $B$ of definition \ref{dfmate}. In other words, if $\T'_0 \in \Add(\Cr {\Gamma})^{\md k[{\Gamma}]}$ and $X, Y \in \Cr$ are indecomposable objects such that
 \begin{itemize}
  \item $X, Y \notin \T'_0$,
  \item $\add(\T'_0, X[{\Gamma}])$ is maximal $\md k[{\Gamma}]$-stable rigid,
  \item $\mu_{\bar X}(\T') = \add(\T'_0, Y[{\Gamma}])$ where $\T' = \add(\T'_0, X[{\Gamma}])$,  
 \end{itemize}
 then
 \begin{equation}P_{\bar X} P_{\bar Y} = \prod_{\ig \in \bar I\,|\, B(\T')_{\bar T'_i \bar X} < 0} P_{\bar T'_i}^{-B(\T')_{\bar T'_i \bar X}} + \prod_{\ig \in \bar I\,|\, B(\T')_{\bar T'_i \bar X} > 0} P_{\bar T'_i}^{B(\T')_{\bar T'_i \bar X}}.\label{mutpx}\end{equation}
\end{cor}

\begin{demo}
 There exists $X_0 \in X[{\Gamma}]$ and $Y_0 \in Y[{\Gamma}]$ which are neighbours. Hence, one can suppose that $\Ext^1_\Cr(X, Y) = 1$ up to replacing $Y$ by a different representative of its equivalence class modulo ${\Gamma}$. Let 
 $$0 \rightarrow X \xrightarrow{f} Z \xrightarrow{g} Y \rightarrow 0$$
 be a non split admissible short exact sequence. In this case, $f$ is a minimal left $\T'_0$-approximation and $g$ is a minimal right $\T'_0$-approximation using proposition \ref{neighapp}. Hence, for $i \in \llbracket 1, n \rrbracket$, the number of $T'_i$ appearing in $Z$ is
  $$\#\{q \in Q'_1 \,|\, s(q) = X, t(q) = T'_i\}$$
  where $Q'$ is the Auslander-Reiten quiver of $\T'$. Thus, for $i \in \bar I$, the number of $T'_i$ with $i \in \ig$ which appear in $Z$ is
  $$\#\{q \in Q'_1 \,|\, s(q) = X, t(q) \in \bar T'_i\}$$
  and, as $\T'$ has no ${\Gamma}$-loops, if this number is strictly positive, it is equal by definition to $-B(\T')_{\bar T'_i \bar X}$, which permits to conclude. The second term of the right-hand side of (\ref{mutpx}) can be handled in the same way. \cqfd
\end{demo}

Denote by $\Ar(\Cr, {\Gamma}, \T)$ the subalgebra of $\Q(x_{\ig})_{\ig \in \bar I}$ generated by the $P_{\bar X}$ where $\bar X$ goes over all ${\Gamma}$-orbits of objects of $\Cr$ such that $\bigoplus_{X \in \bar X} X$ is rigid. Denote by $\Ar(\Cr, \T)$ the subalgebra of $\Q(x_i)_{i \in \llbracket 1, n \rrbracket}$ generated by the $X'_X$ where $X$ goes over the rigid objects of $\Cr$.

Denote by $\Ar_0(\Cr, {\Gamma}, \T)$ the subalgebra of $\Q(x_{\ig})_{\ig \in \bar I}$ generated by the $P_{\bar X}$ where $\bar X$ goes over the ${\Gamma}$-orbits of objects of $\Cr$. Denote by $\Ar_0(\Cr, \T)$ the subalgebra of $\Q(x_i)_{i \in \llbracket 1, n \rrbracket}$ generated by the $P'_X$ where $X$ goes over $\Cr$.

\begin{cor} \label{qalga}
 There is a commutative diagram of inclusions
 $$\xymatrix{
  \Ar(B(\T)) \ar@{^{(}->}[r] \ar@{^{(}->}[d] & \pi\left(\Ar(\tilde B(\T))\right) \ar@{^{(}->}[d] \\
  \Ar(\Cr, {\Gamma}, \T) \ar@{^{(}->}[r] \ar@{^{(}->}[d] & \pi\left(\Ar(\Cr, \T)\right) \ar@{^{(}->}[d] \\
  \Ar_0(\Cr, {\Gamma}, \T) \ar@{=}[r] & \pi\left(\Ar_0(\Cr, \T)\right).
 }$$
\end{cor}

\begin{demo}
 First of all, the inclusions $\Ar(B(\T)) \subset \Ar(\Cr, {\Gamma}, \T)$ and $\Ar(\tilde B(\T)) \subset \Ar(\Cr, \T)$ come from corollary \ref{pcfmutfz}. The bottom equality is clear using the definition of $P_{\bar X}$ in terms of $X'_X$. The horizontal middle inclusion comes from the fact that for all $P_{\bar X}$ where $\bigoplus_{X \in \bar X} X$ is rigid, $X \in \bar X$ is rigid and therefore $P'_X \in \Ar(\Cr, \T)$. The upper horizontal inclusion comes from the fact that if $P_{\bar X}$ and $P_{\bar Y}$ are linked by a string of mutations in $\Ar(\Cr, {\Gamma}, \T)$, then $P'_X$ and $P'_Y$ are also in $\Ar(\Cr, \T)$ according to proposition \ref{relmut}. \cqfd
\end{demo}

\begin{rem}
 In general, it seems to be a difficult problem to understand which of the inclusions in the previous diagram are isomorphisms.
\end{rem}

Let $\Delta$ (resp. $\Delta^\circ$, $\tilde \Delta$) be the non-oriented version of the graph whose adjacency matrix is the upper square submatrix of $B(\T)$ (resp. $B^\circ(\T)$, $\tilde B(\T)$). According to proposition \ref{relmut}, $\tilde \Delta$ and $\Delta$ are related by a classical folding process. On the other hand, $\Delta^\circ$ and $\Delta$ are related by a folding process deformed by some positive integer coefficients (the $n_X = \ell(X) / \# X$).

\begin{lem} \label{eqDyn}
 The following are equivalent:
 \begin{enumerate}
  \item $\Delta$ is a Dynkin diagram; \label{eqDyn1}
  \item $\Delta^\circ$ is a Dynkin diagram. \label{eqDyn2}
 \end{enumerate}
 Moreover, under these assumptions, $\tilde \Delta$ is also a Dynkin diagram.
\end{lem}

\begin{demo}
 Every diagram can be supposed to be connected without loss of generality. The proof that (\ref{eqDyn2}) implies (\ref{eqDyn1}) and that $\tilde \Delta$ is a Dynkin diagram can be done by a finite number of computations (see table of page \pageref{subqj}). 

 Let us now show that (\ref{eqDyn1}) implies (\ref{eqDyn2}). Suppose that $\Delta$ is a Dynkin diagram. Let us call critical point of $\Delta$ or $\tilde \Delta$ every vertex of valuation at least $3$ or every non simple edge. As $\Delta$ is a Dynkin diagram, it has at most one critical point. If $\tilde \Delta$ has a cycle, then it induces in $\Delta$ a cycle or at least to critical points, hence $\tilde \Delta$ is a tree. Moreover, an orbit under ${\Gamma}$ of critical points in $\tilde \Delta$ yields a critical point in $\Delta$. Thus, there is at most such an orbit. Suppose that there is two distinct points $A$ and $B$ in this orbit. As $\tilde \Delta$ is a tree, there is a unique shortest path between these two points, which is ``folded'' by the element of ${\Gamma}$ which sends $A$ on $B$. Hence it is easy to see that the middle of this path gives rise to a loop or a second critical point in $\Delta$, which is not possible. Finally, $\tilde \Delta$ has at most one critical point, which leads to the conclusion using a case by case proof. \cqfd
\end{demo}

\begin{prop}
 \label{typefini}
 The following are equivalent:
 \begin{enumerate}
  \item $\Ar(B(\T))$ has a finite number of cluster variables; \label{typefini1}
  \item $\Ar(\tilde B(\T))$ has a finite number of cluster variables; \label{typefini2}
  \item $\Ar(B^\circ(\T))$ has a finite number of cluster variables; \label{typefini3}
  \item $\Cr$ has a finite number of isomorphism classes of rigid indecomposable objects. \label{typefini4}
 \end{enumerate}
\end{prop}

\begin{demo}
 First of all, it is clear that (\ref{typefini2}) or (\ref{typefini3}) imply (\ref{typefini1}) using proposition \ref{relmut}. Now, if (\ref{typefini1}) is true, using the results of \cite{FoZe03}, $B(\T)$ is mutation-equivalent to a matrix, whose principal square submatrix encodes a Dynkin diagram. Up to changing $\T$, one can suppose that the principal square submatrix of $B(\T)$ encodes a Dynkin diagram. By lemma \ref{eqDyn}, the principal square submatrices of $\tilde B(\T)$ and $B^\circ(\T)$ encode Dynkin diagrams, which leads to (\ref{typefini2}) and (\ref{typefini3}) using again \cite{FoZe03}. It is clear that (\ref{typefini4}) implies the three others. If (\ref{typefini2}) is satisfied, the principal square submatrix of $\tilde B(\T))$ can be supposed to encode a Dynkin diagram. As a consequence, using a theorem of Keller and Reiten \cite{KeRe08}, the stable category of $\Cr$ is equivalent to a cluster category, which yields (\ref{typefini4}). \cqfd
\end{demo}

\subsection{Linear independence of cluster monomials}

One retains the notation of section \ref{caramas}. Again, this section generalizes results of \cite{FuKe}.

\begin{df}
 Two objects $X$, $Y$ of $\Cr$ are said to be \emph{congruent modulo ${\Gamma}$} if there exists two decompositions into indecomposable direct summands $X = X_1 \oplus X_2 \oplus \dots \oplus X_m$ and $Y = Y_1 \oplus Y_2 \oplus \dots \oplus Y_m$ such that for all $i \in \llbracket 1, m \rrbracket$, $X_i$ and $Y_i$ are equivalent modulo ${\Gamma}$ (that is, if there exists $g \in {\Gamma}$ such that $\g \tens X_i \simeq Y_i$). Equivalently, $X$ and $Y$ are congruent modulo ${\Gamma}$ if
 $$\bigoplus_{g \in {\Gamma}} \g \tens X \simeq \bigoplus_{g \in {\Gamma}} \g \tens Y.$$
\end{df}

If $X \in \Cr$, there exists an admissible short exact sequence
$$0 \rightarrow T^1_X \rightarrow T^2_X \xrightarrow{f} X \rightarrow 0$$
where $f$ is a minimal right $\T$-approximation. Using lemma \ref{cokapp} and the maximality of $\T$, one gets that $T^1_X \in \T$. Following \cite{Pa08} and \cite{FuKe}, write
$$\ind_\T(X) = [T^0_X] - [T^1_X] \in K_0(\T).$$
Denote also by $\ind'_\T(X)$ the image of $\ind_\T(X)$ in $K_0(\T) / {\Gamma}$.

The following lemma generalizes \cite[lemma 2.1]{DeKe08}:
\begin{lem}
\label{caracgeo}
 If $X$ is rigid and
 $$0 \rightarrow T^1_X \xrightarrow{h} T^2_X \xrightarrow{f} X \rightarrow 0$$
 is the previous admissible short exact sequence, then the orbit of $h$ under the action of $\Aut_\Cr(T^1_X) \times \Aut_\Cr(T^2_X)$ is a dense open subset of $\Hom_\Cr(T^1_X, T^2_X)$.
\end{lem}

\begin{demo}
 Let $h': T^1_X \rightarrow T^2_X$ be a morphism. Applying $\Hom_\Cr(-, X)$ to the admissible short exact sequence leads to the following long exact sequence:
 $$0 \rightarrow \Hom_\Cr(X, X) \rightarrow \Hom_\Cr(T^2_X, X) \rightarrow \Hom_\Cr(T^1_X, X) \rightarrow \Ext^1_\Cr(X, X) = 0$$
 which shows that there exists $\alpha: T^2_X \rightarrow X$ such that $\alpha h = f h'$. Moreover, as $f$ is a right  $\T$-approximation, there exists $\beta: T^2_X \rightarrow T^2_X$ such that $\alpha = f \beta$. Hence $f(\beta h - h') = 0$ and therefore, as $h$ is a kernel of $f$, there exists ${\Gamma}: T^1_X \rightarrow T^1_X$ such that $\beta h - h' = h {\Gamma}$. In other words, $h'$ is in the image of
 \begin{align*}
  \End_\Cr(T^1_X) \times \End_\Cr(T^2_X) &\rightarrow \Hom_\Cr(T^1_X, T^2_X) \\
  ({\Gamma}, \beta) &\mapsto \beta h - h {\Gamma}
 \end{align*}
 which is the differential of the application 
 \begin{align*}
  \Aut_\Cr(T^1_X) \times \Aut_\Cr(T^2_X) &\rightarrow \Hom_\Cr(T^1_X, T^2_X) \\
  (g_1, g_2) &\mapsto g_2 h g_1^{-1}.
 \end{align*}
 Here, we use the two identifications
  $$\Lie(\Aut_\Cr(T^1_X) \times \Aut_\Cr(T^2_X)) \simeq \End_\Cr(T^1_X) \times \End_\Cr(T^2_X)$$
   and
   $$\Lie(\Hom_\Cr(T^1_X, T^2_X)) \simeq \Hom_\Cr(T^1_X, T^2_X).$$
   Hence, one gets the conclusion. \cqfd
\end{demo}

The following lemma is inspired from \cite[lemma 2.2]{DeKe08}:
\begin{lem} \label{intvide}
 If $X$ is rigid and
 $$0 \rightarrow T^1_X \xrightarrow{h} T^2_X \xrightarrow{f} X \rightarrow 0$$
 is the admissible short exact sequence defined as before, then $T^1_X$ and $T^2_X$ have no common direct summand.
\end{lem}

\begin{demo}
 Suppose that $T^1_X \simeq T^0 \oplus \tilde T^1_X$ and $T^2_X \simeq T^0 \oplus \tilde T^2_X$ with $T^0 \neq 0$. Using previous lemma, the orbit of $h$ under the action of $\Aut_\Cr(T^1_X) \times \Aut_\Cr(T^2_X)$ is a dense open subset of $\Hom_\Cr(T^1_X, T^2_X)$. As a consequence, up to the action of $\Aut_\Cr(T^1_X) \times \Aut_\Cr(T^2_X)$, $h$ can be supposed to be of the form
 $$\mat{h_{11}}{h_{12}}{h_{21}}{h_{22}}$$
 where $h_{11}$ is an automorphism of $T^0$. Then, using the action of
 $$\left(\mat{\id_{T^0}}{-h_{11}^{-1} h_{12}}{0}{\id_{T^1_X}}, \mat{\id_{T^0}}{0}{-h_{21} h_{11}^{-1}}{\id_{T^2_X}}\right) \in \Aut_\Cr(T^1_X) \times \Aut_\Cr(T^2_X)$$
 $h$ can be decomposed as a direct sum and therefore $f$ is not minimal. \cqfd
\end{demo}

The following lemma is inspired from \cite[lemma 3.2]{FuKe} and \cite[theorem 2.3]{DeKe08}:
\begin{lem}
 If $X \in \Cr$ is rigid, then the congruence class of $X$ modulo ${\Gamma}$ is determined by $\ind'_\T(X)$. In other words, if $Y \in \Cr$ is rigid and if $\ind'_\T(Y) = \ind'_\T(X)$, then $X$ and $Y$ are congruent modulo ${\Gamma}$.
\end{lem}

\begin{demo}
 Let  
 $$0 \rightarrow T^1_X \rightarrow T^2_X \xrightarrow{f} X \rightarrow 0 \quad \text{and} \quad 0 \rightarrow T^1_Y \rightarrow T^2_Y \xrightarrow{f'} Y \rightarrow 0$$
 be the admissible short exact sequences which define $\ind_\T(X)$ and $\ind_\T(Y)$.
 
 Using lemma \ref{intvide}, $T^1_X$ and $T^2_X$ on one hand and $T^1_Y$ and $T^2_Y$ on the other hand have no common indecomposable summand. As a consequence, $\ind_\T(X)$ and $\ind_\T(Y)$ fully determine them and $\ind'_\T(X) = \ind'_\T(Y)$ determine their congruence classes modulo ${\Gamma}$. Summing up,
 $$\bigoplus_{g \in {\Gamma}} \g \tens T^1_X \simeq \bigoplus_{g \in {\Gamma}} \g \tens T^1_Y \quad \text{and} \quad \bigoplus_{g \in {\Gamma}} \g \tens T^2_X \simeq \bigoplus_{g \in {\Gamma}} \g \tens T^2_Y.$$
 If one denotes $T^1 = T^1_X$ and $T^2 = T^2_X$, one gets the two following admissible short exact sequences:
 $$0 \rightarrow \bigoplus_{g \in {\Gamma}} \g \tens T^1 \xrightarrow{h} \bigoplus_{g \in {\Gamma}} \g \tens T^2 \rightarrow \bigoplus_{g \in {\Gamma}} \g \tens X \rightarrow 0$$
 and
 $$0 \rightarrow \bigoplus_{g \in {\Gamma}} \g \tens T^1 \xrightarrow{h'} \bigoplus_{g \in {\Gamma}} \g \tens T^2 \rightarrow \bigoplus_{g \in {\Gamma}} \g \tens Y \rightarrow 0.$$
 
 As $\bigoplus_{g \in {\Gamma}} \g \tens X$ and $\bigoplus_{g \in {\Gamma}} \g \tens Y$ are rigid, using lemma \ref{caracgeo}, $h$ and $h'$ are in the same orbit under the action of $\Aut_\Cr(T^1) \times \Aut_\Cr(T^2)$ and, as a consequence,
 $$\bigoplus_{g \in {\Gamma}} \g \tens X \simeq \bigoplus_{g \in {\Gamma}} \g \tens Y$$
 which implies that $X$ and $Y$ are congruent modulo ${\Gamma}$. \cqfd
\end{demo}

Let us now adapt \cite[corollary 4.4]{FuKe}:
\begin{prop} \label{ppfk1}
 If the rank of $\tilde B(\T)$ is full, then
 \begin{enumerate}
  \item $\bar X \mapsto P_{\bar X}$ induces an injection from the set of indecomposable classes of $\Cr$ modulo ${\Gamma}$ such that $\bigoplus_{X \in \bar X} X$ is rigid.
  \item Suppose that $E \subset \Add(\Cr)^{\Gamma}$ is a finite multiset of maximal rigid ${\Gamma}$-stable categories, and fix for every $\T' \in E$, an object $X_{\T'} \in \T'$. If the $X_{\T'}$ are not congruent modulo ${\Gamma}$ pairwise, then the $P_{\bar X_{\T'}}$ are linearly independent.
 \end{enumerate}
\end{prop}

\begin{demo}
 The first point is a direct consequence of the second. Suppose that for some $c_{\T'} \in \Q$,
 $$\sum_{\T' \in E} c_{\T'} P_{\bar X_{\T'}} = 0.$$
 Using the proof of \cite[corollary 4.4]{FuKe}, one deduces the relation
 $$\sum_{\T' \in E} c_{\T'} \prod_{i = 1}^n x_{\ig}^{[\ind_{\T}(X_{\T'}):[T_i]]} = 0$$
 where $[\ind_{\T}(X_{\T'}):[T_i]]$ is the coefficient of $[T_i]$ in the decomposition of $\ind_{\T}(X_{\T'})$ in the basis $\{[T_j]\}_{j \in \llbracket 1, n \rrbracket}$. Hence
 $$\sum_{\T' \in E} c_{\T'} \prod_{\ig \in \bar I} x_{\ig}^{[\ind'_{\T}(X_{\T'}):[\bar T_i]]} = 0$$
 where $[\bar T_i]$ is the orbit of $[T_i]$ in $K_0(\T)/{\Gamma}$. Using previous lemma, the $\ind'_{\T}(X_{\T'})$ are distinct. Hence, every $c_{\T'}$ vanishes. \cqfd
\end{demo}

\begin{cor}
 \label{indlin}
 If $\tilde B(\T)$ is of full rank then the $P_{\bar X}$, where $\bar X$ runs over the equivalence classes of $\Cr$ modulo ${\Gamma}$ such that 
 $$\bigoplus_{X \in \bar X} X$$
 is rigid, are linearly independent over $\Q$. The cluster monomials with coefficients are linearly independent over $\Q$. Equivalently, the cluster monomials without coefficients are linearly independent over the ground ring. 
\end{cor}

\begin{demo}
 The second part is an immediate consequence of the first one because the cluster monomials come from such $\bar X$ through the inclusion $\Ar(B(\T)) \subset \Ar(\Cr, {\Gamma}, \T)$. The first point is a clear consequence of proposition \ref{ppfk1}. \cqfd
\end{demo}

\begin{rem}
 Of course, in proposition \ref{ppfk1} and corollary \ref{indlin}, one can replace the assumption that $\tilde B(\T)$ has maximal rank by the stronger hypothesis that $B(\T)$ has maximal rank.
\end{rem}

\section{Applications}

\subsection{Reminder about root systems and enveloping algebras}

\label{rac}
For more details about this section, see for example \cite[\S 1.4]{Hu02}. Let $\tilde C$ be a symmetrizable generalized Cartan matrix with rows indexed by $\tilde \Delta_0$ and let $\tilde \Delta$ be the bi-valued unoriented graph with vertex-set $\tilde \Delta_0$ such that if $i, j \in \tilde \Delta_0$, there is an edge between $i$ and $j$ if $C_{ij} < 0$ and its valuation is $(-C_{ij}, -C_{ji})$. Let ${\Gamma}$ be a group acting on $\tilde \Delta$ in such a way that $\tilde \Delta$ has no edge between any two vertices of the same ${\Gamma}$-orbit (the action will be said to be admissible).

\begin{df}
 One will denote by $\Delta$ the unoriented graph with vertex-set $\tilde \Delta_0 / {\Gamma}$ and Cartan matrix defined by
 $$C_{\ig \jg} = \frac{1}{\# \jg} \sum_{(i, j) \in \ig \times \jg} \tilde C_{ij}.$$
\end{df}

\begin{lem} \label{conssym}
 Every symmetrizable Cartan matrix can be obtained by this method from a symmetric Cartan matrix and a cyclic group.
\end{lem}

\begin{demo}
 Suppose that $C'$ is a symmetrizable matrix of order $n$. Let $(d_i)_{i \in \llbracket 1, n \rrbracket}$ be the positive integer entries of a diagonal matrix $D$ such that $DC'$ is symmetric. For $i \in \llbracket 1, n \rrbracket$ put $n_i = \prod_{j \neq i} d_j$. Let 
 $$I = \bigcup_{i \in \llbracket 1, n \rrbracket} \{i\} \times \Z / n_i \Z.$$
 One denotes by $a$ the automorphism of $I$ defined by $(i, j) \mapsto (i, j+1)$. For $(i, j)$ and $(i', j')$ in $I$ define
 $$C''_{(i,j), (i', j')} = \left\{
 \begin{array}{ll}
  \frac{d_i C'_{ii'}}{d_i \vee d_{i'}} = \frac{d_{i'} C'_{i'i}}{d_i \vee d_{i'}}  & \quad \text{if } n_i \wedge n_{i'} | j-j' \\
  0 & \quad \text{else.}
 \end{array}
 \right.$$
 It is easy to check that the group generated by $a$ acts on the diagram associated to the Cartan matrix $C''$. Moreover, it is an easy computation to check that the symmetric Cartan matrix obtained from $C''$ and $a$ is $C'$. \cqfd
\end{demo}

\begin{rem}
 An other proof of lemma \ref{conssym} is given in \cite[proposition 14.1.2]{Lu93}.
\end{rem}

Let $\ggo$ and $\tilde \ggo$ be the Kac-Moody Lie algebras associated to $C$ and $\tilde C$. One denotes by $(\tilde e_i)_{i \in Q_0}$ and $(\tilde f_i)_{i \in Q_0}$ (resp. $(e_{\ig})_{\ig \in Q_0 / \Gamma}$ and $(f_{\ig})_{\ig \in Q_0 / \Gamma}$) the Chevalley generators of $\tilde \ggo$ (resp. $\ggo$). One sets $h_i = [e_i, f_i]$ and $\tilde h_i = [\tilde e_i, \tilde f_i]$. Let $\ngo$ (resp. $\ngo_-$) be the nilpotent subalgebra (resp. opposite nilpotent subalgebra) of $\ggo$ generated by the $(e_{\ig})_{\ig \in Q_0 / \Gamma}$ (resp. by the $(f_{\ig})_{\ig \in Q_0 / \Gamma}$). Let $\bgo$ (resp. $\bgo_-$) be the Borel subalgebra (resp. opposite Borel subalgebra) of $\ggo$ generated by the $(h_{\ig})_{\ig \in Q_0 / \Gamma}$ and $\ngo$ (resp. $\ngo_-$). Let $\tilde \ngo$ (resp. $\tilde \ngo_-$) be the nilpotent subalgebra (resp. opposite nilpotent subalgebra) of $\tilde \ggo$ generated by the $(\tilde e_i)_{i \in Q_0}$ (resp. by the $(\tilde f_i)_{i \in Q_0}$.  Let $\tilde \bgo$ (resp. $\tilde \bgo_-$) be the Borel subalgebra (resp. opposite Borel subalgebra) of $\tilde \ggo$ generated by the $(\tilde h_i)_{i \in Q_0}$ and $\tilde \ngo$ (resp. $\tilde \ngo_-$). 

% Let $G$, $\tilde G$, $B$, $\tilde B$, $B_-$, $\tilde B_-$, $N$ and $\tilde N_-$ be the Kac-Moody group associated to these Lie algebras (see \cite{Ku02} and \cite{GeLeSc}).

As ${\Gamma}$ acts on $\tilde \Delta$, it acts also on $\tilde \ggo$ by extending its action on the Chevalley generators.

\begin{prop}[\citeb{theorem 7.1.5}{Hu02}] \label{hub}
 There is a monomorphism of Lie algebras $\ggo \inj \tilde \ggo^{\Gamma}$
 $$e_{\ig} \mapsto \sum_{i \in \ig} \tilde e_i \quad f_{\ig} \mapsto \sum_{i \in \ig} \tilde f_i \quad  h_{\ig} \mapsto \sum_{i \in \ig} \tilde h_i $$
 which can be restricted to a monomorphism $\ngo \inj \tilde \ngo^{\Gamma}$. If $C$, or equivalently $\tilde C$, is of Dynkin type, this monomorphism is an isomorphism.
\end{prop}

\begin{cor}
 \label{defkappa}
 There is an epimorphism
 $$\kappa : U(\tilde \ngo)^*_{\gr} / {\Gamma} \surj U(\ngo)^*_{\gr}$$
  where the quotient by ${\Gamma}$ has to be understood as the quotient by the ideal generated by the elements of the form $(g f - f)$ for $g \in {\Gamma}$ and $f \in U(\tilde \ngo)^*_{\gr}$. Here, $U(\tilde \ngo)^*_{\gr}$ and $U(\ngo)^*_{\gr}$ denote graded dual spaces. If $C$ is of Dynkin type, $\kappa$ is an isomorphim.
\end{cor}

\begin{demo}
 It is a clear translation of proposition \ref{hub}. \cqfd %The two isomorphisms are classical (see for example \cite[proposition 5.1]{GeLeSc08-1}).  \cqfd
\end{demo}

\begin{lem}[\citeb{proposition 4}{Hu04}]
 \label{racns}
 Let $R$ be a root system of type $\tilde \Delta$. Let $V$ be the lattice generated by $R$. Then the linear map
 \begin{align*}
  \alpha: V &\rightarrow V \\
  v &\mapsto \sum_{g \in {\Gamma}} g v
 \end{align*}
 maps $R$ to a root system of type $\Delta$. 
\end{lem}

\subsection{$\Sub I_J$ and partial flag varieties}

\label{subqj}
This application generalizes \cite{GeLeSc08}. Let $Q$ be a quiver such that the underlying unoriented graph $\tilde \Delta$ is a Dynkin diagram of type $A$, $D$ or $E$. Let ${\Gamma}$ be a group acting on $Q$ in such a way that $Q$ has no arrow between any two vertices of the same orbit (the action will be said to be admissible). It induces an action on $\tilde \Delta$. We denote by $Q_{\Gamma}$ the same quiver as in section \ref{acc} and by $\tilde \Delta_{\Gamma}$ the underlying unoriented graph. We denote by $\Delta$ the diagram defined in section \ref{rac}.

\newcommand{\lmp}[1]
{\begin{minipage}{4.5cm}#1\end{minipage}} 
Here is the list of all possible cases, where ${\Gamma}$ acts faithfully on a Dynkin diagram $\tilde \Delta$.

\newcommand{\mpc}[1]{\begin{minipage}{2cm} \begin{center} \vspace*{.1cm} #1 \end{center} \end{minipage}}
\newcommand{\mpcd}[1]{\begin{minipage}{.7cm} \begin{center} \vspace*{.1cm} #1 \end{center} \end{minipage}}

\begin{center}
  \tiny
  \begin{longtable}{|c|c|c|c|}
    \hline
    {\normalsize \mpc{$\tilde \Delta$}} & {\normalsize \mpcd{${\Gamma}$}} & {\normalsize \mpc{$\Delta$}} & {\normalsize \mpc{$\tilde \Delta_{\Gamma}$}} \\
    \hline
   \endfirsthead
    \hline
    {\normalsize \mpc{$\tilde \Delta$}} & {\normalsize \mpcd{${\Gamma}$}} & {\normalsize \mpc{$\Delta$}} & {\normalsize \mpc{$\tilde \Delta_{\Gamma}$}} \\
    \hline
   \endhead
    \lmp{
    $$\xymatrix@R=.27cm@C=.27cm{
     1 \ar@{-}[r] & 2 \ar@{-}[r] & \dots \ar@{-}[r] & n-1 \ar@{-}[rd] &  \\
      & & & & n \\
     1' \ar@{-}[r] & 2' \ar@{-}[r] & \dots \ar@{-}[r] & (n-1)' \ar@{-}[ru] & 
    }$$}
    &
    $\Z/2\Z$
    &
    \lmp{$$\xymatrix@R=.27cm@C=.27cm{
     1 \ar@{-}[r] & 2 \ar@{-}[r] & \dots \ar@{-}[r] & n-1 \ar@{=}[rr] & \text{\LARGE <} & n  \\
    }$$}
    &
    \lmp{$$\xymatrix@R=.27cm@C=.27cm{
     & & & & n_+ \\
     1 \ar@{-}[r] & 2 \ar@{-}[r] & \dots \ar@{-}[r] & n-1 \ar@{-}[ur] \ar@{-}[dr] &  \\
     & & & & n_- \\
    }$$} \\ 
    \hline
    \lmp{$$\xymatrix@R=.27cm@C=.27cm{
     & & & & n \\
     1 \ar@{-}[r] & 2 \ar@{-}[r] & \dots \ar@{-}[r] & n-1 \ar@{-}[ur] \ar@{-}[dr] &  \\
     & & & & n' \\
    }$$}
    &
    $\Z/2\Z$
    &
    \lmp{$$\xymatrix@R=.27cm@C=.27cm{
     1 \ar@{-}[r] & 2 \ar@{-}[r] & \dots \ar@{-}[r] & n-1 \ar@{=}[rr] & \text{\LARGE >} & n  \\
    }$$}
    &
    \lmp{
    $$\xymatrix@R=.27cm@C=.27cm{
     1_+ \ar@{-}[r] & 2_+ \ar@{-}[r] & \dots \ar@{-}[r] & (n-1)_+ \ar@{-}[rd] &  \\
      & & & & n \\
     1_- \ar@{-}[r] & 2_- \ar@{-}[r] & \dots \ar@{-}[r] & (n-1)_- \ar@{-}[ru] & 
    }$$} \\ 
    \hline
    \lmp{$$\xymatrix@R=.27cm@C=.27cm{
     1 \ar@{-}[dr] &  \\
     1' \ar@{-}[r] & 2 \\
     1'' \ar@{-}[ur] &  \\
    }$$}
    &
    $\Z/3\Z$
    &
    \lmp{$$\xymatrix@R=.27cm@C=.27cm{
     1 \ar@3{-}[rr] & \text{\LARGE <} & 2 \\
    }$$}
    &
    \lmp{$$\xymatrix@R=.27cm@C=.27cm{
       & 2_1 \\
     1 \ar@{-}[r] \ar@{-}[ur] \ar@{-}[dr] & 2_{e^{\frac{2i\pi}{3}}} \\
       & 2_{e^{-\frac{2i\pi}{3}}}
    }$$} \\ 
    \hline
    \lmp{$$\xymatrix@R=.27cm@C=.27cm{
     1 \ar@{-}[dr] &  \\
     1' \ar@{-}[r] & 2 \\
     1'' \ar@{-}[ur] &  \\
    }$$}
    &
    $\mathfrak{S}_3$
    &
    \lmp{$$\xymatrix@R=.27cm@C=.27cm{
     1 \ar@3{-}[rr] & \text{\LARGE <} & 2 \\
    }$$}
    &
    \lmp{$$\xymatrix@R=.27cm@C=.27cm{
       & 2_+ \\
     1_+ \ar@{-}[ur] \ar@{-}[dr] &  \\
       & 2_2 \\
     1_- \ar@{-}[ur] \ar@{-}[dr] &  \\
       & 2_-
    }$$} \\ 
    \hline
    \lmp{$$\xymatrix@R=.27cm@C=.27cm{
     1 \ar@{-}[dr] & & & \\
     & 2 \ar@{-}[dr] & & \\
     & & 3 \ar@{-}[r] & 4 \\
     & 2' \ar@{-}[ur] & & \\
     1' \ar@{-}[ur] & & &
    }$$}
    &
    $\Z/2\Z$
    &
    \lmp{$$\xymatrix@R=.27cm@C=.27cm{
     1 \ar@{-}[r] & 2 \ar@{=}[rr] & \text{\LARGE <} & 3 \ar@{-}[r] & 4
    }$$}
    &
    \lmp{$$\xymatrix@R=.27cm@C=.27cm{
     & & & 4_+ \\
     & & 3_+ \ar@{-}[ur] & \\
     1 \ar@{-}[r] & 2 \ar@{-}[ur] \ar@{-}[dr]& & \\
     & & 3_- \ar@{-}[dr]& \\
     & & & 4_-
    }$$} \\ 
    \hline 
  \end{longtable}
\end{center}

\begin{rem}
  Observe that all non simply-laced Dynkin diagrams can be realized for appropriate $\tilde \Delta$ and ${\Gamma}$.
\end{rem}

One retains the notation of section \ref{rac}. Let $N$ and $\tilde N$ be the Lie groups associated with $\ngo$ and $\tilde \ngo$.

\begin{nt}
 If $i \in Q_0$, $x_i$ denotes the one-parameter subgroup of $\tilde N$ defined by
 $$x_i(t) = \exp(t e_i).$$
\end{nt}

%\begin{df}
% For every finite sequence $(i_1, i_2, \dots, i_n)$ of $\tilde \Delta_0$, one gets, for all $\phi \in \C[\tilde N]$ and $g \in {\Gamma}$,
% $$(g \phi)(x_{g i_1}(t_1) x_{g i_2}(t_2) \dots x_{g i_n}(t_n)) = \phi(x_{i_1}(t_1) x_{i_2}(t_2) \dots x_{i_n}(t_n))$$
% and this property defines entierly an action of ${\Gamma}$ on $\C[\tilde N]$. 
%\end{df}

Let $X \in \md \Lambda_Q$ and $(i) = i_1 i_2 \dots i_n$ be a word on $Q_0$. One will denote by $\Phi_{X, (i)}$ the (closed) subvariety of
$$\Gr_0(X) \times \Gr_1(X) \times \dots \times \Gr_{n-2}(X) \times \Gr_{n-1}(X)$$
consisting of the $(X_0, X_1, \dots, X_{n-2}, X_{n-1})$ such that
\begin{itemize}
 \item for all $j \in \llbracket 1, n-1 \rrbracket$, $X_{j-1} \subset X_j$;
 \item for all $j \in \llbracket 1, n-1 \rrbracket$, $X_j / X_{j-1} \simeq S_{i_j}$;
 \item $X/X_{n-1} \simeq S_{i_n}$.
\end{itemize}

The following result is obtained by duality from the Lagrangian construction of $U(\ngo)$ by Lusztig \cite{Lu91}, \cite{Lu00}.

\begin{thm}[\citeb{\S 9}{GeLeSc06-1}]
 \label{defphi}
 For all $X \in \md \Lambda_Q$, there exists a unique $\phi_X \in \C[N]$ such that for every word $(i)$ on $Q_0$, and every $\tg \in \C^n$ (where $n$ is the length of $(i)$), 
 $$\phi_X(x_{i_1}(t_1) x_{i_2}(t_2) \dots x_{i_n}(t_n)) = \sum_{\ag \in \N^n} \chi\left(\Phi_{X, (i)^\ag}\right) \frac{\tg^\ag}{\ag !}.$$
 Here, $\displaystyle (i)^\ag = \underbrace{i_1 \dots i_1}_{a_1} \underbrace{i_2 \dots i_2}_{a_2} \dots \underbrace{i_n \dots i_n}_{a_n}$, $\displaystyle \tg^\ag = t_1^{a_1} t_2^{a_2} \dots t_n^{a_n}$ and $\displaystyle \ag ! = a_1 ! a_2 ! \dots a_n !$.
\end{thm}

Using the dualities $\C[N] \simeq U(\ngo)^*_{\gr}$ and $\C[\tilde N] \simeq U(\tilde \ngo)^*_{\gr}$, we can lift the action of $\Gamma$ on $U(\tilde \ngo)^*_{\gr}$ to an action of $\Gamma$ on $C[\tilde N]$. Therefore the isomorphism $\kappa$ defined in corollary \ref{defkappa} can be lifted to an isomorphism $\kappa : \C[\tilde N] / \Gamma \simeq \C[N]$ where the quotient by $\Gamma$ is the quotient by the ideal generated by the elements of the form $(gf - f)$ for $g \in \Gamma$ and $f \in \C[\tilde N]$. 

\begin{nt}
 \label{defpsi}
 For $X \in \md \Lambda_Q$, let $\psi_X = \kappa(\pi(\phi_X))$ where $\pi : \C[\tilde N] \rightarrow \C[\tilde N]/\Gamma$ is the canonical projection.
\end{nt}

The action of ${\Gamma}$ on $Q$ induces an action of ${\Gamma}$ on $kQ$, then on $\Lambda_Q$, then on $\md \Lambda_Q$ (see section \ref{acc}). Following the proof of \cite[theorem 3, \S 8]{GeLeSc07-1}, it is easy to see that the following diagram commutes:
 $$\xymatrix{
  \Ext^1_{\Lambda_Q}(-_1, -_2) \ar[d]^c \ar[rr]^{\g \tens -} & & \Ext^1_{\Lambda_Q}(\g \tens -_1, \g \tens -_2) \ar[d]^c \\
  \Ext^1_{\Lambda_Q}(-_2, -_1)^* & & \Ext^1_{\Lambda_Q}(\g \tens -_2, \g \tens -_1)^* \ar[ll]^{(\g \tens -)^*}
 }$$
where $c$ is the functorial isomorphism from $\Ext^1_{\Lambda_Q}(-_1, -_2)$ to $\Ext^1_{\Lambda_Q}(-_2, -_1)^*$. In other terms, the action of ${\Gamma}$ on $\md \Lambda_Q$ is $2$-Calabi-Yau in the sense of definition \ref{act2CY}.

\begin{lem}
 For $X \in \md \Lambda_Q$ and $g \in {\Gamma}$, 
 $$\psi_{\g \tens X} = \psi_X.$$
\end{lem}

\begin{demo}
 One has
 \begin{align*}
  \phi_{\g \tens X}(x_{i_1}(t_1) x_{i_2}(t_2) \dots x_{i_n}(t_n)) &= \sum_{\ag \in \N^n} \chi\left(\Phi_{\g \tens X, (i)^\ag}\right) \frac{\tg^\ag}{\ag !} \\
  &= \sum_{\ag \in \N^n} \chi\left(\Phi_{X, {g^{-1} \cdot (i)}^\ag}\right)  \frac{\tg^\ag}{\ag !}\\
  &= \phi_{X}(x_{g^{-1} \cdot i_1}(t_1) x_{g^{-1} \cdot i_2}(t_2) \dots x_{g^{-1} \cdot i_n}(t_n))
 \end{align*}
 which implies the result. \cqfd
\end{demo}

\begin{nt}
 One will denote
 $$\psi_{\bar X} = \psi_X$$
 where $\bar X$ is the ${\Gamma}$-orbit of $X$.
\end{nt}

\begin{thm}[\citeb{lemma 7.3}{GeLeSc05} and \citeb{theorem 9.2}{GeLeSc06-1}]
 \begin{enumerate}
  \item If $X, Y \in \md \Lambda_Q$ then $\phi_{X \oplus Y} = \phi_{X} \phi_{Y}$.
  \item If $X, Y \in \md \Lambda_Q$ and $\dim \Ext^1_{\Lambda_Q}(X, Y) = 1$, and if one considers two non-split short exact sequences 
  $$0 \rightarrow X \rightarrow Z \rightarrow Y \rightarrow 0 \quad \text{and} \quad 0 \rightarrow Y \rightarrow Z' \rightarrow X \rightarrow 0$$
  then $\phi_{X} \phi_{Y} = \phi_{Z} + \phi_{Z'}$.
 \end{enumerate}
\end{thm}

\begin{cor}
 \label{fmut2}
 \begin{enumerate}
  \item If $X, Y \in \md \Lambda_Q$, $$\psi_{\bar X \oplus \bar Y} = \psi_{\bar X} \psi_{\bar Y}.$$
  \item If $X, Y \in \md \Lambda_Q$ and $\dim \Ext^1_{\Lambda_Q}(X, Y) = 1$, and if one considers two non-split short exact sequences 
  $$0 \rightarrow X \rightarrow Z \rightarrow Y \rightarrow 0 \quad \text{and} \quad 0 \rightarrow Y \rightarrow Z' \rightarrow X \rightarrow 0$$
  then $\psi_{\bar X} \psi_{\bar Y} = \psi_{\bar Z} + \psi_{\bar Z'}$.
 \end{enumerate}
\end{cor}

\begin{demo}
 This follows immediately from the fact that $\psi_{\bar X}$ is the image of $\phi_X$ under a ring homomorphism. \cqfd
\end{demo}

Let now $J \subset Q_0$ be non-empty, ${\Gamma}$-stable and $K = Q_0 \setminus J$. 

\begin{df}
 For $j \in Q_0$, one denotes by $I_j$ the injective $\Lambda_Q$-module of socle $S_j$. Put
 $$I_J = \bigoplus_{j \in J} I_j$$
 and denote by $\Sub I_J$ the full subcategory of $\md \Lambda_Q$ whose objects are isomorphic to submodules of $I_J^{\oplus n}$ for some $n \in \N$.
\end{df}

\begin{prop}[\citeb{\S 3}{GeLeSc08}]
 The category $\Sub I_J$ is an exact, $\Hom$-finite, Krull-Schmidt, Frobenius, and $2$-Calabi-Yau subcategory of $\md \Lambda_Q$.
\end{prop}

\begin{lem}
 All projective objects of $\Sub I_J$ have right rigid quasi-approximations.
\end{lem}

\begin{demo}
 First of all, for every simple $\Lambda_Q$-module $S$, 
 $$\bigoplus_{g \in {\Gamma}} \g \tens S$$
 is rigid because the action of ${\Gamma}$ on $Q$ is admissible. Moreover, the injective objects $L_i$ of $\Sub I_J$ constructed explicitly in \cite[\S 3.3]{GeLeSc08} have simple heads. It is now clear that lemma \ref{qapsc} can be applied. \cqfd
\end{demo}

\begin{nt}
 If $i \in Q_0$, one denotes by $\Er_i^\dag$ the functor from $\md \Lambda_Q$ to itself that maps a module to its quotient by the largest possible power of $S_i$.
\end{nt}

\begin{prop}[\citeb{proposition 5.1}{GeLeSc08}]
 The $\Er_i^\dag$ satisfy the following relations
 \begin{enumerate}
  \item $\Er_i^\dag \Er_i^\dag = \Er_i^\dag$;
  \item $\Er_i^\dag \Er_j^\dag = \Er_j^\dag \Er_i^\dag$ if there is no edge between $i$ and $j$ in $\Delta$;
  \item $\Er_i^\dag \Er_j^\dag \Er_i^\dag = \Er_j^\dag \Er_i^\dag \Er_j^\dag$ if there is an edge between $i$ and $j$ in $\Delta$.
 \end{enumerate}
\end{prop}

\begin{nt}
 If $\ig \in Q_0/{\Gamma}$, one denotes
 $$\Erg_{\ig}^\dag = \prod_{i \in \ig} \Er_i^\dag$$
 which is well-defined because the factors $\Er_i^\dag$ in the product commute. This functor maps the ${\Gamma}$-stable objects of $\md \Lambda_Q$ to ${\Gamma}$-stable objects.
\end{nt}

Denote $\tilde \Delta_K$ (resp. $\Delta_K$) the restriction of the diagram $\tilde \Delta$ to $K$ (resp. of the diagram $\Delta$ to $K/{\Gamma}$). Denote by $W$, $\tilde W$, $W_K$ and $\tilde W_K$ the Weyl groups of $\Delta$, $\tilde \Delta$, $\Delta_K$ and $\tilde \Delta_K$.

One denotes by $(\sigma_i)_{i \in Q_0}$ (resp. $(\sigma_{\ig})_{i \in Q_0/{\Gamma}}$) the generators of $\tilde W$ (resp. $W$). One gets an injective morphism
\begin{align*}
 W &\rightarrow \tilde W \\
 \sigma_{\ig} &\mapsto \prod_{i \in \ig} \sigma_i
\end{align*}
which restricts to a morphism from $W_K$ to $\tilde W_K$.

\begin{prop} \label{propTi}
 Let $(\ig)$ be a reduced expression of the longest element of $W$ and $\ell$ be its length. Assume that $(\ig)$ has a left factor which is a reduced expression of the longest element of $W_K$. Then
 $$T_{(\ig)} = \bigoplus_{m = 1}^\ell \Erg_{\ig_1}^\dag \Erg_{\ig_2}^\dag \dots \Erg_{\ig_m}^\dag\left(\bigoplus_{i \in \ig_m} I_i\right) \oplus \bigoplus_{i \in Q_0} I_i$$
 has a direct summand $T_{(\ig), K}$ which is maximal rigid and ${\Gamma}$-stable in $\Sub I_J$.
\end{prop}

\begin{demo}
 The only thing to add to \cite[proposition 7.3]{GeLeSc08} is that $T_{(\ig), K}$ is ${\Gamma}$-stable. It is clear by definition of the functors $\Erg_{\ig}^\dag$. \cqfd
\end{demo}

\newcommand{\Pa}{
\begin{minipage}{76.2pt}\tiny 
$$\xymatrix@R=.2cm@C=.2cm{
 & & & & a' \ar[dl] \\
 & & & b' \ar[dl] & \\
 & & c \ar[dl] & & \\
 & b \ar[dl] & & & \\
 a & & & & 
}$$
\end{minipage}}

\newcommand{\Pb}{ 
\begin{minipage}{76.2pt}\tiny
$$\xymatrix@R=.2cm@C=.2cm{
 & & & b' \ar[dl] \ar[dr] & \\
 & & c \ar[dl] \ar[dr] & & a' \ar[dl] \\
 & b \ar[dl] \ar[dr] & & b' \ar[dl] & \\
 a \ar[dr] & & c \ar[dl] & & \\
 & b & & & 
}$$
\end{minipage}}

\newcommand{\Pc}{ 
\begin{minipage}{76.2pt}\tiny
$$\xymatrix@R=.2cm@C=.2cm{
 & & c \ar[dl] \ar[dr] & & \\
 & b \ar[dl] \ar[dr] & & b' \ar[dl] \ar[dr] & \\
 a \ar[dr] & & c \ar[dl] \ar[dr] & & a' \ar[dl] \\
 & b \ar[dr] & & b' \ar[dl] & \\
 & & c & &
}$$
\end{minipage}}

\newcommand{\Pd}{ 
\begin{minipage}{76.2pt}\tiny
$$\xymatrix@R=.2cm@C=.2cm{
 & b \ar[dl] \ar[dr] & & & \\
 a \ar[dr] & & c \ar[dl] \ar[dr] & & \\
 & b \ar[dr] & & b' \ar[dl] \ar[dr] &  \\
 & & c \ar[dr] & & a' \ar[dl] \\
 & & & b' &
}$$
\end{minipage}}

\newcommand{\Pe}{ 
\begin{minipage}{76.2pt}\tiny
$$\xymatrix@R=.2cm@C=.2cm{
 a \ar[dr] & & & & \\
 & b \ar[dr] & & & \\
 & & c \ar[dr] & & \\
 & & & b' \ar[dr] & \\
 & & & & c' 
 }$$
\end{minipage}}

\newcommand{\MA}{ 
\begin{minipage}{9.4pt}\tiny
$$\xymatrix@R=.2cm@C=.2cm{
 c
 }$$
\end{minipage}}

\newcommand{\MB}{ 
\begin{minipage}{27.1pt}\tiny
$$\xymatrix@R=.2cm@C=.2cm{
 & b' \ar[dl] \\
 c &
 }$$
\end{minipage}}

\newcommand{\MC}{ 
\begin{minipage}{27.1pt}\tiny
$$\xymatrix@R=.2cm@C=.2cm{
 b \ar[dr] & \\
 & c 
 }$$
\end{minipage}}

\newcommand{\MD}{ 
\begin{minipage}{45.5pt}\tiny
$$\xymatrix@R=.2cm@C=.2cm{
 & c \ar[dl] \ar[dr] & \\
 b \ar[dr] & & b' \ar[dl] \\
 & c & 
 }$$
\end{minipage}}

\newcommand{\ME}{ 
\begin{minipage}{45.5pt}\tiny
$$\xymatrix@R=.2cm@C=.2cm{
 & & a' \ar[dl] \\
 & b' \ar[dl] & \\
 c & &
 }$$
\end{minipage}}

\newcommand{\MF}{ 
\begin{minipage}{45.5pt}\tiny
$$\xymatrix@R=.2cm@C=.2cm{
 a \ar[dr] & & \\
 & b \ar[dr] & \\
 & & c 
 }$$
\end{minipage}}

\newcommand{\MG}{ 
\begin{minipage}{60.5pt}\tiny
$$\xymatrix@R=.2cm@C=.2cm{
 & & b' \ar[dl] \ar[dr] & \\
 & c \ar[dl] \ar[dr] & & a' \ar[dl] \\
 b \ar[dr] & & b' \ar[dl] & \\
 & c & & \\
}$$
\end{minipage}}

\newcommand{\MH}{ 
\begin{minipage}{60.5pt}\tiny
$$\xymatrix@R=.2cm@C=.2cm{
 & b \ar[dl] \ar[dr] & & \\
 a \ar[dr] & & c \ar[dl] \ar[dr] & \\
 & b \ar[dr] & & b' \ar[dl]  \\
 & & c & \\
}$$
\end{minipage}}

\newcommand{\MI}{
\begin{minipage}{60.5pt}\tiny 
$$\xymatrix@R=.2cm@C=.2cm{
 & & & a' \ar[dl] \\
 & & b' \ar[dl] & \\
 & c \ar[dl] & & \\
 b & & & \\
}$$
\end{minipage}}

\newcommand{\MJ}{ 
\begin{minipage}{60.5pt}\tiny
$$\xymatrix@R=.2cm@C=.2cm{
 a \ar[dr] & & & \\
 & b \ar[dr] & & \\
 & & c \ar[dr] & \\
 & & & b'  \\
 }$$
\end{minipage}}

\begin{ex}
 \label{exTi}
 Suppose here that $\tilde \Delta = A_5$ is indexed in the following way:
 $$\xymatrix@R=.5cm@C=.5cm{
     a \ar@{-}[r] & b \ar@{-}[r] &  c \ar@{-}[r]& b' \ar@{-}[r] & a'
 }$$
 on which ${\Gamma} = \Z/2\Z$ acts in the only non-trivial possible way. Hence $\Delta = C_3$ indexed in the following way:
 $$\xymatrix@R=.5cm@C=.5cm{
     a \ar@{-}[rr] & & b \ar@{=}[rr] & \text{\LARGE <} & c
 }$$
 Let $J/{\Gamma} = \{c\}$. Then $(\ig) = (a,b,a,c,b,a,c,b,c)$ is suitable. Then, it is easy to compute $T_{(\ig)}$ and therefore $T_{(\ig), K}$ (see figure \ref{figTiTiK}).
 \begin{figure}
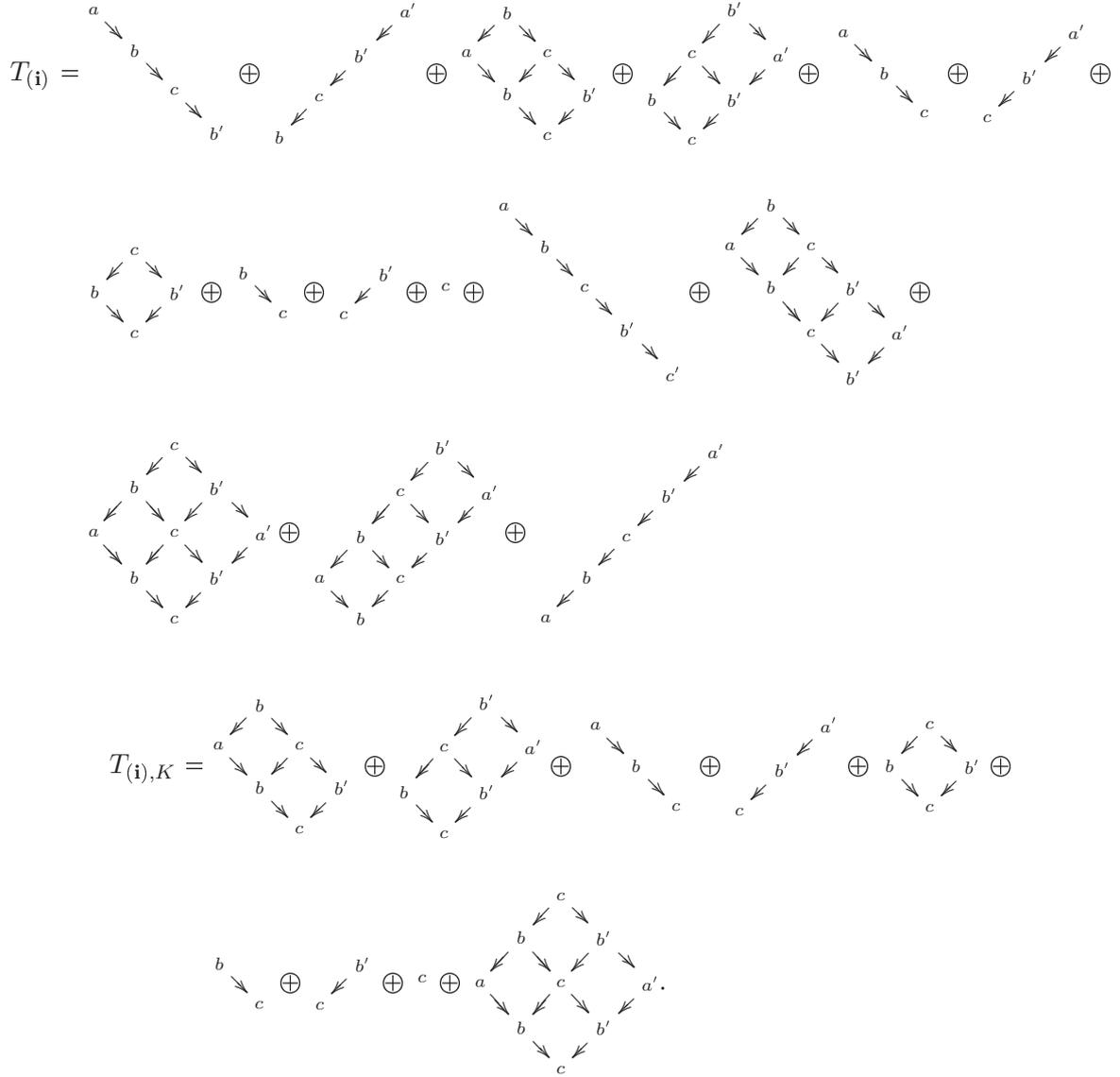

 \begin{align*}
  T_{(\ig)} =& \MJ \oplus \MI \oplus \MH \oplus \MG \oplus \MF \oplus \ME \oplus \\ &\MD \oplus \MC \oplus \MB \oplus \MA \oplus \Pe \oplus \Pd \oplus \\ &\Pc \oplus \Pb \oplus \Pa
 \end{align*} 
  \begin{align*}
  T_{(\ig), K} =& \MH \oplus \MG \oplus \MF \oplus \ME \oplus \MD \oplus \\ &\MC \oplus \MB \oplus \MA \oplus \Pc.
 \end{align*}
 \caption{Explicit computation of $T_{(\ig)}$ and $T_{(\ig), K}$} \label{figTiTiK} 
 \end{figure}
\end{ex}

Fix now $T = T_{(\ig), K}$ as in proposition \ref{propTi}. Let $\T = \add(T)$.

\begin{lem}
 The category $\T$ has no ${\Gamma}$-loops nor ${\Gamma}$-$2$-cycles.
\end{lem}

\begin{demo}
 By \cite{GeLeSc08}, $\T$ is cluster-tilting and $\End_{\Sub I_J}(T)$ is of finite global dimension. The result follows by theorem \ref{recapend}. \cqfd
\end{demo}

Hence, we can apply the results of section \ref{caramas}. We retain the notation of section \ref{caramas}.

\begin{figure}
\begin{center}
 \setlength{\unitlength}{1pt}
 \begin{picture}(403,406)
  \put(161.95,22.75){\MF}
  \put(81.15,113.05){\MC}
  \put(244.45,113.05){\MH}
  \put(0,202.75){\MA}
  \put(161.95,202.75){\MD}
  \put(326.6,202.75){\Pc}
  \put(81.15,292.75){\MB}
  \put(244.45,292.75){\MG}
  \put(161.95,382.75){\ME}
  \put(18.4,208.45){\vector(1,1){61.75}}
  \put(108.4,118.45){\vector(1,1){61.75}}
  \put(198.4,28.45){\vector(1,1){61.75}}
  \put(117.25,307.3){\vector(1,1){43.7}}
  \put(207.25,217.3){\vector(1,1){43.7}}
  \put(297.25,127.3){\vector(1,1){43.7}}
  \put(18.4,189.05){\vector(1,-1){61.75}}
  \put(108.4,279.05){\vector(1,-1){61.75}} 
  \put(198.4,369.05){\vector(1,-1){61.75}}   
  \put(117.25,90.2){\vector(1,-1){43.7}}
  \put(207.25,180.2){\vector(1,-1){43.7}}
  \put(297.25,270.2){\vector(1,-1){43.7}}
  \put(160.95,198.75){\vector(-1,0){142.55}}
  \put(243.45,288.75){\vector(-1,0){126.2}}
  \put(243.45,109.05){\vector(-1,0){126.2}}
  \put(325.6,198.75){\vector(-1,0){109.15}}
 \end{picture}
 \caption{Auslander-Reiten quiver of $\T$}
 \label{ART}
 \end{center}
\end{figure}
 
\begin{ex}
 Continue with example \ref{exTi}. The Auslander-Reiten quiver of $\T$ is displayed in figure \ref{ART}. The action of $\Gamma$ corresponds to the reflection in the middle horizontal $\tau$-orbit. The projective-injective objects are the five rightmost ones. Indexing the lines of the exchange matrix by the $\Gamma$-orbits of
 $$\MD, \MB, \MA, \ME, \MG, \Pc$$
 in this order, one gets
 $$B(\T) = \left(\begin{matrix}
                  0 & -1 & 1 \\
                  2 & 0 & -2 \\
                  -1 & 1 & 0 \\
                  0 & -1 & 0 \\
                  -2 & 1 & 0 \\
                  1 & 0 & 0
                 \end{matrix}\right).$$
\end{ex}

For $X \in \md \Lambda_Q$ and $\ig \in Q_0/{\Gamma}$, let $k_{\ig}(X)$ be the total dimension of the maximal submodule of $X$ supported by $\ig$. Let $\mathcal{R}$ be the set of isomorphism classes of $\Lambda_Q$-modules $X$ such that
 $$\bigoplus_{g \in {\Gamma}} \g \tens X$$
 is rigid. If $\dg \in \N^{Q_0 / {\Gamma}}$, $\mathcal{R}_\dg$ is the set of elements of $\mathcal{R}$ of dimension vector (summed on ${\Gamma}$-orbits) $\dg$ and if $k \in \N$, $\mathcal{R}_{\dg, \ig, k} = \{X \in \mathcal{R}_\dg \,|\, k_{\ig}(X) = k \}$. One will denote by $\sim$ the equivalence relation on $\mathcal{R}$ identifying $X$ and $Y$ if 
 $$\bigoplus_{g \in {\Gamma}} \g \tens X \simeq \bigoplus_{g \in {\Gamma}} \g \tens Y.$$

\begin{lem}
 With the previous notation, $\Erg_{\ig}^\dag$ induces an injection from $\mathcal{R}_{\dg, \ig, k} / \sim$ to $\mathcal{R}_{\dg - k \ig, \ig, 0} / \sim$.
\end{lem}

\begin{demo}
 First of all, the map is clearly defined. Let now $X, Y \in \mathcal{R}_{\dg, \ig, k}$ such that $\Erg_{\ig}^\dag(X) \sim \Erg_{\ig}^\dag(Y)$. Define $\tilde X = \bigoplus_{g \in {\Gamma}} \g \tens X$ and $\tilde Y = \bigoplus_{g \in {\Gamma}} \g \tens Y$. These two rigid modules have the same dimension vector. Moreover, at each vertex of $\ig$, they have the same socle. Finally $\Erg_{\ig}^\dag(\tilde X) \simeq \Erg_{\ig}^\dag(\tilde Y)$. Using \cite[lemma 12.5 (e)]{Lu91} together with the result stating that the orbit of a rigid module is a dense open subset of its irreducible component in the module variety \cite[corollary 3.15]{GeLeSc06-1}, it implies that $\tilde X \simeq \tilde Y$. Hence $X \sim Y$ which is the claimed result. \cqfd
\end{demo}

For $X \in \md \Lambda_Q$, denotes by $\epsilon_{\bar X}$ the element of $U(\mathfrak n)^*$ corresponding to $\psi_{\bar X}$. If $f \in U(\mathfrak n)$, $\ig \in Q_0 / \Gamma$ and $n \in \N$, one gets
$$\epsilon_{\bar X}(f e_\ig^n) = \sum_{k \in \Z} k \chi\left(\left\{Y \in \Gr_n(X) \, | \, \udim Y = \ig^n \text{ and } \epsilon_{\widebar {X/Y}}(f) = k \right\}\right)$$
where $\udim Y$ denotes the dimension vector (summed on $\Gamma$-orbits) of $Y$ and $\chi$ denotes the Euler characteristic (it follows for example from \cite[\S 5.2]{GeLeSc05}).

\begin{prop}
 \begin{enumerate} \label{partscan}
  \item \label{partscan1}For all $X \in \mathcal{R}$, there exists $f_{\bar X} \in U(\mathfrak n)$ homogeneous of degree $\dg$ such that for every $Y \in \mathcal{R}$, $$\epsilon_{\bar Y} (f_{\bar X}) = \left\{\begin{array}{ll}
   1 & \text{if } X \sim Y; \\
   0 & \text{else.}
  \end{array}  
  \right.$$
  \item \label{partscan2}The $\psi_{\bar X}$, where $\bar X$ runs over the equivalence classes of indecomposable objects of $\T$, are algebraically independent.
  \item \label{partscan3}Cluster monomials of $\C[N]$ are linearly independent.
 \end{enumerate}
\end{prop}

\begin{rem}
 The following proof is an adaptation of the proof of existence of the semicanonical basis by Lusztig \cite{Lu00}. We suggested a proof in \cite{De08} which was very close to this one, with a dual description of the cluster character. Unfortunately, we are not able to construct in this way an analogue of the dual semicanonical basis in the non simply-laced case, but only the set of cluster monomials which should be a part of it.
\end{rem}

\begin{demo}
 First, (\ref{partscan3}) is an easy consequence of (\ref{partscan1}) because cluster monomials of $\C[N]$ are of the form $\psi_{\bar X}$ where $X \in \mathcal{R}$. Moreover (\ref{partscan2}) is a particular case of (\ref{partscan3}). Let $\dg$ be such that $X \in \mathcal{R}_{\dg}$. Let us construct $f_X$ by induction on $\dg$. For $\dg = 0$, $f_0 = 1$. Suppose that $\dg \neq 0$. As $X$ is nilpotent, there exists $\ig$ such that $k_{\ig}(X) > 0$. Let us now argue by decreasing induction on $k_{\ig}(X)$. Suppose the result is proved for every $X'$ of dimension $\dg$ such that $k_{\ig}(X) < k_{\ig}(X') \leq \dg_{\ig}$. Let $f_\circ = f_{\Erg^\dag_{\ig}(\bar X)} e_{\ig}^{k_{\ig}(X)}$. For $Y \in \mathcal{R}$, one gets
  \begin{align*}
   \epsilon_{\bar Y}(f_\circ) &= \sum_{k \in \Z} k \chi\left(\left\{Z \in \Gr_{k_{\ig}(X)}(Y) \, | \, \udim Z = \ig^{k_{\ig}(X)} \text{ and } \epsilon_{\widebar {Y/Z}}\left( f_{\Erg^\dag_{\ig}(\bar X)}\right) = k \right\}\right) \\
   &= \left\{
    \begin{array}{ll}
     0 & \text{if } k_\ig(Y) < k_\ig(X) \\
     \epsilon_{\Erg_{\ig}^\dag(\bar Y)}\left( f_{\Erg^\dag_{\ig}(\bar X)}\right) & \text{if } k_\ig(Y) = k_\ig(X) \\
     \text{indeterminate} & \text{if } k_\ig(Y) > k_\ig(X)
    \end{array}
   \right. \\
   &= \left\{
    \begin{array}{ll}
     0 & \text{if } k_\ig(Y) < k_\ig(X) \\
     0 & \text{if } k_\ig(Y) = k_\ig(X) \text{ and } \Erg_{\ig}^\dag(Y) \nsim \Erg^\dag_{\ig}(X) \\
     1 & \text{if } k_\ig(Y) = k_\ig(X) \text{ and } \Erg_{\ig}^\dag(Y) \sim \Erg^\dag_{\ig}(X) \\
     \text{indeterminate} & \text{if } k_\ig(Y) > k_\ig(X)
    \end{array}
   \right. \\
   &= \left\{
    \begin{array}{ll}
     0 & \text{if } k_\ig(Y) < k_\ig(X) \\
     0 & \text{if } k_\ig(Y) = k_\ig(X) \text{ and } Y \nsim X \\
     1 & \text{if } Y \sim X \\
     \text{indeterminate} & \text{if } k_\ig(Y) > k_\ig(X).
    \end{array}
   \right.
  \end{align*}
  Hence
  $$f_{\bar X} = f_\circ - \sum_{\stackrel{\bar Y \in \mathcal{R}_\dg / \Gamma}{k_{\ig}(Y) > k_{\ig}(X)}} \epsilon_{\bar Y} (f_\circ) f_{\bar Y}$$
  is suitable (the $f_{\bar Y}$ in the sum exist by induction). \cqfd
\end{demo}

The second point of proposition \ref{partscan} together with corollary \ref{fmut2} proves that $((\psi_{\bar X}), B(\T))$ is the initial seed of a cluster algebra where $\bar X$ runs over the equivalence classes of indecomposable objects of $\T$.

Using \cite[9.3.2]{GeLeSc08}, one can now explicitly compute $B(\T)$. One denotes by $r$ the number of positive roots of $\Delta$ and $r_K$ the number of positive roots of $\Delta \restr{K/{\Gamma}}$.

Consider $(\ig)$ and $\ell$ defined as before. Let $\mathcal{I} = Q_0 / \Gamma \coprod \llbracket 1, \ell \rrbracket$ ordered in the following way: the order on $Q_0 / \Gamma$ does not matter, the order on $\llbracket 1, \ell \rrbracket$ is the natural one and, if $\jg \in Q_0/{\Gamma}$ and $n \in \llbracket 1, \ell \rrbracket$ then $\jg < n$. We now extend the word $(\ig)$ defined on $\llbracket 1, \ell \rrbracket$ to a word defined on $\mathcal{I}$ by seting, for $\jg \in Q_0/{\Gamma}$, $\ig_{\jg} = \jg$. Let $e((\ig)) = \{n \in \llbracket 1, \ell \rrbracket\,|\, \exists m \in \llbracket 1, \ell \rrbracket, m > n \text{ and } \ig_m = \ig_n \}$. If $n \in e((\ig)) \coprod Q_0/{\Gamma}$, then one denotes $n^+ = \min \{m \in \llbracket 1, \ell \rrbracket\,|\, m > n \text{ and } \ig_m = \ig_n\}$. For $\jg \in K/{\Gamma}$, let $t_{\jg} = \max \{t \leq r_K \,|\, \ig_t = \jg\}$. For $\jg \in J/{\Gamma}$, let $t_{\jg} = \jg$. Then, one constructs a matrix $B((\ig), K)$ whose lines are indexed by $(\rrbracket r_K, r \rrbracket \inter e((\ig))) \coprod \{t_{\jg}\,|\, \jg \in Q_0/{\Gamma}\}$ and columns by $\rrbracket r_K, r \rrbracket \inter e((\ig))$:
$$B((\ig), K)_{m n} = \left\{\begin{array}{ll}
 1 & \quad \text{if } m^+ = n  \\
 -1 & \quad \text{if } n^+ = m  \\
 -C_{\ig_m \ig_n} & \quad \text{if } n < m < n^+ < m^+ \\
 C_{\ig_m \ig_n} & \quad \text{if } m < n < m^+ < n^+ \\
 0 & \quad \text{else.}
\end{array}
\right.$$
 where $C$ is the Cartan matrix of $\Delta$.
 
\begin{ex}
 Continue with example \ref{exTi}. Then $e((\ig)) = \{1, 2, 3, 4, 5, 7 \}$. One computes $a^+ = 1$, $b^+ = 2$, $c^+ = 4$, $1^+ = 3$, $2^+ = 5$, $3^+ = 6$, $4^+ = 7$, $5^+ = 8$ and $7^+ = 9$. And also $t_a = 3$, $t_b = 2$ and $t_c = c$. One deduces that $B((\ig), K)$ has lines indexed by $\{4, 5, 7, 3, 2, c\}$ and columns indexed by $\{4, 5, 7\}$. The Cartan matrix of $\Delta$ is
 $$C = \left( \begin{matrix}
               2 & -1 & 0 \\
               -1 & 2 & -2 \\
               0 & -1 & 2
              \end{matrix} \right). $$
 Hence,
 $$B((\ig), K) = \left( \begin{matrix}
 		0 & -1 & 1 \\
		2 & 0 & -2 \\
		-1 & 1 & 0 \\
		0 & -1 & 0 \\
		-2 & 1 & 0 \\
		1 & 0 & 0		 
                \end{matrix} \right).$$
\end{ex}

\begin{prop}
 There is an indexation of the ${\Gamma}$-orbits of indecomposable direct summands of $T = T_{(\ig), K}$ by $\{t_n\,|\, n \in Q_0/{\Gamma}\} \coprod (\rrbracket r_K, r \rrbracket \inter e((\ig)))$ such that, via this identification,
 $$B(\T) = B((\ig), K).$$
\end{prop}

\begin{demo}
 Let $(i)$ be the image of $(\ig)$ by the substitution
 $$\ig \in Q_0 / {\Gamma} \mapsto \prod_{i \in \ig} i$$
 where the order of the letters in the product does not matter. Hence, one gets a representative of the longest element of $\tilde W$.
 
 First of all, according to \cite[9.3.2]{GeLeSc08}, $\tilde B(\T) = \tilde B((i), K)$ where $\tilde B((i), K)$ is the analogous of $B((\ig), K)$ for the action of the trivial group. To $n \in \llbracket 1, r \rrbracket$, one assigns the set $n' \subset \llbracket 1, \tilde r \rrbracket$ where $\tilde r$  is the number of positive roots of $\tilde \Delta$ in such a way that the set of letters at positions $n'$ in $(i)$ comes from the letter at position $n$ in $(\ig)$ through the above substitution. For $\jg \in Q_0/{\Gamma}$, let $\jg' = \jg$. Finally, for $n \in  Q_0 / {\Gamma} \coprod \llbracket 1, r \rrbracket $, let $n^\circ \in n'$. Then, it is enough to do the following computation using proposition \ref{relmut}, where $\tilde C$ is the Cartan matrix of $\tilde \Delta$:
 \begin{align*}
  B(\T)_{m  n} &= \sum_{\tilde m \in m'} \tilde B(\T)_{\tilde m n^\circ} = \sum_{\tilde m \in m'} \tilde B((i), K)_{\tilde m n^\circ} \\
  &= \sum_{\tilde m \in m'} \left\{\begin{array}{ll}
 1 & \quad \text{if } \tilde m^+ = n^\circ  \\
 -1 & \quad \text{if } n^{\circ+} = \tilde m  \\
 -\tilde C_{i_{\tilde m} i_{n^\circ}} & \quad \text{if } n^\circ < \tilde m < n^{\circ+} < \tilde m^+ \\
 \tilde C_{i_{\tilde m} i_{n^\circ}} & \quad \text{if } \tilde m < n^\circ  < \tilde m^+ < n^{\circ+} \\
 0 & \quad \text{else.}
 \end{array}
 \right. \\ 
 &= \left\{\begin{array}{ll}
  1 & \quad \text{if } m^+ = n  \\
 -1 & \quad \text{if } n^+ = m  \\
 \sum_{\tilde m \in m'} -\tilde C_{i_{\tilde m} i_{n^\circ}} & \quad \text{if } n < m < n^+ < m^+ \\
 \sum_{\tilde m \in m'} \tilde C_{i_{\tilde m} i_{n^\circ}} & \quad \text{if } m < n  < m^+ < n^+ \\
 0 & \quad \text{else.}
 \end{array}
 \right. \\ 
 &= \left\{\begin{array}{ll}
  1 & \quad \text{if } m^+ = n  \\
 -1 & \quad \text{if } n^+ = m  \\
 -C_{\ig_m \ig_n} & \quad \text{if } n < m < n^+ < m^+ \\
 C_{\ig_m \ig_n} & \quad \text{if } m < n  < m^+ < n^+ \\
 0 & \quad \text{else.}
 \end{array} \right.
 \end{align*}
 which is the expected result. \cqfd
\end{demo}

\begin{cor}
 The matrix $B(\T)$ is of full rank.
\end{cor}

\begin{demo}
 It is clear that for any column index $n$ of the matrix $B((\ig), K)$, there is a unique line index $n^-$ such that $(n^-)^+ = n$. Hence, for all column indices $n$, $B((\ig), K)_{n^- n} = 1$ by definition of $B((\ig), K)$. Moreover, if $m < n$ are two column indices, then $B((\ig), K)_{m^- n} = 0$ as $(m^-)^+ = m < n$. For summarize, the submatrix of $B((\ig), K)$ whose lines are the $n^-$ in the same order as the columns is lower triangular with diagonal $1$. \cqfd
\end{demo}

% \begin{prop}
%  Pour toute ${\Gamma}$-orbite $\bar X$ de classes d'isomorphisme de $\Sub Q_J$ telle que
%  $$\bigoplus_{X \in \bar X} X$$
%  soit rigide et qui est dans un amas que l'on peut atteindre par mutations à partir de $\T$, on a
%  $$\psi_{\bar X} = P_{\bar X}(\psi_{\bar T_i})_{\ig \in Q_0/{\Gamma}}.$$
% \end{prop}
% 
% \begin{demo}
%  C'est une conséquence du théorème \ref{fmut1} et du corollaire \ref{fmut2} par récurrence dans le graphe des mutations. \cqfd
% \end{demo}
% 
% On reprend les notations de la section \ref{echange}. Voici une réponse à une conjecture de Fomin et Zelevinsky dans ce cas:
% \begin{prop}
%  Les monômes d'amas de $\Ar(B(\T))$ sont linéairement indépendants.
% \end{prop}
% 
% \begin{demo}
%  On le déduit du corollaire \ref{indlin}. \cqfd
% \end{demo}

\begin{rems}
 \label{cnj}
 \begin{itemize}
  \item One conjectures that for every ${\Gamma}$-orbit $\bar X$ of isomorphism classes of $\Sub I_J$, one has
 $$\psi_{\bar X} = P_{\bar X}(\psi_{\bar T_i})_{i \in Q_0/{\Gamma}}.$$  
   This is clearly true for the rigid $\bar X$ which can be reached from $\T$ by mutations.
%  \item In simply-laced case, the conjecture of the first point holds according to an unpublished work of Schr\"oer.
  \item Using corollary \ref{indlin}, it gives another proof of the linear independence of cluster monomials of $\C[N]$ than proposition \ref{partscan}.
 \end{itemize}
\end{rems}

Let $G$ (resp. $\tilde G$) be the connected and simply-connected simple Lie group corresponding to the Lie algebra $\ggo$ (resp. $\tilde \ggo$). Let $B$ (resp. $\tilde B$) be its Borel subgroup corresponding to the Lie algebra $\bgo$ (resp. $\tilde \bgo$). Now, $N$ and $\tilde N$ are considered to be unipotent subgroups of $B$ and $\tilde B$. Let $\tilde B_K$ be the parabolic subgroup of $\tilde G$ generated by $\tilde B$ and the one-parameter subgroups $\tr{x_i}(t)$ for $i \in K$ and $t \in \C$. Let also $B_K$ be the parabolic subgroup of $G$ generated by $B$ and the images in $G$ of the one-parameters subgroups $\tr{x_i}(t)$ for $i \in K$ and $t \in \C$. Let $\tilde N_K$ be the unipotent radical of $\tilde B_K$. Let $N_K$ be the unipotent radical of $B_K$.

Let $\Ar'(\Sub I_J, {\Gamma}, \T)$ be the subalgebra of $\C[N]$ generated by the $\psi_{\bar X}$, where $\bar X$ runs over the ${\Gamma}$-orbits of isomorphism classes of $\Sub I_J$ such that
$$\bigoplus_{X \in \bar X} X$$
is rigid. 

Let $\Ar_0'(\Sub I_J, {\Gamma}, \T)$ be the subalgebra of $\C[N]$ generated by the $\psi_{\bar X}$, where $\bar X$ runs over the ${\Gamma}$-orbits of $\Sub I_J$. The algebras $\Ar(\Sub I_J, {\Gamma}, \T)$ and $\Ar_0(\Sub I_J, {\Gamma}, \T)$ were defined in section \ref{caramas}.

\begin{prop} \label{diagaa}
 If the conjecture of remark \ref{cnj} holds true then there is a commutative diagram
 $$\xymatrix{
  & \Ar(B(\T)) \ar@{^(->}[dl] \ar@{^(->}[dr] & \\
  \Ar(\Sub I_J, {\Gamma}, \T) \ar@{^(->}[d] \ar[rr]^{\sim} & & \Ar'(\Sub I_J, {\Gamma}, \T) \ar@{^(->}[d] \\
  \Ar_0(\Sub I_J, {\Gamma}, \T) \ar[rr]^{\sim} & & \Ar_0'(\Sub I_J, {\Gamma}, \T). \\
 }$$ 
\end{prop}

\begin{demo}
 It is clear because the $(\psi_{\bar T_i})_{\ig \in Q_0 / {\Gamma}}$ are algebraically independent. \cqfd
\end{demo}

\begin{rem}
 In proposition \ref{diagaa}, the four inclusions exist even if the conjecture of remark \ref{cnj} does not hold. It is a hard problem to understand when these inclusions are isomorphisms. 
\end{rem}

\begin{prop}
 One has $\C[N_K] \simeq \Ar_0'(\Sub I_J, {\Gamma}, \T)$.
\end{prop}

\begin{demo}
 This is the immediate translation of \cite[proposition 9.1]{GeLeSc08} together with the fact that $\kappa : \C[\tilde N]/\Gamma \simeq \C[N]$ defined just before notation \ref{defpsi} restricts to an isomorphism $\C[\tilde N_K]/\Gamma \simeq \C[N_K]$. \cqfd
\end{demo}

\begin{conj}
 One has $\C[N_K] \simeq \Ar(B(\T))$.
\end{conj}

\begin{lem}
 \label{rk}
 The clusters of $\Ar(B(\T))$ have
 \begin{enumerate}
  \item \label{rk1}$r - r_K$ cluster variables;
  \item \label{rk2}$\# Q_0 / {\Gamma}$ coefficients.
 \end{enumerate}
\end{lem}

\begin{demo}
 The point (\ref{rk2}) is proved using \cite[proposition 3.2]{GeLeSc08}. The proof of (\ref{rk1}) starts with the particular case $K = \emptyset$. In this case, the upper bound by $r$ is found as in \cite{GeSc05}, by counting the ${\Gamma}$-stable components of the module variety, and using the description of the roots of $\Delta$ of lemma \ref{racns}. Let now $T$ be a basic maximal $\Gamma$-stable rigid $\Lambda_Q$-module. As we have seen before, it is cluster tilting and therefore, according to \cite[proposition 7.3]{GeLeSc08}, it has $\tilde r$ indecomposable direct summands where $\tilde r$ is the number of positive roots of $\tilde \Delta$. One see that the $\Gamma$-orbits of these summands correspond to $\Gamma$-orbits of roots in the description of \cite{GeSc05}. Hence $T$ has exactly $r$ $\Gamma$-orbits of indecomposable summands using lemma \ref{racns}. After that, if $K \neq \emptyset$, the proof of \cite[proposition 7.1]{GeLeSc08} works exactly in the same way and therefore, there is at most $r - r_K$ cluster variables. The fact that $r - r_K$ is reached is the same as above for proving that $r$ is reached. \cqfd
\end{demo}

\renewcommand{\lmp}[1]
{\begin{minipage}{6cm}#1\end{minipage}}
One can now prove the following result, a part of which is proved in \cite[\S 11.4]{GeLeSc08} and the other part is conjectured in  \cite[\S 14.2]{GeLeSc08}:
\begin{prop}
 The cluster algebra $\Ar(B(\T))$ has a finite number of clusters exactly in the following cases (the circled vertices are those of $J$ and $n$ is the number of vertices):
\begin{center}
  \begin{longtable}{|c|c|}
    \hline
    Type of $G$ & Type of $\Ar(B(\T))$ \\
    \hline
   \endfirsthead
    \hline
    Type of $G$ & Type of $\Ar(B(\T))$ \\
    \hline
   \endhead
    \lmp{
    $$\xymatrix@R=.5cm@C=.5cm{
     \obullet \ar@{-}[r] & \bullet \ar@{-}[r] & \dots \ar@{-}[r] & \bullet
    }$$}
    & $A_0$ \\
    \hline
    \lmp{
    $$\xymatrix@R=.5cm@C=.5cm{
     \bullet \ar@{-}[r] & \obullet \ar@{-}[r] & \bullet \ar@{-}[r] & \dots \ar@{-}[r] & \bullet
    }$$}
    & $A_{n-2}$ \\
    \hline
    \lmp{
    $$\xymatrix@R=.5cm@C=.5cm{
     \obullet \ar@{-}[r] & \obullet \ar@{-}[r] & \bullet \ar@{-}[r] & \dots \ar@{-}[r] & \bullet
    }$$}
    & $A_{n-1}$ \\
    \hline
    \lmp{
    $$\xymatrix@R=.5cm@C=.5cm{
     \obullet \ar@{-}[r] & \bullet \ar@{-}[r] & \dots \ar@{-}[r] & \bullet \ar@{-}[r] & \obullet
    }$$}
    & $(A_1)^{n-1}$ \\
    \hline
    \lmp{
    $$\xymatrix@R=.5cm@C=.5cm{
     \obullet \ar@{-}[r] & \bullet \ar@{-}[r] & \dots \ar@{-}[r] & \bullet \ar@{-}[r] & \obullet \ar@{-}[r] & \bullet
    }$$}
    & $A_{2n-4}$ \\
    \hline
    \lmp{
    $$\xymatrix@R=.5cm@C=.5cm{
     \obullet \ar@{-}[r] & \obullet \ar@{-}[r] & \bullet \ar@{-}[r] & \dots \ar@{-}[r] & \bullet \ar@{-}[r] & \obullet
    }$$}
    & $A_{2n-3}$ \\
    \hline
    \lmp{
    $$\xymatrix@R=.5cm@C=.5cm{
     \obullet \ar@{-}[r] & \bullet \ar@{-}[r] & \dots \ar@{-}[r] & \bullet \ar@{-}[r] & \obullet \ar@{-}[r] & \bullet
    }$$}
    & $A_{2n-4}$ \\
    \hline
    \lmp{
    $$\xymatrix@R=.5cm@C=.5cm{
     \bullet \ar@{-}[r] & \obullet \ar@{-}[r] & \obullet \ar@{-}[r] & \bullet 
    }$$}
    & $D_4$ \\
    \hline
    \lmp{
    $$\xymatrix@R=.5cm@C=.5cm{
     \obullet \ar@{-}[r] & \obullet \ar@{-}[r] & \obullet \ar@{-}[r] & \bullet 
    }$$}
    & $D_5$ \\
    \hline
    \lmp{
    $$\xymatrix@R=.5cm@C=.5cm{
     \obullet \ar@{-}[r] & \obullet \ar@{-}[r] & \obullet \ar@{-}[r] & \obullet 
    }$$}
    & $D_6$ \\
    \hline
    \lmp{
    $$\xymatrix@R=.5cm@C=.5cm{
     \bullet \ar@{-}[r] & \bullet \ar@{-}[r] & \obullet \ar@{-}[r] & \bullet \ar@{-}[r] & \bullet 
    }$$}
    & $D_4$ \\
    \hline
    \lmp{
    $$\xymatrix@R=.5cm@C=.5cm{
     \obullet \ar@{-}[r] & \bullet \ar@{-}[r] & \obullet \ar@{-}[r] & \bullet \ar@{-}[r] & \bullet 
    }$$}
    & $E_6$ \\
    \hline
    \lmp{
    $$\xymatrix@R=.5cm@C=.5cm{
     \bullet \ar@{-}[r] & \obullet \ar@{-}[r] & \obullet \ar@{-}[r] & \bullet \ar@{-}[r] & \bullet 
    }$$}
    & $E_6$ \\
    \hline
    \lmp{
    $$\xymatrix@R=.5cm@C=.5cm{
     \obullet \ar@{-}[r] & \obullet \ar@{-}[r] & \obullet \ar@{-}[r] & \bullet \ar@{-}[r] & \bullet 
    }$$}
    & $E_7$ \\
    \hline
    \lmp{
    $$\xymatrix@R=.5cm@C=.5cm{
     \bullet \ar@{-}[r] & \bullet \ar@{-}[r] & \obullet \ar@{-}[r] & \bullet \ar@{-}[r] & \bullet \ar@{-}[r] & \bullet
    }$$}
    & $E_6$ \\
    \hline
    \lmp{
    $$\xymatrix@R=.5cm@C=.5cm{
     \bullet \ar@{-}[r] & \obullet \ar@{-}[r] & \obullet \ar@{-}[r] & \bullet \ar@{-}[r] & \bullet \ar@{-}[r] & \bullet
    }$$}
    & $E_8$ \\
    \hline
    \lmp{
    $$\xymatrix@R=.5cm@C=.5cm{
     \bullet \ar@{-}[r] & \bullet \ar@{-}[r] & \obullet \ar@{-}[r] & \bullet \ar@{-}[r] & \bullet \ar@{-}[r] & \bullet \ar@{-}[r] & \bullet
    }$$}
    & $E_8$ \\
    \hline
    \lmp{
    $$\xymatrix@R=.5cm@C=.5cm{
     & & & & \bullet \\
     \obullet \ar@{-}[r] & \bullet \ar@{-}[r] & \dots \ar@{-}[r] & \bullet \ar@{-}[ur] \ar@{-}[dr] &  \\
     & & & & \bullet
    }$$}
    & $(A_1)^{n-2}$ \\
    \hline
    \lmp{
    $$\xymatrix@R=.5cm@C=.5cm{
     & & \obullet \\
     \bullet \ar@{-}[r] & \bullet \ar@{-}[ur] \ar@{-}[dr] &  \\
     & & \obullet
    }$$}
    & $A_5$ \\
    \hline
    \lmp{
    $$\xymatrix@R=.5cm@C=.5cm{
     & & & \obullet \\
     \bullet \ar@{-}[r] & \bullet \ar@{-}[r] & \bullet \ar@{-}[ur] \ar@{-}[dr] &  \\
     & & & \bullet
    }$$}
    & $A_5$ \\
    \hline
    \lmp{$$\xymatrix@R=.5cm@C=.5cm{
     \obullet \ar@{-}[r] & \bullet \ar@{-}[r] & \dots \ar@{-}[r] & \bullet \ar@{=}[rr] & \text{\LARGE <} & \bullet 
    }$$}
    & $(A_1)^{n-1}$ \\
    \hline
    \lmp{$$\xymatrix@R=.5cm@C=.5cm{
     \obullet \ar@{-}[r] & \bullet \ar@{-}[r] & \dots \ar@{-}[r] & \bullet \ar@{=}[rr] & \text{\LARGE >} & \bullet 
    }$$}
    & $(A_1)^{n-1}$ \\
    \hline
    \lmp{$$\xymatrix@R=.5cm@C=.5cm{
     \obullet \ar@{=}[rr] & \text{\LARGE >} & \obullet  
    }$$}
    & $B_2 = C_2$ \\
    \hline
    \lmp{$$\xymatrix@R=.5cm@C=.5cm{
     \bullet \ar@{-}[r] & \bullet \ar@{=}[rr] & \text{\LARGE <} & \obullet 
    }$$}
    & $B_3$ \\
    \hline
    \lmp{$$\xymatrix@R=.5cm@C=.5cm{
     \bullet \ar@{-}[r] & \bullet \ar@{=}[rr] & \text{\LARGE >} & \obullet 
    }$$}
    & $C_3$ \\
    \hline
  \end{longtable}
 \end{center}
\end{prop}

\begin{demo}
 All simply laced cases are proved in \cite[11.4]{GeLeSc08}. The other cases must come from a simply laced case endowed with a group action by proposition \ref{typefini}. Thus, one has to look at the automorphisms of each diagram stabilizing $J$. This gives immediately a list of five non simply-laced cases:
  \begin{center}
   \begin{longtable}{|c|c|c|}
    \hline
    \mpc{Type of $\tilde G$} & \mpc{$\Gamma$} & \mpc{Type of $G$} \\
    \hline
    \endfirsthead
    \hline
    \mpc{Type of $\tilde G$} & \mpc{$\Gamma$} & \mpc{Type of $G$} \\
    \hline
    \endhead
    \lmp{
    $$\xymatrix@R=.5cm@C=.5cm{
     \obullet \ar@{-}[r] & \bullet \ar@{-}[r] & \dots \ar@{-}[r] & \bullet \ar@{-}[r] & \obullet
    }$$}
    & $\Z/2\Z$ &
    \lmp{$$\xymatrix@R=.5cm@C=.5cm{
     \obullet \ar@{-}[r] & \bullet \ar@{-}[r] & \dots \ar@{-}[r] & \bullet \ar@{=}[rr] & \text{\LARGE <} & \bullet 
    }$$}    
    \\
    \hline
    \lmp{
    $$\xymatrix@R=.5cm@C=.5cm{
     \obullet \ar@{-}[r] & \obullet \ar@{-}[r] & \obullet
    }$$} & $\Z/2\Z$ &
    \lmp{$$\xymatrix@R=.5cm@C=.5cm{
     \obullet \ar@{=}[rr] & \text{\LARGE >} & \obullet  
    }$$} \\
    \hline
    \lmp{
    $$\xymatrix@R=.5cm@C=.5cm{
     \bullet \ar@{-}[r] & \bullet \ar@{-}[r] & \obullet \ar@{-}[r] & \bullet \ar@{-}[r] & \bullet 
    }$$}
    & $\Z/2\Z$ &
    \lmp{$$\xymatrix@R=.5cm@C=.5cm{
     \bullet \ar@{-}[r] & \bullet \ar@{=}[rr] & \text{\LARGE <} & \obullet 
    }$$}
    \\
    \hline
    \lmp{
    $$\xymatrix@R=.5cm@C=.5cm{
     & & & & \bullet \\
     \obullet \ar@{-}[r] & \bullet \ar@{-}[r] & \dots \ar@{-}[r] & \bullet \ar@{-}[ur] \ar@{-}[dr] &  \\
     & & & & \bullet
    }$$}
    & $\Z/2\Z$ &
    \lmp{$$\xymatrix@R=.5cm@C=.5cm{
     \obullet \ar@{-}[r] & \bullet \ar@{-}[r] & \dots \ar@{-}[r] & \bullet \ar@{=}[rr] & \text{\LARGE >} & \bullet 
    }$$}
    \\
    \hline
    \lmp{
    $$\xymatrix@R=.5cm@C=.5cm{
     & & \obullet \\
     \bullet \ar@{-}[r] & \bullet \ar@{-}[ur] \ar@{-}[dr] &  \\
     & & \obullet
    }$$}
    & $\Z/2\Z$ &
    \lmp{$$\xymatrix@R=.5cm@C=.5cm{
     \bullet \ar@{-}[r] & \bullet \ar@{=}[rr] & \text{\LARGE >} & \obullet 
    }$$}    
    \\
    \hline
  \end{longtable} 
  \end{center}
 
  Compute their cluster type:
%  \begin{itemize}
%   \item
  the diagram 
  $$\xymatrix@R=.5cm@C=.5cm{
     \obullet \ar@{-}[r] & \bullet \ar@{-}[r] & \dots \ar@{-}[r] & \bullet \ar@{=}[rr] & \text{\LARGE <} & \bullet 
  }$$
  with $n$ vertices comes from the diagram
  $$\xymatrix@R=.5cm@C=.5cm{
   \obullet \ar@{-}[r] & \bullet \ar@{-}[r] & \dots \ar@{-}[r] & \bullet \ar@{-}[r] & \obullet
  }$$
  with $2n-1$ vertices endowed with the only non-trivial automorphism of order $2$. Hence, its type must be obtained from $(A_1)^{2n-2}$ with the action of $\Z/2\Z$; thus, its cluster type is of the form $(A_1)^k$ for some $k$. Using lemma \ref{rk}, the number of cluster variables in a cluster is $r - r_K = n^2 - (n-1)^2 = 2n-1$ and if one removes the $n$ coefficients, its type has to be of rank $n-1$ which implies the its cluster type is $(A_1)^{n-1}$.
   The other cases can be handled by the same method. \cqfd
\end{demo}

\subsection{Categories $\Cr_M$ and unipotent groups}

\label{cm}
This application is a generalization of \cite{GeLeSc}. Let $Q$ be now an arbitrary quiver without oriented cycles. Let ${\Gamma}$ be a group acting on $Q$ in an admissible way (see previous section). The algebra $kQ$ is naturally identified with a subalgebra of $\Lambda_Q$, one denotes by
$$\pi_Q: \md \Lambda_Q \rightarrow \md k Q$$
the corresponding restriction functor. It is essentially surjective.

\begin{df}
 A module $M \in \md kQ$ is said to be \emph{terminal} if
 \begin{enumerate}
  \item $M$ is preinjective; 
  \item if $X \in \md kQ$ is indecomposable and $\Hom_{k Q}(M, X) \neq 0$, then $X \in \add(M)$;
  \item $\add(M)$ contains all injective $kQ$-modules.
 \end{enumerate}
\end{df}

\begin{df}
 Let $M \in \md kQ$ be a terminal module. Define
 $$\Cr_M = \pi_Q^{-1}(\add(M)).$$
\end{df}

\begin{thm}{\citeb{theorem 2.1}{GeLeSc}}
 Let $M \in \md kQ$ be a terminal module. Then the category $\Cr_M$ is an exact, $\Hom$-finite, Krull-Schmidt, Frobenius, and $2$-Calabi-Yau subcategory of $\md \Lambda_Q$.
\end{thm}

Let $M$ be a terminal, ${\Gamma}$-stable $\md kQ$-module.

Recall the following lemma of Gei\ss, Leclerc and Schröer:

\begin{lem}[\citeb{lemma 5.6}{GeLeSc}]
 \label{stfact}
 The category $\Cr_M$ is a subcategory of $\md \Lambda_Q$ stable by factors. In other terms, for $X \in \Cr_M$ and $Y$ a $\Lambda_Q$-submodule of $X$, then $X / Y \in \Cr_M$.
\end{lem}

\begin{cor}
 All projective objects of $\Cr_M$ have left rigid quasi-approximations.
\end{cor}

\begin{demo}
 In order to prove this, it is enough to see that the hypothesis of lemma \ref{qapsc} are satisfied. The lemmas at the beginning of \cite[\S 8]{GeLeSc} prove that the projective objects of $\Cr_M$ have simple socles in $\md \Lambda_Q$. Moreover, as the action of ${\Gamma}$ on $Q$ is admissible, $\bigoplus_{g \in {\Gamma}} \g \tens S$ is rigid for all simple $\Lambda_Q$-modules $S$. The other hypothesis of lemma \ref{qapsc} are immediate consequences of lemma \ref{stfact}. \cqfd
\end{demo}

The $\Lambda_Q$-modules $T_M$ and $T_M^\vee$ constructed in \cite[\S 7]{GeLeSc} are cluster-tilting and $\End_{\Cr_M}(T_M)$ and $\End_{\Cr_M}(T_M^\vee)$ are of finite global dimension. Moreover, they are ${\Gamma}$-stable as $M$ is. 

The action of ${\Gamma}$ is $2$-Calabi-Yau in the sense of definition \ref{act2CY} for the same reasons as in the previous section. Hence, one can apply the results of section \ref{caramas} in $\Cr_M$ with $\T = \add(T_M)$ and $\T^\vee = \add(T_M^\vee)$. 

Let $\Theta$ be the Gabriel quiver of $\End_{kQ}(M)$, in which one adds an arrow $x \rightarrow \tau(x)$ for every vertex $x$ such that $\tau(x)$ correspond to an indecomposable object of $\add(M)$ (where $\tau$ is the Auslander-Reiten translation). A vertex $i$ of $\Theta$ is called to be \emph{frozen} if $\tau(i) \notin \Theta$.

\begin{prop}[\citeb{\S 7.2}{GeLeSc}]
 The matrices $\tilde B(\T)$ and $\tilde B(\T^\vee)$ are equal, with an appropriate indexation, to the adjacency matrix of $\Theta$, from which the columns corresponding to frozen vertices are removed.
\end{prop}

\begin{nt}
 Let $\N Q^{\op}_0 = \N \times Q_0$ and $\N Q^{\op}_1 = \N \times Q^{\op}_1 \coprod \N \times Q_1$ endowed with maps $s, t: \N Q^{\op}_1 \rightarrow \N Q^{\op}_0$ defined by 
 \begin{align*}
 s(n, q) &= \left\{\begin{array}{ll} (n, s(q)) & \quad \text{if } (n, q) \in \N \times Q^{\op}_1 \\ (n+1, s(q)) & \quad \text{if } (n, q) \in \N \times Q_1 \end{array} \right. \\
 t(n, q) &= (n, t(q)). \end{align*}
 Thus, one defined a quiver $\N Q^{\op}$. Let also $\tilde \Theta$ be the quiver $\N Q^{\op}$, on which one adds the arrows $(n, i) \rightarrow (n+1, i)$ for $(n, i) \in \N Q^{\op}_0$.
 
 As $Q$ has no cycles, an order on $Q_0$ can be fixed in such a way that every arrow of $Q$ has a larger target than its source and in such a way that ${\Gamma}$ acts on $Q_0$ by increasing maps. Thus, this order induces an order on $Q_0 / {\Gamma}$. Then, one endows $\N Q^{\op}_0 = \N \times Q^{\op}$ with the lexicographic order.
\end{nt}

It is classical that the Auslander-Reiten quiver of $\add(M)$ is a subquiver of $\N Q^{\op}$ and that $\Theta$ is a subquiver of $\tilde \Theta$, the Auslander-Reiten translation $\tau$ of $\add(M)$ being given by the arrows $(n, i) \rightarrow (n+1, i)$ in $\tilde \Theta$.

\begin{cor}
 The matrix $B(\T)$ has full rank. 
\end{cor}

\begin{demo}
 One will use the structure of $\tilde \Theta$. For every non-frozen ${\Gamma}$-orbit $\bar X$ of $\Theta$, $\tau(\bar X)$ is a ${\Gamma}$-orbit of $\Theta$. One restricts the order on vertices of $\tilde \Theta$ to an order on vertices of $\Theta$. Then, one gets for every non-frozen ${\Gamma}$-orbit $\bar X$ of $\Theta$, $B(\T)_{\tau(\bar X), \bar X} = -1$ and, if one considers another non-frozen ${\Gamma}$-orbit $\bar Y < \bar X$ of $\Theta$, $B(\T)_{\tau(\bar X), \bar Y} = 0$ by construction of the order on $\tilde \Theta$. In other words, the submatrix of $B(\T)$ whose lines are the $\tau(\bar X)$ is upper triangular of diagonal $-1$. \cqfd
\end{demo}

\begin{cor}
 The cluster monomial of $\Ar(\Cr_M, {\Gamma}, \T)$ are linearly independent. 
\end{cor}

\begin{demo}
 It is a direct application of corollary \ref{indlin}. \cqfd
\end{demo}

\begin{thm} \label{acyc}
 For every acyclic cluster algebra without coefficient $\Ar$, there is a quiver $Q$, a finite group ${\Gamma}$ acting on $Q$, a terminal module $M$ of $\md k Q$ and a cluster-tilting subcategory $\T$ of $\Cr_M$ such that the cluster algebra $\Ar(\Cr_M, {\Gamma}, \T)$ with coefficients specialized to $1$ is isomorphic to $\Ar$. This holds in particular for cluster algebras of finite type.
\end{thm}

\begin{demo}
 Using the same proof as for lemma \ref{conssym} for an acyclic exchange matrix $B$ of $\Ar$, there is a skew-symmetric matrix $\tilde B$ and an action of a finite group ${\Gamma}$ on it, such that $B$ is given by proposition \ref{relmut} (\ref{relmut2}). It is easy to see that $\tilde B$ is acyclic. Let $Q$ be the quiver of adjacency matrix $\tilde B$. If $\md k Q$ has no projective-injective object (i.e., $Q$ is not of type $A_n$ oriented in one direction), then, $M = kQ^{\op} \oplus \tau(kQ^{\op})$ is suitable. If $Q$ is of type $A_n$ oriented in one direction, $\Gamma$ acts trivially. In this case, let $Q'$ be a quiver of type $A_{n+1}$. Then $Q'$ with $M = kQ'^{\op} \oplus \tau(kQ^{\op})$ is suitable (where $Q$ is considered to be a subquiver of $Q'$ with the same source). \cqfd
\end{demo}

Here are two results which are immediate consequences of \cite{GeLeSc}.

\begin{prop}
 There exists in $\Ar(\Cr_M, {\Gamma}, \T)$ a sequence of mutations going from $\T$ to $\T^\vee$.
\end{prop}

\begin{demo}
 Gei\ss, Leclerc and Schröer described in \cite[\S 18]{GeLeSc} a sequence of mutations going from $\T$ to $\T^\vee$ in $\Ar(\Cr_M, \T)$. It is easy to see, following their algorithm, that one can sort these mutations in such a way that the mutations of a ${\Gamma}$-orbit of $\Cr$ are consecutive (it is enough to permute mutations which commute). Then, by proposition \ref{relmut}, it is clear that this sequence of mutations comes from a sequence of mutations in $\Ar(\Cr_M, {\Gamma}, \T)$. \cqfd
\end{demo}

\begin{prop}
 \label{cluspoly}
 The algebra $\Ar(B(\T))$ is a polynomial ring. More precisely,
 $$\Ar(B(\T)) = \C[P_{\bar X }]_{\bar X \in \tilde M}$$
 where $\tilde M$ is the set of ${\Gamma}$-orbits of isomorphism classes of indecomposable objects of $\add(M)$.
\end{prop}

\begin{demo}
 By corollary \ref{qalga}, $\Ar(B(\T)) \subset \Ar(\tilde B(\T))/{\Gamma}$ so, using \cite[theorem 3.4]{GeLeSc}, there is an inclusion $\Ar(B(\T)) \subset \C[P_{\bar X }]_{\bar X \in \tilde M}$. For the converse inclusion, one uses the same technique as in \cite[\S 20.2]{GeLeSc}: every $\bar X \in \tilde M$ appears in the sequence of mutations of the previous proposition. \cqfd
\end{demo}

The end of this section deals with the case of general Kac-Moody groups. It works in par\-ti\-cu\-lar for semisimple Lie groups. For more details about the infinite dimensinal case, in particular about the constructions of the finite dimensional subgroups $N(w)$ and $N^w$, one refers to \cite{GeLeSc}. One retains the notation of section \ref{rac}. Let $N$ (resp. $\tilde N$) be the pro-unipotent pro-group defined from $\ngo$ (resp. $\tilde \ngo$) as in \cite[\S 22]{GeLeSc} (see also \cite[\S 4.4]{Ku02}) in such a way that $\C[N] \simeq U(\ngo)^*_{\gr}$ and $\C[\tilde N] \simeq U(\tilde \ngo)^*_{\gr}$.  

One denotes by $\tilde \Delta^+_M$ the set of dimension vectors of indecomposable direct summands of $M$. As $M$ is stable under the action of ${\Gamma}$, $\tilde \Delta^+_M$ also and therefore, one can denote by $\Delta^+_M$ the image of $\tilde \Delta^+_M$ in the real part of a root system of type $\Delta$. As in \cite[3.7]{GeLeSc}, there exists a unique $w$ in the Weyl group $W$ of $\Delta$ such that $\Delta^+_M = \{\alpha \in \Delta^+ \,|\, w(\alpha) < 0\}$ where $\Delta^+$ is the subset of real positive roots of $\Delta$. Let $\tilde w$ be the image of $w$ in $\tilde W$.

The subalgebras $\ngo(w)$ (resp. $\tilde \ngo(w)$) of $\ngo$ (resp. $\tilde \ngo$) are defined as in \cite[\S 19.3]{GeLeSc}. One retains the definitions of the finite dimensional subgroups $N(w)$ and $N^w$ (resp. $\tilde N(w)$ and $\tilde N^w$) of $N$ (resp. $\tilde N$) given in \cite[\S 22]{GeLeSc}. Then $\C[N(w)] \simeq U(\ngo(w))^*_{\gr}$ and $\C[\tilde N(w)] \simeq U(\tilde \ngo(w))^*_{\gr}$. In the particular case where $Q$ is of Dynkin type, we are in the classical Lie framework. In this case, let $G$ (resp. $\tilde G$) be the connected and simply-connected Lie group associated to $\ggo$ (resp. $\tilde \ggo$). Then
$$N(w) = N \inter (w^{-1} N_- w) \quad \text{and} \quad N^w = N \inter (B_- w B_-)\;;$$
$$\tilde N(\tilde w) = \tilde N \inter (\tilde w^{-1} \tilde N_- \tilde w) \quad \text{and} \quad \tilde N^{\tilde w} = \tilde N \inter (\tilde B_- \tilde w \tilde B_-)$$
where $B_-$ (resp. $\tilde B_-$) denotes the Borel subgroup of $G$ (resp. $\tilde G$) associated to $\bgo_-$ (resp. $\tilde \bgo_-$) and $N_-$ (resp. $\tilde N_-$) denotes the unipotent subgroup of $G$ (resp. $\tilde G$) associated to $\ngo_-$ (resp. $\tilde \ngo_-$).

Here is the analogous of \cite[theorem 3.5]{GeLeSc}:
\begin{thm}
 The cluster algebra $\Ar(B(\T))$ is a cluster algebra structure on $\C[N(w)]$. The cluster algebra $\tilde \Ar(B(\T))$ obtained by  inversing the coefficients is a cluster algebra structure on $\C[N^w]$.
\end{thm}

\begin{demo}
 According to \cite[theorem 3.5]{GeLeSc}, $\Ar(\tilde B(\T))$ is a cluster algebra structure on $\C[\tilde N(\tilde w)]$. It is easy to see that the epimorphism $\kappa : U(\tilde \ngo)^*_{\gr}/\Gamma \surj U(\ngo)^*_{\gr}$ (see section \ref{rac}) restricts to a morphism $U(\tilde \ngo(\tilde w))^*_{\gr}/\Gamma \rightarrow U(\ngo(w))^*_{\gr}$. Hence, one can complete the following commutative diagram :
 $$\xymatrix{
  \Ar(B(\T)) \ar@{^(->}[r] \ar@{-->}[d]_\alpha & \Ar(\tilde B(\T))/\Gamma \ar[r]^\sim & \C[\tilde N(\tilde w)]/\Gamma \ar[r]^\sim & U(\tilde \ngo(\tilde w))^*_{\gr}/\Gamma \ar[d] \ar@{^(->}[r] & U(\tilde \ngo)^*_{\gr}/\Gamma \ar@{->>}[d]_\kappa \\
  \C[N(w)] \ar[rrr]^\sim & & & U(\ngo(w))^*_{\gr} \ar@{^(->}[r] & U(\ngo)^*_{\gr} 
 }$$
 Moreover, using an easy adaptation of proposition \ref{partscan}, one gets that $\alpha$ is a monomorphism. Using proposition \ref{cluspoly}, one concludes that $\alpha$ is an epimorphism, because, by the same method as in \cite[proposition 22.2]{GeLeSc}, $\C[P_{\bar X }]_{\bar X \in \tilde M}$ is the whole algebra $\C[N(w)]$.
 
 The case of $\C[N^w]$ is handled by the same method. \cqfd
\end{demo}

\newcommand{\QA}{ 
\begin{minipage}{12.8pt}\tiny
$$\xymatrix@R=.2cm@C=.2cm{
 a
}$$
\end{minipage}}

\newcommand{\QAP}{ 
\begin{minipage}{76.2pt}\tiny
$$\xymatrix@R=.2cm@C=.2cm{
 a \ar[dr] & & & & \\
 & b \ar[dr]& & & \\
 & & c \ar[dr] & & \\
 & & & b' \ar[dr] & \\
 & & & & a'
}$$
\end{minipage}}

\newcommand{\QB}{ 
\begin{minipage}{12.8pt}\tiny
$$\xymatrix@R=.2cm@C=.2cm{
 b
}$$
\end{minipage}}

\newcommand{\QBP}{ 
\begin{minipage}{60.5pt}\tiny
$$\xymatrix@R=.2cm@C=.2cm{
 b \ar[dr]& & & \\
 & c \ar[dr] & & \\
 & & b' \ar[dr] & \\
 & & & a'
}$$
\end{minipage}}

\newcommand{\QC}{ 
\begin{minipage}{45.5pt}\tiny
$$\xymatrix@R=.2cm@C=.2cm{
 c \ar[dr] & & \\
 & b' \ar[dr] & \\
 & & a'
}$$
\end{minipage}}

\newcommand{\QD}{ 
\begin{minipage}{30.5pt}\tiny
$$\xymatrix@R=.2cm@C=.2cm{
 & b \ar[dl] \\
 a & 
}$$
\end{minipage}}

\newcommand{\QDP}{ 
\begin{minipage}{76.2pt}\tiny
$$\xymatrix@R=.2cm@C=.2cm{
 & b \ar[dl] \ar[dr] & & & \\
 a \ar[dr] & & c \ar[dl] \ar[dr] & & \\
 & b \ar[dr] & & b' \ar[dl] \ar[dr] & \\
 & & c \ar[dr] & & a' \ar[dl] \\
 & & & b' &
}$$
\end{minipage}}

\newcommand{\QE}{ 
\begin{minipage}{60.5pt}\tiny
$$\xymatrix@R=.2cm@C=.2cm{
 & c \ar[dl] \ar[dr] & & \\
 b & & b' \ar[dr] & \\
 & & & a' 
}$$
\end{minipage}}

\newcommand{\QEP}{ 
\begin{minipage}{60.5pt}\tiny
$$\xymatrix@R=.2cm@C=.2cm{
 & c \ar[dl] \ar[dr] & & \\
 b \ar[dr] & & b' \ar[dl] \ar[dr] & \\
 & c \ar[dr] & & a' \ar[dl] \\
 & & b' &
}$$
\end{minipage}}

\newcommand{\QF}{ 
\begin{minipage}{30.5pt}\tiny
$$\xymatrix@R=.2cm@C=.2cm{
 c \ar[dr] & \\
 & b' 
}$$
\end{minipage}}

\newcommand{\QG}{ 
\begin{minipage}{76.2pt}\tiny
$$\xymatrix@R=.2cm@C=.2cm{
 & & c \ar[dl] \ar[dr] & & \\
 & b \ar[dl] & & b' \ar[dr] & \\
 a & & & & a'
}$$
\end{minipage}}

\newcommand{\QGP}{ 
\begin{minipage}{76.2pt}\tiny
$$\xymatrix@R=.2cm@C=.2cm{
 & & c \ar[dl] \ar[dr] & & \\
 & b \ar[dl] \ar[dr] & & b' \ar[dr] \ar[dl] & \\
 a \ar[dr] & & c \ar[dl] \ar[dr] & & a' \ar[dl] \\
 & b \ar[dr] & & b' \ar[dl] & \\
 & & c & &
}$$
\end{minipage}}

\newcommand{\QH}{ 
\begin{minipage}{45.5pt}\tiny
$$\xymatrix@R=.2cm@C=.2cm{
 & c \ar[dl] \ar[dr] &\\
 b & & b'
 }$$
\end{minipage}}

\newcommand{\QHP}{ 
\begin{minipage}{45.5pt}\tiny
$$\xymatrix@R=.2cm@C=.2cm{
 & c \ar[dl] \ar[dr] &\\
 b \ar[dr] & & b' \ar[dl] \\
 & c &
 }$$
\end{minipage}}

\newcommand{\QI}{ 
\begin{minipage}{12.8pt}\tiny
$$\xymatrix@R=.2cm@C=.2cm{
 c
}$$
\end{minipage}}

\newcommand{\QJ}{ 
\begin{minipage}{30.5pt}\tiny
$$\xymatrix@R=.2cm@C=.2cm{
 b' \ar[dr] & \\
 & a'
}$$
\end{minipage}}

\newcommand{\QJP}{ 
\begin{minipage}{76.2pt}\tiny
$$\xymatrix@R=.2cm@C=.2cm{
 & & & b' \ar[dl] \ar[dr] & \\
 & & c \ar[dl] \ar[dr] & & a' \ar[dl]  \\
 & b \ar[dl] \ar[dr] & & b' \ar[dl] & \\
 a \ar[dr] & & c \ar[dl]& & \\
  & b & & &
}$$
\end{minipage}}

\newcommand{\QK}{ 
\begin{minipage}{60.5pt}\tiny
$$\xymatrix@R=.2cm@C=.2cm{
 & & c \ar[dl] \ar[dr] & \\
 & b \ar[dl] & & b' \\
 a & & & 
}$$
\end{minipage}}

\newcommand{\QKP}{ 
\begin{minipage}{60.5pt}\tiny
$$\xymatrix@R=.2cm@C=.2cm{
 & & c \ar[dl] \ar[dr] & \\
 & b \ar[dl] \ar[dr] & & b' \ar[dl] \\
 a \ar[dr] & & c \ar[dl]& \\
  & b & & 
}$$
\end{minipage}}

\newcommand{\QL}{ 
\begin{minipage}{30.5pt}\tiny
$$\xymatrix@R=.2cm@C=.2cm{
 & c \ar[dl] \\
 b &
}$$
\end{minipage}}

\newcommand{\QM}{ 
\begin{minipage}{12.8pt}\tiny
$$\xymatrix@R=.2cm@C=.2cm{
 a'
}$$
\end{minipage}}

\newcommand{\QMP}{ 
\begin{minipage}{76.2pt}\tiny
$$\xymatrix@R=.2cm@C=.2cm{
 & & & & a' \ar[dl] \\
 & & & b' \ar[dl] & \\
 & & c \ar[dl] & & \\
 & b \ar[dl] & & & \\
 a & & & &
}$$
\end{minipage}}

\newcommand{\QN}{ 
\begin{minipage}{12.8pt}\tiny
$$\xymatrix@R=.2cm@C=.2cm{
 b'
}$$
\end{minipage}}

\newcommand{\QNP}{ 
\begin{minipage}{60.5pt}\tiny
$$\xymatrix@R=.2cm@C=.2cm{
 & & & b' \ar[dl] \\
 & & c \ar[dl] &\\
 & b \ar[dl] & & \\
 a & & &
}$$
\end{minipage}}

\newcommand{\QO}{ 
\begin{minipage}{45.5pt}\tiny
$$\xymatrix@R=.2cm@C=.2cm{
 & & c \ar[dl]  \\
 & b \ar[dl] &  \\
 a & & 
}$$
\end{minipage}}

\begin{figure}
\begin{center}
 \setlength{\unitlength}{1pt}
 \begin{picture}(373,286)
  \put(0,262.8){\QA}
  \put(120,262.8){\QB}
  \put(223.65,262.8){\QC}
  \put(51.15,202.8){\QD}
  \put(156.15,202.8){\QE}
  \put(291.15,202.8){\QF}
  \put(88.3,142.8){\QG}
  \put(223.65,142.8){\QH}
  \put(360,142.8){\QI}
  \put(51.15,82.8){\QJ}
  \put(156.15,82.8){\QK}
  \put(291.15,82.8){\QL}
  \put(0,22.8){\QM}
  \put(120,22.8){\QN}
  \put(223.65,22.8){\QO}
  \put(275.15,235.05){\vector(1,-1){17}}
  \put(215.15,175.05){\vector(1,-1){17}}
  \put(155.15,115.05){\vector(1,-1){17}}
  \put(95.15,55.05){\vector(1,-1){17}}
  \put(275.15,45.55){\vector(1,1){17}}
  \put(215.15,105.55){\vector(1,1){17}}
  \put(155.15,165.55){\vector(1,1){17}}
  \put(95.15,225.55){\vector(1,1){17}}
  \put(328.15,98.55){\vector(1,1){31}}
  \put(268.15,158.55){\vector(1,1){31}}
  \put(208.15,218.55){\vector(1,1){31}}
  \put(328.15,182.05){\vector(1,-1){31}}
  \put(268.15,122.05){\vector(1,-1){31}}
  \put(208.15,62.05){\vector(1,-1){31}}
  \put(140.15,250.05){\vector(1,-1){27}}
  \put(80.15,190.05){\vector(1,-1){27}}
  \put(140.15,30.55){\vector(1,1){27}}
  \put(80.15,90.55){\vector(1,1){27}}
  \put(15.15,255.05){\vector(1,-1){46}}
  \put(15.15,25.55){\vector(1,1){46}}
  \put(66.4,17.8){\line(1,1){60}}
  \put(126.4,77.8){\line(-1,1){60}}
  \put(66.4,137.8){\line(1,1){60}}
  \put(126.4,197.8){\line(-1,1){60}}
  \put(69.4,17.8){\line(1,1){60}}
  \put(129.4,77.8){\line(-1,1){60}}
  \put(69.4,137.8){\line(1,1){60}}
  \put(129.4,197.8){\line(-1,1){60}}
 \end{picture}
 \end{center}
 \caption{Auslander-Reiten quiver of $Q$}
 \label{ARkQ}
\end{figure}

\begin{ex}
 In this example, one will denote by $Q$ the quiver
 $$\xymatrix{
 a & b \ar[l] & c \ar[l] \ar[r] & b' \ar[r] & a'
 }$$
 endowed with the non-trivial action of ${\Gamma} = \Z/2\Z$. The Auslander-Reiten quiver of $\md \C Q$ is displayed in figure \ref{ARkQ}. Denote by $M_0$ the direct sum of the indecomposable $\C Q$-modules so that $\Cr_{M_0} = \md \Lambda_Q$. Let $M$ be the direct sum of the indecomposable $\C Q$-modules which are situated on the right of the double line in figure \ref{ARkQ}. 
 
 As seen before, the group $N$ is of type $C_3$. It can be realized as a subgroup of $\tilde N$, which is seen as a subgroup of the subgroup of $\GL_6(\C)$ consisting of the upper unitriangular matrices. More precisely, $N$ is the subgroup of $\GL_6(\C)$ generated by the one-parameter subgroups
  $$\left(\begin{matrix} 
           1 & t & 0 & 0 & 0 & 0 \\
	   0 & 1 & 0 & 0 & 0 & 0 \\
	   0 & 0 & 1 & 0 & 0 & 0 \\
	   0 & 0 & 0 & 1 & 0 & 0 \\
	   0 & 0 & 0 & 0 & 1 & t \\
	   0 & 0 & 0 & 0 & 0 & 1 \\
          \end{matrix}\right) \quad 
    \left(\begin{matrix} 
           1 & 0 & 0 & 0 & 0 & 0 \\
	   0 & 1 & t & 0 & 0 & 0 \\
	   0 & 0 & 1 & 0 & 0 & 0 \\
	   0 & 0 & 0 & 1 & t & 0 \\
	   0 & 0 & 0 & 0 & 1 & 0 \\
	   0 & 0 & 0 & 0 & 0 & 1 \\
          \end{matrix}\right) \quad
    \left(\begin{matrix} 
           1 & 0 & 0 & 0 & 0 & 0 \\
	   0 & 1 & 0 & 0 & 0 & 0 \\
	   0 & 0 & 1 & t & 0 & 0 \\
	   0 & 0 & 0 & 1 & 0 & 0 \\
	   0 & 0 & 0 & 0 & 1 & 0 \\
	   0 & 0 & 0 & 0 & 0 & 1 \\
          \end{matrix}\right).
	  $$  
  As a consequence, $\C[N]$ is a quotient of $\C[\tilde N]$. In the following table, to avoid cumbersome indices, we denote minors of a matrix $x \in N$ by indicating with solid dots the entries of the corresponding submatrix of $x$.
 
 Here is the list of the ${\Gamma}$-orbits of the isomorphism classes of indecomposable direct summands $X$ of $T_{M_0}^\vee$ and the realization of $\psi_X$ (see notation \ref{defpsi}) as a minor:
 \begin{center}
  \begin{longtable}{|cc|c|}
    \hline
    \multicolumn{2}{|c|}{${\Gamma}$-orbit} & $\psi_X$ \\
    \hline
   \endfirsthead
    \hline
    \multicolumn{2}{|c|}{${\Gamma}$-orbit} & $\psi_X$ \\
    \hline
   \endhead
    \QBP & \QNP & $\dets{2}{5}{3}{6} = \dets{1}{1}{5}{5}$ \\
    \hline
    \QC & \QO & $\dets{3}{5}{4}{6} = \dets{1}{1}{4}{4}$ \\
    \hline
    \QEP & \QKP & $\dets{2}{4}{4}{6} = \dets{1}{2}{4}{5}$\\
    \hline
    \QF & \QL & $\dets{3}{4}{4}{5} = \dets{2}{2}{4}{4}$\\
    \hline
    \multicolumn{2}{|c|}{\QGP} & $\dets{1}{3}{4}{6}$\\
    \hline
    \multicolumn{2}{|c|}{\QHP} & $\dets{2}{3}{4}{5}$ \\
    \hline
    \multicolumn{2}{|c|}{\QI} & $\dets{3}{3}{4}{4}$ \\
    \hline
    \hline
    \QAP & \QMP & $\dets{1}{5}{2}{6} = \dets{1}{1}{6}{6}$ \\
    \hline
    \QDP & \QJP & $\dets{1}{4}{3}{6} = \dets{1}{2}{5}{6}$ \\
    \hline
  \end{longtable}
 \end{center}
 The last two orbits are those which do not appear in $T_M^\vee$. A simple computation shows that the element $w$ corresponding to the object $M$ is $w = \sigma_c \sigma_b \sigma_c \sigma_a \sigma_b \sigma_c \sigma_b$. As a consequence, the group $N(w)$ is the subgroup of $N$ consisting of the matrices of the form
 $$\left(\begin{matrix}
          1 & & & * & * & * \\
           & 1 & * & * & * & * \\
           & & 1 & * & * & * \\
           & & & 1 & * &  \\
           & & & & 1 &  \\
           & & & & & 1
         \end{matrix}\right)$$
 and $N^w$ is the subvariety of $N$ consisting of the matrices which satisfy
 $$\dets{1}{1}{6}{6} = \dets{1}{2}{5}{6} = 0$$
 and
 $$\dets{1}{1}{5}{5} \neq 0 \quad \text{and} \quad \dets{1}{2}{4}{5} \neq 0 \quad \text{and} \quad \dets{1}{3}{4}{6} \neq 0.$$
\end{ex}

\section*{Acknowledgments}
The author would like to thank his PhD advisor Bernard Leclerc for his advices and corrections. He would also like to thank Christof Gei\ss, Bernhard Keller, Idun Reiten and Jan Schr\"oer for interesting discussions and comments on the topic.

\bibliographystyle{alphanum}
\bibliography{biblio}

\end{document}